\documentstyle{amsppt}
\magnification=1200
\catcode`\@=11
\redefine\logo@{}
\catcode`\@=13

\define \bn{\Bbb N}
\define \bz{\Bbb Z}
\define \bq{\Bbb Q}
\define \br{\Bbb R}

\define \M{{\Cal M}}
\define\Ha{{\Cal H}}
\define\La{{\Cal L}}

\define\rk{\text{rk}~}


\define\pr{\text{pr}}

\define\Aut{\text{Aut}}

\define\sign{\text{sign\ }}

\TagsOnRight
\NoBlackBoxes

\document

\topmatter
\title
On the classification of hyperbolic root systems of
the rank three. Part I
\endtitle

\author
Viacheslav V. Nikulin \footnote{Supported by
Grant of Russian Fund of Fundamental Research
\hfill\hfill}
\endauthor

\address
Steklov Mathematical Institute,
ul. Gubkina 8, Moscow 117966, GSP-1, Russia
\endaddress
\email
slava\@nikulin.mian.su
\endemail

\abstract
It was recently understood that from the point of view of automorphic
Lorentzian Kac--Moody algebras and some aspects of Mirror Symmetry,
interesting hyperbolic root systems should have restricted arithmetic type
and a generalized lattice Weyl vector. One can consider hyperbolic root
systems with these properties as an appropriate hyperbolic analog of
the classical finite and affine root systems.

This series of papers is devoted to classification of hyperbolic
root systems of restricted arithmetic type and
with a generalized lattice Weyl vector $\rho$, having the rank $3$
(it is the first non-trivial rank).
In the Part I we announce classification of the maximal hyperbolic root
systems of elliptic (i.e. $\rho^2 >0$) and parabolic (i.e. $\rho^2=0$)
type, having the rank $3$. We give sketch of the proof. Details of the
proof and further results (non-maximal cases) will be given in Part II.
Classification for hyperbolic type (i.e. $\rho^2<0$) and applications
will be considered in Part III.
\endabstract

\rightheadtext
{Hyperbolic root systems}
\leftheadtext{V.V. Nikulin}
\endtopmatter

\document

\head
0. Introduction
\endhead

It was recently understood that from the point of
view of the theory of Lorentzian (or hyperbolic) Kac--Moody algebras,
some aspects of Mirror Symmetry and some other theories,
interesting hyperbolic root systems should have some very restricted
properties: they have {\it restricted arithmetic type} and
{\it a generalized lattice Weyl vector}
(see Sect. 1.3).
One can consider hyperbolic root systems with these properties
as an appropriate hyperbolic analog of the classical
finite and affine root systems. These properties:
restricted arithmetic type and existence of a generalized lattice
Weyl vector, are so strong that there is a hope to describe
all hyperbolic root systems with these properties. It was shown
\cite{N4}, \cite{N5}, \cite{N11}, \cite{N13} that for the
fixed rank $\ge 3$ number of such root systems is essentially finite.
For example number of such maximal root systems of a fixed rank
$\ge 3$ is finite. It would be very interesting to classify all
hyperbolic root systems of restricted arithmetic type and with
a generalized lattice Weyl vector. It may be of the same
importance as classification of finite and affine root systems
having many applications in Lie algebras theory and other
mathematical and physical theories.

This series of papers is devoted to classification of hyperbolic
root systems of restricted arithmetic type and with a generalized
lattice Weyl vector, having the rank $3$. For the ranks $1$ and $2$ this
problem is trivial. Thus, the rank $3$ is the first non-trivial
case. On the other hand, the rank $3$ case is probably the most
reach, complicated and important in applications.
For higher rank $\ge 5$ some partial results about
this classification for elliptic type
(see the definition below) have been known:
see \'E.B. Vinberg \cite{V3} for the rank $\ge 30$,
F. Esselmann \cite{E2} for the rank $\ge 20$, C. Walhorn \cite{W} for
the rank $5$, R. Scharlau and C. Walhorn \cite{SW}
about some results for the non-compact case of rank $4$.
The reason was that hyperbolic lattices
of the rank $\ge 5$ represent zero, fundamental domains of
their automorphism groups are not compact. The rank $3$ case is
mainly non-compact.

Below we present results in more details.

Let $S$ be a hyperbolic lattice (i.e. a non-degenerate integral
symmetric bilinear form of signature $(1,n)$). We denote by
$(.\,,\,.)$ the form of the lattice $S$. An element $\alpha \in S$
is called a {\it root} if $\alpha^2=(\alpha,\alpha)<0$ and
$\alpha^2|2(\alpha,S)$.
A root $\alpha\in S$ defines a reflection
$$
s_\alpha:x \mapsto x-(2(x,\alpha)/\alpha^2)\alpha
$$
of the lattice $S$.
A subset $\Delta \subset S$ of roots of $S$ is called a
{\it root system with the root lattice $S$}
if properties a), b) and c) below are valid:

a) $s_\alpha(\Delta)=\Delta$ for any $\alpha \in \Delta$;

b) for any root $\beta \in S$ such that $s_\beta (\Delta)=\Delta$
there exists a root $\beta^\prime \in \Delta$ such that
$\beta^\prime =\lambda \beta$, $\lambda \in \bq$;

c) $\Delta$ generates $S\otimes \bq$.

We denote by $W(\Delta)\subset O(S)$ the group generated by
all reflections in roots of $\Delta$ and by
$$
A(\Delta)=\{\phi \in O^+(S)\ |\ \phi(\Delta)=\Delta\}
$$
the group of symmetries of $\Delta$. Let $\M$ be a fundamental
polyhedron of $W(\Delta)$ in the hyperbolic space defined by $S$
(see Sect. 1) and by
$$
A(\Delta,\,\M)=\{\phi \in A(\Delta)\ |\ \phi(\M)=\M\}
$$
the group of symmetries of $\M$. Then
$A(\Delta)=W(\Delta)\rtimes A(\Delta,\,\M)$. The root system
$\Delta$ has {\it restricted arithmetic type} if $A(\Delta)$ has
finite index in $O(S)$. The root system $\Delta$ has
a {\it generalized lattice Weyl vector} if there exist a non-zero
$\rho \in S$ and a subgroup $A\subset A(\Delta,\,\M)$ of finite index
such that $A(\rho)=\rho$. The $\rho$ is called a {\it generalized
lattice Weyl vector} of $\Delta$ for $\M$.

It is sufficient to describe {\it primitive root systems}
when all elements of $\Delta$ are primitive (non-divisible) elements
of $S$. Any other root system {\it is twisted}
to a primitive one: one can get it by an
appropriate multiplication on restricted
natural numbers of elements of a primitive root system.
Further we suppose that $\Delta$ is primitive. Then
$\Delta$ is equivalent to a reflection subgroup
$W=W(\Delta)\subset O(S)$. Really, $\Delta$ is the set of all
primitive roots $\alpha$ of $S$ such that the reflection
$s_\alpha\in W$. Therefor, by definition, $W$ has
{\it restricted arithmetic type and a generalized lattice
Weyl vector} if $\Delta$ does. It is easy to prove that if
$W\subset W^\prime$ are two reflection subgroups of $O(S)$,
both groups $W$ and $W^\prime$ have restricted arithmetic type and
the group $W$ has a generalized lattice Weyl vector,
then $W^\prime$ also has a generalized
lattice Weyl vector. It follows that the full reflection group
$W(S)$ (generated by reflections in all roots of $S$) has restricted
arithmetic type and a generalized lattice Weyl vector if there exists
a reflection subgroup $W\subset W(S)$ having restricted
arithmetic type and a generalized lattice Weyl vector. A hyperbolic
lattice $S$ having this property is called {\it reflective}.
Depending on the type of the generalized lattice Weyl vector
$\rho$ of $W(S)$: $(\rho,\rho)>0$, $(\rho,\rho)=0$
or $(\rho, \rho)<0$, the reflective lattice $S$ is called reflective
of {\it elliptic, parabolic or hyperbolic type}.
Therefor, any primitive root system $\Delta$ of restricted
arithmetic type and with a generalized lattice Weyl vector is
a subsystem of the primitive root system $\Delta (S)$ of all
primitive roots of the reflective lattice $S$. Root systems
$\Delta(S)$ of reflective lattices $S$ are especially important:
the root system $\Delta (S)$ is the maximal hyperbolic root system of
restricted arithmetic type and with a generalized lattice Weyl vector
having the root lattice $S$.

Obviously, $S$ has
elliptic type if and only if $[O(S):W(S)]<\infty$, and this type of
reflective lattices was considered long time ago in papers of
\'E.B. Vinberg and the author (see references). Reflective lattices
of parabolic and hyperbolic type were recently studied in
\cite{N11} and \cite{N13}. We have the following general result

\proclaim{Theorem 1 (\cite{N4}, \cite{N5}, \cite{N11}, \cite{N13})}
For the fixed rank $\rk S\ge 3$ the set of primitive (see Sect. 1.4)
reflective (of any type) hyperbolic lattices $S$ is finite.
In particular, for the fixed rank $\rk S\ge 3$, the set of primitive
root lattices $S$ of hyperbolic root systems of
restricted arithmetic type and with a generalized
lattice Weyl vector is finite.
\endproclaim

To make a further reduction,
we say that an embedding of lattices
$S_1\subset S$ of the same rank is called {\it equivariant} if
it defines the embedding of automorphism groups
$O(S_1)\subset O(S)$. For an equivariant embedding
$S_1\subset S$, the lattice $S$ is reflective if $S_1$ does.
It is easy to prove (and is well-known)
that any hyperbolic lattice $S_1$ has an equivariant
(and actually canonical) embedding into an elementary
hyperbolic lattice $S$. (See Proposition 2.2.2.)
Here a lattice $L$ is called
{\it elementary} if its discriminant group $A_L=L^\ast/L$ is
elementary: any its $p$-component $A_{L_p}\cong (\bz/p\bz)^{r_p}$
is a $p$-elementary Abelian group. To describe all reflective hyperbolic
lattices, it is sufficient to describe all reflective primitive
elementary hyperbolic lattices $S$ and to describe their primitive
reflective equivariant sublattices  $S_1\subset S$ (of finite index).

In this part I, we give the complete list of all
elliptically or parabolically reflective elementary hyperbolic
lattices of rank $3$. See Sect. 2.3, basic Theorem 2.3.2.1 and
Theorems 2.3.3.1, 2.3.4.1. We give idea of the proof (see the
proof of basic Theorem 2.3.2.1).
We use the same method: existing of narrow
parts of the fundamental polyhedron $\M$, which we used to
prove the finiteness Theorem 1. Moreover, we use the additional
arithmetic argument: number of classes of central symmetries of
the fundamental polyhedron $\M$ (see Lemma 2.3.1.2).
Details of the proof will be given
in Part II. In Part II we will also estimate
indexes $[S:S_1]$ of elliptically or parabolically
reflective primitive equivariant sublattices $S_1\subset S$ of
the primitive elementary hyperbolic lattices $S$ of the rank $3$.
Using this estimate, one can find all elliptically or parabolically
reflective hyperbolic lattices of the rank $3$ if one wants.
These estimates are very important in some applications.
In Part III  we will extend these results to the  hyperbolic type of
rank $3$ and will give some applications of our classification.

We hope that our classification will be important in some theories.
We mention that from the point of view of the theory of
Lorentzian Kac--Moody algebras, the rank $3$ case is
the hyperbolic analog of $sl_2$. Moreover the rank $3$ case is more
visual than the rank $\ge 4$ case. For the rank $3$ case one
has to work with polygons (i.e. fundamental polygons of
reflection groups) which are much simpler and much more visual than
polyhedra of dimension $\ge 3$.

\head
1.Hyperbolic reflection groups and hyperbolic root systems
\endhead

\subhead
1.1. Klein model of hyperbolic space
\endsubhead
Let $\Phi$ be a {\it hyperbolic linear space} (i.e. $\Phi$ is a $\br$-linear
space equipped with a non-degenerate symmetric $\br$-bilinear form of
signature $(1,n)$ where $n+1=\dim \Phi$). For $x,y \in \Phi$ we denote by
$(x,y)$ the value at the pair $(x,y)$ of the bilinear form
corresponding to $\Phi$, we also denote $x^2=(x,x)$.
One can associate to $\Phi$ an open cone
$$
V=V(\Phi)=\{x \in \Phi\ |\ x^2>0\}.
\tag{1.1.1}
$$
The $V$ is the disjoint union of two open convex half-cones. We denote
by $V^+=V^+(\Phi)$ one of them.
The {\it hyperbolic space} corresponding to $\Phi$ is
the set of rays
$$
\La=\La (\Phi)=V^+/\br_{++}=\{\br_{++}x\ |\ x \in V^+\}
\tag{1.1.2}
$$
where $\br_{++}$ denote the set of positive real numbers.
The distance $\rho$ in $\La$ is defined by the formula
$$
\cosh{\rho(\br_{++}x,\br_{++}y)}=(x,y)/\sqrt{x^2y^2}.
\tag{1.1.3}
$$
Then the curvature of $\La$ is equal to $-1$. The dimension
$\dim \La =\dim \Phi -1$.

Let $O(\Phi)$ be the group of automorphisms of the hyperbolic
linear space $\Phi$. Then the group of motions of $\La$ is
naturally identified with the subgroup
$O^+(\Phi)\subset O(\Phi)$ of index $2$ of automorphisms
which fix the half-cone $V^+$.

A {\it half-space} of $\La$ is the set
$$
\Ha_\delta^+=\{\br_{++}x \in \La |\ (x,\delta)\ge 0\},
\tag{1.1.4}
$$
where $\delta \in \Phi$ and $\delta^2<0$. The half-space $\Ha_\delta^+$
is bounded by the {\it hyperplane}
$$
\Ha_\delta=\{\br_{++}x \in \La\ |\ (x,\delta)=0\}.
\tag{1.1.5}
$$
The element $\delta \in \Phi$ is defined by the half-space
$\Ha_\delta^+$ (respectively, the hyperplane $\Ha_\delta$) up
to multiplication by elements of $\br_{++}$ (respectively, by
elements of the set $\br^\ast$ of non-zero real numbers). The
$\delta$ is called the {\it orthogonal vector} to the half-space
$\Ha_\delta^+$ (respectively to the hyperplane $\Ha_\delta$).

Consider two half-spaces $\Ha_{\delta_1}^+$,
$\Ha_{\delta_2}^+$ orthogonal
to $\delta_1,\,\delta_2 \in \Phi$.
Assume that $\Ha_{\delta_1}^+\not\subset \Ha_{\delta_2}^+$,
$\Ha_{\delta_2}^+\not\subset \Ha_{\delta_1}^+$, and
intersection $\Ha_{\delta_1}^+\cap \Ha_{\delta_2}^+$ contains a non-empty
open subset of $\La (\Phi)$. Then we have two cases:

a) $\Ha_{\delta_1}^+\cap \Ha_{\delta_2}^+$ is an angle of a value
$\phi$ where
$$
\cos{\phi}=(\delta_1,\delta_2)/\sqrt{\delta_1^2\delta_2^2}
\tag{1.1.6}
$$
if $-1\le (\delta_1,\delta_2)/\sqrt{\delta_1^2\delta_2^2}\le 1$.

b) Intersection of the hyperplanes $\Ha_{\delta_1}$, $\Ha_{\delta_2}$ is
empty (they are hyperparallel) and the distance $\rho$ between them
is equal to
$$
\cosh \rho =(\delta_1,\delta_2)/\sqrt{\delta_1^2\delta_2^2}
\tag{1.1.7}
$$
if $1\le (\delta_1,\delta_2)/\sqrt{\delta_1^2\delta_2^2}$.

As usual, we complete $\La$ by the set of {\it infinite points}
$\br_{++}c$ where $c\in \Phi$, $c^2=0$ and
$(c,V^+)>0$. Then $\La$ with its infinite points
is the closed ball
$$
\overline{\La}=(\overline{V^+}-\{0\})/\br_{++},
\tag{1.1.8}
$$
where $\overline{V^+}$ is the closure of $V^+$ in
the linear space $\Phi$. The boundary
$\La_\infty=\overline\La - \La$ is the sphere of dimension $n-1$.
It is called $\infty${\it-sphere} or {\it infinity} of $\La$.
We also add infinite points to half-spaces,
hyperplanes and other subsets of $\La$ considering
their closure in $\overline{\La}$.

{\it A non-degenerate convex locally finite polyhedron} $\M$ in a
hyperbolic space $\La$ is the intersection of a set of half-spaces:
$$
\M=\bigcap_{\delta \in P(\M)}{\Ha_\delta^+}.
\tag{1.1.9}
$$
The $\M$ is called {\it non-degenerate} if $\M$ contains a non-empty open
subset of $\La$. The $\M$ is called {\it locally finite} if for any
$X\in \La$, there exists its open neighbourhood $U\subset \La$ such that
$\M\cap U$ is the intersection with $U$ of a finite set of half-spaces
in $\La$. Here $P(\M)\subset \Phi$ is a finite or countable subset of
elements with negative square. We assume that no two of these elements
are proportional and every half-space $\Ha_\delta^+$,
$\delta \in P(\M)$, defines a face of $\M$
of the highest dimension $\dim \La -1$
(it is clear what this means, since $\M$ is locally finite). Then, the
polyhedron $\M$ defines the subset $P(\M)\subset \Phi$ uniquely up to
multiplication by $\br_{++}$ elements of $P(\M)$. The set $P(\M)$ is
called {\it the set of vectors orthogonal to $\M$
(or to faces of highest dimension of $\M$)}.
Below, we only consider locally finite and convex
polyhedra.

A polyhedron $\M$ is called {\it elliptic} (or it has
{\it geometrically finite volume}) if $\M$ is a convex
envelope of a finite set of points in $\La$, some of them
at infinity of $\La$.

\subhead
1.2. Reflection groups of hyperbolic lattices, and hyperbolic root systems
\endsubhead
Let $S$ be a hyperbolic (that is, of signature $(1,n)$) integral
symmetric bilinear form over $\bz$. This means that $S$ is a
free $\bz$-module of finite rank equipped with an integral
valued symmetric $\bz$-bilinear form and the corresponding real
form $\Phi=S\otimes \br$ is hyperbolic (it has
signature $(1,n)$ like in Sect. 1.1). For brevity,
$S$ will be referred to as the {\it hyperbolic lattice}. Like
in Sect. 1.1, we correspond to a hyperbolic lattice $S$ the
hyperbolic space $\La=\La (S)=\La(S\otimes \br)$. From the
theory of integral quadratic forms (more generally, from
the theory of arithmetic discrete groups), it is known that the
automorphism group $O^+(S)$ is discrete in $\La$; if
$\rk S\ge 3$, the group $O^+(S)$ has a fundamental domain of
finite volume which is
an elliptic convex polyhedron
(e. g. see \cite{C}, \cite{R}, \cite{VSh}).

An element $\alpha \in S$ is called {\it root} if
$\alpha^2<0$ and $\alpha^2 | 2(S,\alpha)$. If $\alpha$ is
a root, then $-\alpha$ is also a root.
A non-zero element $x \in S$
is called {\it primitive} if $x/m \notin S$ for any natural $m>1$.
Obviously, any root $\alpha$ is multiple to a primitive root which
is called {\it the primitive root of $\alpha$.} If
$\alpha_0$ is a primitive root, then $\lambda\alpha_0$ is a root only
for a finite number of $\lambda \in \bn$. Really,
${\alpha_0\over (\lambda/2)} \in S^\ast=Hom(S, \bz)$, and
$\lambda\le 2a(S)$ where $a(S)$ is the exponent (the maximal order of
elements) of the discriminant group $A_S=S^\ast /S$ of $S$. In
particular, $\lambda \le 2|\det(S)|$ where the {\it determinant of
$S$} is equal to $\det(S)=\det((e_i,e_j))$ for a bases
$e_1,\dots, e_n$ of $S$. Thus, to find roots, it is sufficient
to find primitive roots and then multiply them on restricted
multiplicities.

Any root $\alpha \in S$
defines the {\it reflection $s_\alpha \in O^+(S)$}
which is given by the formula
$$
x\mapsto x-{2(x,\alpha)\over \alpha^2}\alpha,\ \ x \in S.
\tag{1.2.1}
$$
The reflection depends only from $\pm\alpha_0$ where
$\alpha_0$ is the primitive root of $\alpha$.
The reflection $s_\alpha$ acts in $\La=\La (S)$ as the
reflection in the hyperplane $\Ha_\alpha$ orthogonal to the root
$\alpha$ (i.e. $s_\alpha$ is identical on $\Ha_\alpha$ and
changes two half-spaces bounded by $\Ha_\alpha$).
It is easy to see that an automorphism
$\phi\in O^+(S)$ gives a reflection in a hyperplane of
$\La$ if and only if $\phi$ is
a reflection in some root of the lattice $S$.

All reflections $s_\alpha$ in all roots of $S$ generate
the {\it reflection group} $W(S)\subset O^+(S)$ of the lattice $S$.
We shall consider subgroups $W\subset W(S)$ generated by some set of
reflections of $S$. We shall call $W$ by {\it reflection subgroup of
$S$}.

Let $W$ be a reflection subgroup of $S$.
From the theory of groups generated by reflections (e.g. see
\cite{V1}, \cite{V4}, \cite{VSh}), it follows that $W$ has a fundamental
polyhedron $\M$ in $\La$ which is the closure of a
connected component of $\La$ without hyperplanes of all
reflections from $W$. The polyhedron $\M$
is a convex locally finite polyhedron. We always assume that the
set $P(\M)$ of orthogonal vectors to $\M$ is a set of roots of $S$.
It is also called the {\it set of simple roots}.
The group $W$ is generated by reflections $s_\alpha$, $\alpha \in P(\M)$.
The set $\{s_\alpha\ |\  \alpha \in W(P(\M))\}$ gives then all
reflections of the reflection subgroup $W$. The set
$\Delta=W(P(\M))$ is called a {\it (hyperbolic) root system} of the lattice
$S$. One can then look at the lattice $S$ as a {\it root lattice}.
Hyperbolic root systems are very important for some theories,
for example, for the theory of Lorentzian (or hyperbolic) Kac--Moody
algebras (see \cite{Kac}).
It is clear that to find all possible root systems, it is
sufficient to restrict by {\it primitive root systems}
containing only primitive roots of $S$. All other
root systems are {\it twisted} to primitive root systems: there roots
are multiple (by bounded multiplicities)
to roots of a primitive root system. To find a root
system $\Delta$ it is sufficient to find a system of its simple roots
$P(\M)$ since $\Delta=W(P(\M))$ where $W$ is generated by reflections in
$P(\M)$. To find a system of simple roots $P(\M)$, it is sufficient to
find a system of simple primitive roots $P(\M)_{\pr}$ which is equivalent
to the reflection subgroup $W\subset W(S)$. Arbitrary system of
simple roots $P(\M)$ is twisted to $P(\M)_{\pr}$: one can get it
multiplying elements of $P(\M)_{\pr}$ by some bounded natural numbers
({\it twisting coefficients}).

\subhead
1.3. Restricted arithmetic type and a generalized lattice Weyl vector
\endsubhead
We denote by
$$
A(P(\M))=\{g \in O^+(S)\ |\ g(P(\M))=P(\M)\}
\tag{1.3.1}
$$
the {\it group of symmetries of the fundamental polyhedron (together with
the set $P(\M)$ of its orthogonal roots).}
The group $W$ with the set $P(\M)$ of simple roots (equivalently,
the root system $\Delta=W(P(\M))$) have
{\it restricted arithmetic type} (see \cite{N11}, \cite{N10}),
if the semi-direct product
$\Aut(\Delta)=
W\rtimes A(P(\M))$ has finite index in $O^+(S)$. Obviously,
$A(P(\M))\subset A(P(\M)_{\pr})$, and $P(\M)_{\pr}$ has restricted
arithmetic type if $P(\M)$ does. If $W$ and $P(\M)_{\pr}$ have
restricted arithmetic type, we say that the reflection subgroup
$W$ {\it has restricted arithmetic type}.
A non-zero element $\rho \in S$ is called a
{\it generalized lattice Weyl vector} if $A(\rho)=\rho$ for
a subgroup $A\subset A(P(\M))$ of finite index (equivalently,
the orbit $A(P(\M))(\rho)$ is finite). Reflection
subgroups $W$ with $P(\M)$ (and corresponding root systems
$\Delta =W(P(\M))$) of restricted arithmetic type and with
a generalized lattice Weyl vector are important for the
theory of K3 and other surfaces (see \cite{P-\u S\u S},
\cite{N1}---\cite{N9}, \cite{N12}, \cite{N13}, \cite{AN1}, \cite{AN2}),
so called Borcherds type automorphic products and
theory of automorphic Lorentzian Kac--Moody algebras
(see \cite{B2}---\cite{B6}, \cite{GN1}---\cite{GN7},
\cite{N10}, \cite{N11}) and other theories
(e. g. see \cite{HM1}, \cite{HM2}, \cite{M},
\cite{Kaw1}, \cite{Kaw2}, \cite{CCL} on applications in Physics).
It seems, hyperbolic root systems satisfying these
conditions are of the same importance as well-known finite
and affine root systems.

It is an interesting but
difficult problem to find all hyperbolic root systems $\Delta$
(or corresponding systems of simple roots $P(\M)$) of restricted
arithmetic type and with a generalized lattice Weyl vector.

\smallpagebreak

Below we shall make some general remarks about this problem.

First, root systems, simple root systems
and reflection subgroups $W$ of restricted arithmetic type
with a generalized lattice Weyl vector are divided on three types:

{\it Elliptic type:} There exists a generalized lattice
Weyl vector $\rho$ with positive square $\rho^2>0$.

{\it Parabolic type:} There exists a generalized lattice
Weyl vector $\rho$ with zero square $\rho^2=0$ and there
does not exist a generalized lattice Weyl vector with
positive square.

{\it Hyperbolic type:} There exists a generalized lattice
Weyl vector $\rho$ with negative square $\rho^2<0$ and
there does not exist a generalized lattice Weyl vector with
positive or zero square.

Obviously, $W$ and $P(\M)$ of restricted arithmetic type have
a generalized lattice Weyl vector $\rho$ if and only if
$W$ and $P(\M)_{\pr}$ have restricted arithmetic type and
the same generalized lattice Weyl vector $\rho$. Thus, to find
all $W$ and $P(\M)$ of restricted arithmetic type and
with a generalized lattice Weyl vector, we can restrict by
primitive sets $P(\M)_{\pr}$ which are equivalent to
reflection subgroups $W\subset W(S)$.
One should find all reflection subgroups $W\subset W(S)$
of restricted arithmetic type and with a generalized lattice
Weyl vector for $P(\M)_{\pr}$ where $\M$ is a fundamental polyhedron
for $W$. Here the set $P(\M)_{\pr}$ is defined by the reflection
subgroup $W\subset W(S)$
(up to choosing the fundamental polyhedron $\M$).

\subhead
1.4. Reflective hyperbolic lattices
\endsubhead
Assume that a reflection subgroup
$W\subset W(S)$ has restricted arithmetic type and a generalized
lattice Weyl vector $\rho$ (for $P(\M)_{\pr}$ where $\M$ is
a fundamental polyhedron of $W$). Suppose that
$W\subset W^\prime \subset W(S)$ where
$W^\prime$ is a reflection subgroup of restricted arithmetic type.
We can choose a fundamental polyhedron $\M^\prime$ for $W^\prime$
such that $\M^\prime\subset \M$. If $\rho$ is a generalized
lattice Weyl vector for $P(\M)_{\pr}$, then $\rho$ is a
generalized lattice Weyl vector for $P(\M^\prime)_{\pr}$
(one should use that $A(P(\M^\prime )_{\pr})\cap A(P(\M)_{\pr})$ has
finite index in $A(P(\M^\prime)_{\pr})$ since $W$ has
restricted arithmetic type). Thus,
$W^\prime$ also has a generalized
lattice Weyl vector. If $W$ has elliptic type, then
$W^\prime$ also has elliptic type. If $W$ has parabolic type,
then $W^\prime$ has parabolic or elliptic type. If $W$ has
hyperbolic type, then $W^\prime$ has either hyperbolic,
or parabolic, or elliptic type.
The full reflection group $W(S)$ of the lattice $S$ obviously
has restricted arithmetic type: it is normal in $O^+(S)$.
Thus, $W(S)$ itself has restricted arithmetic type and
a generalized lattice Weyl vector if $W(S)$
has a reflection subgroup $W\subset W(S)$ of restricted
arithmetic type and with a generalized lattice Weyl vector.
Hyperbolic lattices $S$ with this property are called
{\it reflective}. Thus, (equivalently) a hyperbolic lattice $S$ is
{\it reflective} if it has a generalized lattice Weyl vector for
its reflection group $W(S)$. Depending on the type
(elliptic, parabolic or hyperbolic) of $W(S)$
the reflective lattice $S$ is called {\it reflective of
elliptic, parabolic or hyperbolic type} or {\it
elliptically, parabolically or hyperbolically reflective}.
It seems that these maximal reflection groups $W(S)$ having a
generalized lattice Weyl vector are especially interesting.
Moreover, for reflective lattices $S$, one can try to find out
all reflection subgroups $W\subset W(S)$ of restricted arithmetic
type and with a generalized lattice Weyl vector.
If the lattice $S$ is not reflective, any
reflection subgroup $W\subset W(S)$ of restricted arithmetic type
does not have a generalizes lattice Weyl vector.

Multiplication of the form of a lattice $S$ does not change
anything. Thus, we can always suppose that the lattice $S$
is {\it primitive}, i.e. $S(1/k)$ is not a lattice for any
natural $k>1$. Here $K(r)$ denote a lattice which is obtained by
multiplication of the form of a lattice $K$ on $r\in \bq$.

We have the following general result:

\proclaim{Theorem 1.4.1 (\cite{N4}, \cite{N5}, \cite{N11},
\cite{N13})} For the fixed rank $\rk S\ge 3$ the set of
primitive reflective hyperbolic lattices $S$ is finite.
\endproclaim

Cases $\rk S=1$ or $2$ are trivial. Any hyperbolic lattice $S$
of $\rk S=1$ is reflective since $O^+(S)$ is trivial. If $\rk S=2$,
the lattice $S$ is reflective if and only if either $S$ has a non-zero
element with square $0$ (then $O^+(S)$ is finite and $S$ is
reflective of elliptic type) or $S$ has at least
one root (then $S$ is also reflective of elliptic type).

\smallpagebreak

The main subject of these series of papers is to describe hyperbolic
reflective lattices of rank $3$. We apply the same method (of
narrow parts of polyhedra) which was applied to prove Theorem 1.
But to find all possibilities, we need to improve this method
and add some additional arithmetic arguments.
Certainly, classification is more delicate problem than
proving finiteness.

\head
2. Classification of reflective hyperbolic lattices of the rank $3$
and of elliptic or parabolic type: formulations
\endhead

\subhead
2.1. The principle of classification
\endsubhead
Let $S_1$ be a hyperbolic lattice and
$S_1\subset S_2$ its overlattice. Here and in what follows we
consider only embeddings of lattices of the same rank. If
the embedding $S_1\subset S_2$ induces an embedding of
the automorphism groups $O(S_1)\subset O(S_2)$ (in general,
these groups are commensurable), we say that the embedding
$S_1\subset S_2$ is an {\it equivariant embedding.} We remark that
index $[O(S_2):O(S_1)]$ is always finite because both groups
$O^+(S_2)$ and $O^+(S_1)$ have a fundamental domain of finite volume
in naturally identified hyperbolic spaces $\La (S_1\otimes \br)=
\La (S_2\otimes \br)$. It follows that for the equivariant embedding,
$W(S_1)\subset W(S_2)$ and both groups $W(S_1)$ and $W(S_2)$
have restricted arithmetic type. It follows that the lattice
$S_2$ is reflective if the lattice $S_1$ is reflective. Obviously,
$S_2$ is primitive if $S_1$ does. A lattice $S$ is called
{\it maximal} if it does not have equivariant embeddings to other
lattices of the same rank except the lattice $S$ itself.

Our principle of classification of reflective hyperbolic lattices of
the rank $3$ is as follows: We define {\it elementary lattices}
such that any lattice has an equivariant embedding into an elementary
one. We give the list of
{\it reflective elementary hyperbolic lattices of the rank $3$.}
Thus, any reflective hyperbolic lattice of the rank $3$ has an equivariant
embedding into a lattice of this list.
For any reflective elementary hyperbolic lattice $S$ of
rank $3$ we calculate the set
of simple roots $P(\M)_{\pr}$ of $W(S)$ and the Gram matrix
$G(P(\M)_{\pr})=((\alpha,\beta))$, $\alpha,\,\beta \in P(\M)_{\pr}$.
Moreover, we estimate indexes $[S:S_1]$ of its primitive
reflective sublattices $S_1\subset S$, and one can find all of
them if one needs (their number is obviously finite).
It seems, reflective elementary hyperbolic lattices are especially
interesting (e. g. they contain maximal reflective hyperbolic
lattices), and finding all of them is especially important.
On the other hand, the mere fact that a
 lattice $S_1$ has an equivariant embedding $S_1\subset S$ into a
known lattice $S$ contains a lot of information about the lattice $S_1$
itself. For example, $\det(S_1)=det(S)[S:S_1]^2$,
a primitive root $\alpha_1$ of $S_1$ is equal to $t\alpha$ where
$\alpha$ is a primitive root of $S$ and $t\le [S:S_1]$, it
follows $\alpha_1^2=t^2\alpha^2$.

\subhead
2.2. Elementary hyperbolic lattices of the rank 3
\endsubhead
We recall some elements of the {\it discriminant form technique}
(see \cite{N1}).

Let $L$ be an arbitrary lattice (equivalently, a non-degenerate
integral symmetric bilinear form). The lattice $L$ is called
{\it  odd} if there exists $x \in L$ with odd square $x^2$.
Otherwise, $L$ is called {\it even}.

The group $A_L=L^\ast/L$ is called the
{\it discriminant group of $L$}. Its order $\sharp A_L=|\det (L)|$.
We extend the form of $L$ to the dual lattice $L\subset L^\ast =
Hom(L,\bz)$. It then takes values in $\bq$ and defines
(considering $\mod L$) a non-degenerated finite symmetric bilinear form
$$
b_L:A_L\times A_L \to \bq/\bz.
\tag{2.2.1}
$$
If the lattice $L$ is even, one
similarly defines a finite quadratic form
$$
q_L:A_L\to \bq/2\bz
\tag{2.2.2}
$$
with the symmetric bilinear form $b_L$. These forms are called
the {\it discriminant bilinear and the discriminant quadratic form} of
the lattice $L$.

An overlattice $L\subset \widetilde{L}$ is defined by
an isotropic (for $b_L$)
subgroup $H_{\widetilde{L}}=\widetilde{L}/L\subset A_L$.
The discriminant bilinear form of $\widetilde{L}$ is then equal to
$$
b_{\widetilde{L}}=b_L|(H_{\widetilde{L}}^\perp/H_{\widetilde{L}}).
\tag{2.2.3}
$$
If the lattice $L$ is even and
the subgroup $H_{\widetilde{L}}$ is also isotropic for the
discriminant quadratic form $q_L$, then $\widetilde{L}$ is even
and
$$
q_{\widetilde{L}}=q_L|(H_{\widetilde{L}}^\perp/H_{\widetilde{L}}).
\tag{2.2.4}
$$
For a prime $p$, we denote as $A_p$ the $p$-component of
$A_L$,
and $b_p$ and $q_p$ restrictions of $b_L$ and $q_L$ on
the subgroup.
Consider the $p$-adic lattice $L_p=L\otimes \bz_p$ where
$\bz_p$ is the ring of $p$-adic integers and $\bq_p$ the field of
$p$-adic numbers. One can similarly define the discriminant group and
the discriminant forms
$A_{L_p}$, $b_{L_p}$, $q_{L_p}$ of $L_p$.
We have natural identifications:
$A_p=A_{L_p}$, $b_p=b_{L_p}$ and $q_p=q_{L_p}$.

\definition{Definition 2.2.1} A lattice $L$ is called {\it elementary}
if for any prime $p\mid det(L)$ the $p$-component $A_{L_p}$
is a $p$-elementary Abelian group:
$A_{L_p}\cong (\bz/p\bz)^{r_p}$. Here $\{r_p=r_p(L)\ |\  p\mid \det(L)\}$
are invariants of the elementary lattice $L$.
We have $|\det (L)|=\prod_{p}{p^{r_p}}$. In particular, the invariants
$r_p(L)$ are defined by $\det(L)$.
\enddefinition

We have the following simple and well-known

\proclaim{Proposition 2.2.2} Any lattice $L$ has an equivariant
embedding into an elementary lattice.
\endproclaim

\demo{Proof} Let $p^{t_p}$ be the exponent of the
group $A_p=A_{L_p}$. It means that
$p^{t_p}A_p=0$ and $p^{t_p-1}A_p\not=0$.
If $L$ is elementary, $L$ has an equivariant embedding into itself.
If $L$ is not elementary, one of exponents $t_p>1$.
Obviously, $b_{L_p}:A_p\times A_p\to {1\over p^{t_p}}\bz/\bz$.
It follows that subgroup $H_p=p^{t_p-1}A_p\subset A_p$ is isotropic.
It is not zero and is invariant with respect to $O(L)$. It follows that
$H_p$ defines an equivariant embedding $L\subset \widetilde{L}$
with $\widetilde{L}/L=H_p$ and index
$\sharp H_p$. If the lattice $\widetilde{L}$ is not elementary, we
repeat the procedure for $\widetilde{L}$. This procedure is
finite because $\det(L_1)=\det(L)[L_1:L]^{-2}$ for
an overlattice $L\subset L_1$ of finite index $[L_1:L]$
where $\det(L)$ and $\det(L_1)$ are integers. It proves Proposition.
\enddemo

Let $L$ be an elementary lattice, $m \in \bn$ and $m$ is
square-free. The lattice
$$
L^{\ast,m}(m)=(L^\ast \cap {{1\over m}L})(m)
\tag{2.2.5}
$$
is called {\it $m$-dual to the lattice $L$}. This lattice is also
elementary and
$$
r_p(L^{\ast,m}(m))=
\cases
r_p(L), &\text{if $p\nmid m$}\\
\rk L-r_p(L) &\text{if $p\mid m$}
\endcases.
\tag{2.2.6}
$$
One can also see
that $L^{\ast,m}(m)\otimes \bz_p=L\otimes \bz_p$ if $p\nmid m$, and
$L^{\ast,m}(m)\otimes \bz_p=(L\otimes \bz_p)^\ast (p)$ if $p\mid m$.
Moreover,
$$
(L^{\ast,m}(m))^{\ast,m}(m)=L
\tag{2.2.7}
$$
which shows that it really
is a duality. More generally, for two square-free numbers $m_1$, $m_2$
one has
$$
(L^{\ast,m_1}(m_1))^{\ast,m_2}(m_2)=L^{\ast,m}(m)
\tag{2.2.7'}
$$
where $m=\text{l.c.m.}(m_1,m_2)/\text{g.c.d}(m_1,m_2)$.
It follows that one has a natural identification
$$
O(L)=O(L^{\ast,m}(m)).
\tag{2.2.8}
$$
If $L$ is a hyperbolic lattice, one can look at this equality as
follows: {\it Hyperbolic spaces and actions of automorphism groups
in these spaces of $L$ and $L^{\ast,m}(m)$ are naturally
identified.}

Now we consider only primitive elementary hyperbolic lattices
of rank $3$. Let $F$ be a primitive elementary hyperbolic lattice of
the rank $3$ and $\det(F)=\prod_p{p^{r_p}}$.
Since $F$ is primitive, $0\le r_p\le 2$. Let
$m$ be the product of all prime $p$ such that $p\mid \det(F)$ and
$r_p=2$. Equivalently, $\det(F)=(\det(F)/m)m$ where both natural numbers
$\det(F)/m$ and $m$ are square-free. By \thetag{2.2.6},
the $m$-dual lattice
$$
S=F^{\ast,m}(m)
\tag{2.2.9}
$$
has then a square-free determinant $\det(S)=\det(F)/m$.
Thus we have

\proclaim{Proposition 2.2.3} Any primitive elementary hyperbolic
lattice $F$ of the rank three is $m$-dual to a
hyperbolic lattice $S$ with a square-free determinant
$d=\det (S)$ where $m \mid d$. Then $\det (F)= dm$. A lattice
with a square-free determinant is (obviously) elementary and primitive.
\endproclaim

By Proposition 2.2.3, it is sufficient to describe
hyperbolic lattices of the rank three with a square-free determinant.
Let $S$ be a hyperbolic lattice with a square-free
determinant $d=\det(S)$. If $d$ is odd, then $S$ is also odd
(otherwise, the $2$-adic lattice $S\otimes \bz_2$ is even unimodular of
odd rank, which is impossible).
If $d$ is even, then the lattice $S$ may be even or odd.
If $S$ is odd (respectively even) we say that $S$ has {\it odd type}
(respectively {\it even type}).
For odd prime $p\mid d$ the discriminant
bilinear form $b_{S_p}\cong b_{\theta_p}(p)$ where
$b_{\theta_p}(p)$ is the discriminant bilinear form of the
$p$-adic $1$-dimensional lattice $\langle \theta_p p\rangle$ where
$\theta_p\in \bz_p^\ast/(\bz_p^\ast)^2\ \cong \{\pm 1\}$ by
isomorphism $\theta_p \to \left({\theta_p\over p}\right)=(-1)^{\eta_p}$
where $\left({.\over p}\right)$ is the Legendre symbol and
$\eta_p\in \{0,1\}$. (Here and in what follows $\langle A \rangle$
denote a lattice with the matrix $A$.) Thus, the
lattice $S$ with a square-free determinant $d$
has the invariant
$$
\eta=\{\eta_p \
| \ \text{odd\ } p\mid d\}\ \text{where}\  \eta_p \in \{0,\,1\}\
\text{and}\ (-1)^{\eta_p}=\left({\theta_p\over p}\right).
\tag{2.2.10}
$$
Below we code this invariant by a non-negative integer
$\eta$ having a binary decomposition
$$
\eta=\eta_{p_t}....\eta_{p_1}
\tag{2.2.10'}
$$
where $p_1,\dots,p_t$ are all odd prime divisors of $d$ in
increasing order. For example, for $d=30=2\cdot 3\cdot 5$ the invariant
$\eta=2$ means that $\left({\theta_3\over 3}\right)=1$,
$\left({\theta_5\over 5}\right)=-1$

We have

\proclaim{Proposition 2.2.4} A hyperbolic lattice $S$ of
rank three and with a square-free determinant is defined up
to isomorphism by the invariants $(d,\,type,\,\eta)$.
\endproclaim

\demo{Proof} Invariants $(d,\,type,\,\eta)$ define
the discriminant bilinear form of $S$. By \cite{N1},
the signature, type and discriminant bilinear form of $S$
define the genus of $S$. By general results of
Eichler and Kneser on spinor genus of indefinite
lattices of rank $\ge 3$ (see \cite{E1}, \cite{Kn} and also
\cite{C}), elementary lattices of rank $\ge 3$ have only one class
in a genus. It proves the statement.
\enddemo

\definition{Definition  2.2.5} A hyperbolic lattice $S$ of
the rank three and with a square-free determinant $d$ is
called {\it main} if $type \equiv d\mod 2$.
In other words, the lattice $S$ should be even if the determinant
$d$ is even. If the determinant $d$ is odd, then the lattice
$S$ will be necessarily odd. In particular, main hyperbolic
lattices of rank three and with a square-free determinant
are defined by the invariants $(d,\eta)$.
\enddefinition

Main hyperbolic lattices of the rank three and with a square-free
determinant are the most important. Really, suppose that $\widetilde{S}$
has rank three and square-free determinant, but it is not main.
Then the determinant
$d(\widetilde{S})=2d$ is even but the lattice
$\widetilde{S}$ is odd. Let
$\widetilde{S}^{ev}\subset \widetilde{S}$
be the maximal even sublattice of $\widetilde{S}$. It has index two.
Considering $\widetilde{S}\otimes \bz_2$, one finds that
$\widetilde{S}^{ev}=S(2)$ where
$S$ has the determinant $d$ and is then main odd
hyperbolic lattice. By the construction,
$O(\widetilde{S})\subset O(S)$.
Thus, {\it main hyperbolic lattices of the rank three and with
a square-free determinant have maximal automorphism groups}.
In particular, if $\widetilde{S}$ is reflective,
then $S$ is also reflective. Moreover, we get that
any non-main hyperbolic lattice $\widetilde{S}$
is the index two odd overlattice
$$
S(2)\subset \widetilde{S}
\tag{2.2.11}
$$
where $S$ is a main odd hyperbolic lattice. If $S$ has invariants
$(d,\text{odd},\eta)$ (where $d$ is square-free and odd),
then $\widetilde{S}$ has invariants
$(2d,\text{odd},\eta+\omega(d))$ where
$\omega(d)=\{\omega(d)_p=(p^2-1)/8\mod 2\ \ |\ \
\text{odd prime\ } p\mid d\}$.
We use that $\left({2\over p}\right)=(-1)^{(p^2-1)/8}$.

For the fixed $S$, all index $2$ odd lattices $\widetilde{S}$
in \thetag{2.2.11} are isomorphic by Proposition 2.2.3.
Thus, \thetag{2.2.11} defines the one-to-one correspondence
$S \leftrightarrow \widetilde{S}$ between odd hyperbolic
lattices of rank three and with a square-free determinant
where $S$ has odd determinant $d$ and $\widetilde{S}$ has
the even determinant $2d$.

We have mentioned that \thetag{2.2.11} defines embedding
$$
O(\widetilde{S})\subset O(S)
$$
since the lattice $S(2)$ is the
maximal even sublattice of $\widetilde{S}$. Further we can consider
two cases.

{\it Case when lattices $S$ and $\widetilde{S}$ are equivariantly
equivalent:} By definition, it means that the odd overlattice
$S(2)\subset \widetilde{S}$ is unique. Then \thetag{2.2.11} gives
$O(\widetilde{S})=O(S)$, and lattices $S$,
$\widetilde{S}$ have the same reflective type.
By the discriminant form technique,
this means that the $2$-component
$A_{S(2)_2}\cong (\bz/2\bz)^3$ of the discriminant group of
$S(2)$ has an unique element $w\in A_{S(2)_2}$ with
$q_{S(2)}(w)=1\mod 2$. Then $\widetilde{S}/S(2)=(\bz/2\bz)w$,
and index two odd overlattice in \thetag{2.2.11} is unique.
The discriminant quadratic form $q_{S(2)_2}$ is defined by its
signature $\sign q_{S(2)_2}\mod 8$ (see \cite{N1} where all
finite quadratic forms on $2$-elementary groups were classified).
One can easily check that the element $w$ is unique if and only if
$\sign q_{S(2)_2}\equiv \pm 1 \mod 8$. We have (e.g. see \cite{N1})
$$
\sum_p{\sign q_{S(2)_p}}\equiv \sum_{p\mid d}
{(1-p+4\eta_p+4\omega(d)_p)}+\sign q_{S(2)_2}\equiv 1-2
\equiv -1\mod 8.
$$
It follows that $S$ and $\widetilde{S}$ are equivariantly
equivalent if and only if
$$
\sum_{p\mid d}
{(1-p+4\eta_p+4\omega(d)_p)}\equiv \pm 1-1\equiv 0\ \text{or}\  6 \mod 8.
$$

{\it Case when $S$ and $\widetilde{S}$ are not equivariantly
equivalent:} by definition it means that the
odd index two overlattice $\widetilde{S}$ in \thetag{2.2.11} is
not unique. By consideration above, it is equivalent to
$$
\sum_{p\mid d}
{(1-p+4\eta_p+4\omega(d)_p)}\mod 8 \notin \{0~\mod 8,\  6~\mod 8\}.
$$
One can check that then there are three different
possibilities for the element $w\in A_{S(2)_2}$ with
$q_{S(2)}(w)=1\mod 2$. Thus, there are three different
overlattices $\widetilde{S}$ in \thetag{2.2.11}. All of them
give isomorphic lattices $\widetilde{S}$ but one cannon expect that
automorphism groups $S$ and $\widetilde{S}$ are the same (and it
is never the case). The \thetag{2.1.11} defines only embedding
$O(\widetilde{S})\subset O(S)$, and we only have that {\it the reflective
type of $\widetilde{S}$ is dominated by the reflective type of $S$.}
To find all reflective lattices $\widetilde{S}$, one should first find
all reflective lattices $S$, and then check if the lattice
$\widetilde{S}$ is also reflective. We summarize these considerations
below.

\proclaim{Proposition 2.2.6} All non-main hyperbolic
lattices $\widetilde{S}$
of rank three and with a square-free determinant
(remind that non-main means that $\widetilde{S}$
is odd but $\det(\widetilde{S})=2d$ is even)
are in one-to-one correspondence $S\leftrightarrow \widetilde{S}$ with
main odd hyperbolic lattices $S$ of rank three and with a square-free
determinant $d=\det(S)=\det(\widetilde{S})/2$. The correspondence is
defined by the embedding of lattices
$$
S(2)\subset \widetilde{S},
\tag{2.2.12}
$$
where $S(2)$ is the maximal even sublattice of $\widetilde{S}$
(it has index two).

If $S$ has invariants $(d,\,odd,\,\eta)$, where $d$ is odd
square-free, then
$\widetilde{S}$ has invariants $(2d,\,odd,$ $\eta+\omega(d))$
where
$$
\omega(d)=\{\omega(d)_p=(p^2-1)/8\mod 2\ \ |\ \
\text{odd prime\ } p\mid d\}.
\tag{2.2.13}
$$

If
$$
\sum_{p\mid d}
{(1-p+4\eta_p+4\omega(d)_p)}\equiv \pm 1-1\equiv 0\ \text{or}\  6 \mod 8,
\tag{2.2.14}
$$
the overlattice $\widetilde{S}$ in \thetag{2.2.12} is unique,
\thetag{2.2.12} defines isomorphism $O(\widetilde{S})=O(S)$ and
lattices $S$ and $\widetilde{S}$ have the same reflective type.
We then say that lattices $S$ and $\widetilde{S}$ are
equivariantly equivalent.

If
$$
\sum_{p\mid d}
{(1-p+4\eta_p+4\omega(d)_p)}\mod 8 \notin \{0~\mod 8,\  6~\mod 8\},
\tag{2.2.15}
$$
\thetag{2.2.12} only defines the embedding $O(\widetilde{S})\subset O(S)$,
and reflective type of $\widetilde{S}$ is dominated by reflective type of
$S$. For this case we say that lattices $S$ and $\widetilde{S}$ are
not equivariantly equivalent.
\endproclaim

\subhead
2.3. The classification of elliptically or parabolically
reflective elementary hyperbolic lattices of the rank 3
\endsubhead

\subsubhead
2.3.1. Notation
\endsubsubhead
According to Sect. 2.2, classification of reflective (of
any type) primitive elementary hyperbolic lattices is
divided in three parts:

1) Classification of main of them with
square-free determinant.

2) Classification of non-main lattices $\widetilde{S}$ of
them with square-free determinant where $S$ belongs to 1) and is odd
(see Proposition 2.2.6).

3) Taking of $m$-dual lattices to lattices $L$ from 1) and 2)
where $m\mid \det(L)$ (all of them will be reflective of the same type as
$L$).

For a hyperbolic lattice $S$, we denote by
$\M$ a fundamental polyhedron of $W(S)$ and by
$P(\M)_{\pr}$ the set of
orthogonal to $\M$ primitive roots from $S$. We denote by
$A(P(\M)_{\pr})$ the group of symmetries of $P(\M)_{\pr}$
(see \thetag{1.3.1}).

\definition{Definition 2.3.1.1} An automorphism
$\phi \in A(P(\M)_{\pr})$ is called a {\it central symmetry}
if it acts as a central symmetry in the hyperbolic space $\La(S)$.
Two central symmetries $\phi_1,\,\phi_2 \in A(P(\M)_{\pr})$ are equivalent
if they are conjugate in $A(P(\M)_{\pr})$ (i.e.
$\phi_2=\xi\phi_1\xi^{-1}$ where $\xi \in A(P(\M)_{\pr})$.
We denote by $h=h(S)$ the
{\it number of classes} of central symmetries in
$A(P(\M)_{\pr})$.
\enddefinition

This invariant is important because we have the following trivial

\proclaim{Lemma 2.3.1.2} If a hyperbolic lattice
$S$ is reflective of elliptic type,
then $h=h(S)\le 1$. If $S$ is reflective of parabolic type, then
$h=h(S)=0$. If $S$ is reflective of hyperbolic type and $\rk S=3$,
then $h=h(S)=0\ \text{or}\ 2$.
\endproclaim

For all reflective elementary hyperbolic lattices $S$,
we give the invariant $h=h(S)$, the set $P(\M)_{\pr}$ and the
Gram matrix $G(P(\M)_{\pr})=((\alpha_i,\alpha_j))$,
$\alpha_i,\,\alpha_j\in P(\M)_{\pr}$.

We denote by $\langle A \rangle$ a lattice with the matrix $A$.
Thus, one has $((e_i,e_j))=A$ for a bases $\{e_i\}$ of the lattice.
This bases is called {\it standard.}
We denote by $\oplus$ the orthogonal sum of lattices. Thus,
$\langle A \rangle \oplus \langle B \rangle$ denote a lattice with
the matrix $A\oplus B$. Let $A$ be a matrix of size $n\times n$
defining a lattice $\langle A \rangle$ with the standard bases
$\{e_i\}$. We denote by $\langle A \rangle
(c_{11},\dots,c_{1n}; \dots ;c_{m1},\dots,c_{mn})$
the overlattice of $\langle A \rangle$ generated by
$\{e_i\}$ and by elements $c_{11}e_1+\,\cdots\, +c_{1n}e_n$, ... ,
$c_{m1}e_1+\,\cdots\, +c_{mn}e_n$ of the dual lattice
$\langle A \rangle^\ast \subset \langle A \rangle \otimes \bq$
where $c_{ij} \in \bq$.

We denote by $U$ the lattice
$$
U=\left\langle \matrix0&1\\1&0\endmatrix\right\rangle.
\tag{2.3.1.1}
$$
It is unimodular even of signature $(1,1)$.

\subsubhead
2.3.2. Classification of elliptically or parabolically reflective
main hyperbolic lattices of rank $3$ and with square-free determinant
\endsubsubhead
We remind that a hyperbolic lattice $S$ of rank three and
with square-free determinant $d=\det(S)$ is called {\it main}
if it is even for even $d=\det(S)$. If $d$ is odd, then $S$ will
be necessarily odd.
By Proposition 2.2.4, a main hyperbolic lattice
$S$ of rank three and with square-free determinant
is defined by invariants $(d,\eta)$. We have the
following basic result:

\proclaim{Basic Theorem 2.3.2.1} Table 1 below gives the complete
list (containing 122 cases numerated by $N$) of main (i.e.
even for even determinant) hyperbolic lattices $S$ of rank $3$ and
with square-free determinant which are reflective of elliptic
or parabolic type. In Table 1, for each lattice $S$
we give invariants $(d,\eta)$ defining $S$ up to isomorphism,
number $h$ of classes of central symmetries,
matrix of $S$ for some standard bases, the set $P(\M)_{\pr}$ in
the standard bases and the Gram matrix $G(P(\M)_{\pr})$.
All these 122 cases have elliptic type.

If the lattice $S$ represents $0$ (i.e. there exists a non-zero
$x\in S$ with $x^2=0$), we give $S$ in the form
$S=U\oplus \langle -d \rangle $.
Then the fundamental polygon $\M$ is not compact
(i.e. $\M$ contains an infinite vertex).
In particular, Table 1
contains the complete list of reflective hyperbolic lattices of elliptic
or parabolic type of the form $S=U\oplus\langle -d \rangle$
where $d$ is square-free. There are $23$ these cases corresponding to
$$
\split 
d=&1,\,2,\,3,\,5,\,6,\,7,\,10,\,11,\,13,\,14,\,15,\,17,\,21,\,22,\,26,\,\\
  &30,\,33,\,34,\,38,\,42,\,66,\,78,\,110.
\endsplit 
$$
If $S$ is given in different form, the lattice $S$
does not represent $0$ and $\M$ is compact. There are $99$ these cases.
\endproclaim

\demo{Idea of the Proof} Here we give a rough scheme of the proof,
details will be given in Part II. First, we prove a formula
for the number $h$ of
classes of central symmetries of main hyperbolic lattices $S$
of rank $3$ with a square-free determinant. This formula uses
class-numbers of imaginary quadratic fields and
Legendre symbol (we shall give the formula in Part II).
By Lemma 2.3.1.2, elliptically or parabolically reflective $S$
have $h\le 1$. Using the formula for $h$,
we find all main hyperbolic lattices $S$ of rank $3$ with
square-free determinant $d\le 100000$ and $h\le 1$.
Their list is given in Table 3, it contains $206$ lattices $S$.
(In Table 3, we give invariants $(d,\eta)$, $h$ and a lattice
$S$ with these invariants.)
Conjectually, Table 3 gives the complete
(without the condition $d\le 100000$) list with the invariant
$h\le 1$ of main hyperbolic lattices $S$ of the rank $3$ and
with square-free determinant. To avoid this conjecture,
we use the same method (studying of narrow parts of fundamental
polyhedra $\M$)
we used to prove finiteness Theorem 1.4.1 (see \cite{N4},
\cite{N5}, \cite{N11} and also \cite{N8}). This method gives
a finite list of all possible invariants $(d,\eta)$
for the reflective lattices $S$.
Calculating their invariant $h$ (or estimating the invariant $d$),
we find that all of them having $h\le 1$ belong to Table 3.
Using Vinberg's algorithm \cite{V2}, we check reflective
type of all $206$ lattices $S$ of Table 3 and calculate their
fundamental polygon $\M$ for $W(S)$ if they are reflective.
For elliptic and parabolic type, the result is given in Table 1.
Actually, the method of narrow parts of polyhedra permits to find out
a big part of non-reflective lattices of Table 3. In Table 3 we
mark by $er$, $pr$, $hr$ and $nr$
elliptically, parabolically, hyperbolically reflective
type and non-reflective type respectively (there are no lattices of
parabolic type). It finishes the proof.
\enddemo

We mention an important

\proclaim{Corollary 2.3.2.2} The reflection groups
$W(S)$ of lattices $S$ of Table 1 contain all maximal
arithmetic reflection groups $W$ over $\bq$ in hyperbolic plane.
\endproclaim

\demo{Proof} By \'E.B. Vinberg \cite{V1}, $W$ is a
reflection subgroup of finite index $W\subset O(S_1)$ of a
hyperbolic lattice $S_1$ of rank $3$.
It follows that the fundamental polygon $\M$ of $W$
is elliptic, the group $A(P(\M)_{\pr})$ is finite and
has a generalized lattice Weyl vector $\rho$ with $(\rho,\rho)>0$.
By considerations in Sect. 2.2 and Theorem 2.3.2.1,
the group $W$ is a reflection subgroup $W\subset W(S)$ of
one of lattices $S$ of Table 1. Since $W$ is maximal, $W=W(S)$.
It proves the statement.
\enddemo

\remark{Remark 2.3.2.3} Reflective type of lattices
$U\oplus \langle -2k\rangle$ for all $k\le 60$ was found in
\cite{N13}. For $k=2k_1$ where
$k_1$ is odd, the lattice $U\oplus\langle -2k\rangle$
is equivariantly equivalent to the lattice
$U\oplus \langle -k_1 \rangle\cong \langle 1 \rangle \oplus
\langle -1 \rangle \oplus \langle -k_1 \rangle$. Thus, between lattices
$S=U\oplus \langle -d \rangle$ of Theorem 2.3.2.1, only the lattice
$S=U\oplus \langle -33 \rangle$ was not considered in \cite{N13}.
\endremark

\subsubhead
2.3.3. Classification of elliptically or parabolically
reflective non-main hyperbolic lattices of rank three and
with square-free determinant
\endsubsubhead

We remind that a lattice $S$ of rank three and with square-free
determinant is called
{\it non-main} if $S$ is odd but its determinant is even.
Using Basic Theorem 2.3.2.1, Proposition 2.2.6 and additional
considerations and calculations, we get

\proclaim{Theorem 2.3.3.1} Table 2 below gives the complete list
(containing $38$ cases numerated by $N^\prime$) of non-main
(i. e. odd with even determinant) hyperbolic lattices $\widetilde{S}$
with square-free determinant and of rank three which are
reflective of elliptic or parabolic type. In Table 2,
for each lattice $\widetilde{S}$
we give its invariants $(d,\text{odd},\eta)$, the number
$h$ of classes of central symmetries, matrix of $\widetilde{S}$ in
a standard bases, invariants $(d,\eta+\omega(d))$
of the main odd lattice $S$ of $\widetilde{S}$ (see Proposition
2.2.6). All these $38$ cases have elliptic type.

If $S$ and $\widetilde{S}$ are equivariantly
equivalent (in Table 2 there are $21$ these cases),
we only give the lattice $S$. For this case,
calculation of $P(\M)_{\pr}$ and
$G(P(\M)_{\pr})$ for $\widetilde{S}$ follows from calculation
of similar sets for $S$ in Table 1
(see Remark 2.3.3.2 below). In particular, all $\widetilde{S}$
representing $0$ (non-compact case) have this type, then
$\widetilde{S}=\langle 1 \rangle \oplus \langle -1 \rangle \oplus
\langle -d \rangle$ where
$$
d=2,\,6,\,10,\,14,\,22,\,26,\,30,\,34,\,42,\,66.
$$
If $S$ and $\widetilde{S}$ are not equivariantly equivalent
(there are 17 these cases), we (like for Table 1) give
$P(\M)_{\pr}$ and $G(P(\M)_{\pr})$ for $\widetilde{S}$ in
the standard bases.
\endproclaim

\remark{Remark 2.3.3.2}  Suppose that $\widetilde{S}$ is equivariantly
equivalent to $S$. Then $S(2)\subset \widetilde{S}$ is the unique
odd overlattice of $S(2)$. Let $P(\M)_{\pr}$ be the set of orthogonal to
$\M$ primitive roots of $S$ with the Gram matrix
$((\alpha, \beta ))$, $\alpha,\,\beta \in P(\M)_{\pr}$.
We denote by $\widetilde{P}(\M)_{\pr}$ the set of
orthogonal to $\M$ primitive roots of the overlattice
$S(2)\subset \widetilde{S}$. It is
$$
\widetilde{P}(\M)_{\pr}=
\{\widetilde{\alpha}=\alpha\ \mid\ \alpha\in P(\M)_{\pr}\  \and \
\alpha^2\equiv 1 \mod 2\} \cup
$$
$$
\cup
\{\widetilde{\alpha}=\alpha/2\ \mid\  \alpha \in P(\M)_{\pr}\ \and
\alpha^2\equiv 2\mod 4\}.
\tag{2.3.3.1}
$$
It follows that for
$\widetilde{\alpha},\,\widetilde{\beta}\in \widetilde{P}(\M)_{\pr}$ of
the lattice $\widetilde{S}$, one has
$$
(\widetilde{\alpha},\widetilde{\beta})=
\cases
2(\alpha, \beta), &\text{if $\alpha^2\beta^2\equiv 1\mod 2$} \\
(\alpha,\beta), &\text{if $\alpha^2\beta^2\equiv 2\mod 4$}\\
(\alpha,\beta)/2,&\text{if $\alpha^2\equiv \beta^2 \equiv 2 \mod 4$}
\endcases
\tag{2.3.3.2}
$$
which describes the Gram matrix of $\widetilde{P}(\M)_{\pr}$ for the
lattice $\widetilde{S}$ using the Gram matrix of $P(\M)_{\pr}$ for the
lattice $S$.
Here on the left hand side of \thetag{2.3.3.2}
we use the form of the lattice $\widetilde{S}$, and on the right hand
side of \thetag{2.3.3.2} we use the form of the lattice $S$.
\endremark

\subsubhead
2.3.4. Classification of elementary hyperbolic lattices
of rank three which are reflective of elliptic or parabolic type
\endsubsubhead

We remind that a lattice $L$ is called {\it elementary} if
its discriminant group $A_L=L^\ast /L$ is elementary: any
its $p$-component is a $p$-elementary Abelian group:
$A_{L_p}\cong (\bz/p\bz)^{r_p}$.

By Sect. 2.2, we get

\proclaim{Theorem 2.3.4.1} A primitive elementary hyperbolic lattice
$F$ is reflective of elliptic or parabolic type if and only if
$$
F\cong S^{\ast,m}(m),\ m\mid d
\tag{2.3.4.1}
$$
(see \thetag{2.2.5}) where $S$ is a lattice of the determinant $d=\det(S)$
(it is square-free) of Tables 1 or 2. Then $\det(F)=dm$.
In particular, for the fixed lattice $S$ the number
of lattices $F$ is equal to $2^t$ where $t$ is
the number of prime divisors of $d$.

Let $P(\M)_{\pr}$ be the set of orthogonal to $\M$ primitive roots of
$S$. Then the set $\widetilde{P}(\M)_{\pr}$ of orthogonal
to $\M$ primitive roots of $F$ is equal to
$$
\widetilde{P}(\M)_{\pr}=
\{\widetilde{\alpha}=\alpha/k_{\alpha,m}\ \mid \
\alpha \in P(\M)_{\pr}\}
\tag{2.3.4.2}
$$
where $k_{\alpha,m}$ is the greatest divisor of $m$ such that
$\alpha/k_{\alpha,m}\in S^\ast$ (equivalently, if
$(\alpha,S)=t\bz$, then $k_{\alpha,m}=\text{g.c.d}(t,m)$).
The Gram matrix of $\widetilde{P}(\M)_{\pr}$ for the lattice $F$
is equal to
$$
\left((\widetilde{\alpha},\widetilde{\beta})\right)=
\left({(\alpha,\beta)m\over k_{\alpha,m}k_{\beta,m}}\right),\ \
\widetilde{\alpha},\,\widetilde{\beta} \in \widetilde{P}(\M)_{\pr}.
\tag{2.3.4.3}
$$
Here, on the left hand side of \thetag{2.3.4.3}
we use the form of the
lattice $F$, and on the right hand side of \thetag{2.3.4.3}
we use the form of the lattice $S$.
\endproclaim

\newpage

\centerline{\bf Table 1.}
\centerline{The list of main hyperbolic lattices}
\centerline{of rank three and with square-free determinant}
\centerline{which are reflective of elliptic or parabolic type.}

\vskip20pt

\vbox{\noindent
$N=1\ $ $d=1,\,\eta=0,\,h=0$: $U\oplus \langle -1 \rangle$
\nobreak
\newline
$\pmatrix{1}&{0}&{-1}\cr{0}&{0}&{1}\cr{-1}&{1}&{0}\cr
\endpmatrix
\hskip20pt
\pmatrix
{-1}&{1}&{1}\cr{1}&{-1}&{0}\cr{1}&{0}&{-2}\cr\endpmatrix$.}

\vbox{\noindent
$N=2\ $ $d=2,\,\eta=0,\,h=0$: $U\oplus \langle -2 \rangle$
\nobreak
\newline
$\pmatrix{1}&{0}&{-1}\cr{0}&{0}&{1}\cr{-1}&{1}&{0}\cr\endpmatrix
\hskip20pt
\pmatrix{-2}&{2}&{1}\cr{2}&{-2}&{0}\cr{1}&{0}&{-2}\cr\endpmatrix$.}

\vbox{\noindent
$N=3\ $ $d=3,\,\eta=0,\,h=0$:
$\langle 3 \rangle\oplus \langle -1\rangle \oplus \langle -1 \rangle $
\nobreak
\newline
$\pmatrix{0}&{-1}&{1}\cr{0}&{1}&{0}\cr
{1}&{0}&{-3}\cr\endpmatrix
\hskip20pt
\pmatrix{-2}&{1}&{3}\cr{1}&{-1}&{0}\cr
{3}&{0}&{-6}\cr\endpmatrix$.}

\vbox{\noindent
$N=4\ $ $d=3,\,\eta=1,\,h=0$: $U\oplus \langle -3 \rangle$
\nobreak
\newline
$\pmatrix{3}&{0}&{-1}\cr{0}&{0}&{1}\cr{-1}&{1}&{0}\cr
{1}&{1}&{-1}\cr\endpmatrix
\hskip20pt
\pmatrix{-3}&{3}&{3}&{0}\cr{3}&{-3}&{0}&{3}\cr{3}&{0}&{-2}&{0}\cr
{0}&{3}&{0}&{-1}\cr\endpmatrix$.}

\vbox{\noindent
$N=5\ $ $d=5,\,\eta=0,\,h=0$: $U\oplus \langle -5 \rangle$
\nobreak
\newline
$\pmatrix{5}&{0}&{-1}\cr{0}&{0}&{1}\cr{-1}&{1}&{0}\cr
{2}&{1}&{-1}\cr\endpmatrix
\hskip20pt
\pmatrix{-5}&{5}&{5}&{0}\cr{5}&{-5}&{0}&{5}\cr
{5}&{0}&{-2}&{1}\cr{0}&{5}&{1}&{-1}\cr\endpmatrix$.}

\vbox{\noindent
$N=6\ $ $d=5,\, \eta=1,\,h=0$:
$\langle 1 \rangle \oplus \langle -10 \rangle \oplus
\langle -2 \rangle (0,1/2,1/2)$
\nobreak
\newline
$\pmatrix{0}&{0}&{1}\cr{0}&{1}&{0}\cr{1}&{0}&{-1}\cr{1}&{{-1}\over{2}}
&{{-1}\over{2}}\cr\endpmatrix
\hskip20pt
\pmatrix{-2}&{0}&{2}&{1}\cr{0}&{-10}&{0}&{5}\cr
{2}&{0}&{-1}&{0}\cr{1}&{5}&{0}&{-2}\cr\endpmatrix$.}

\vbox{\noindent
$N=7\ $ $d=6=2\cdot 3,\,\eta=0,h=0$: $U\oplus \langle -6 \rangle$
\nobreak
\newline
$\pmatrix{3}&{0}&{-1}\cr{0}&{0}&{1}\cr{-1}&{1}&{0}\cr\endpmatrix
\hskip20pt
\pmatrix{-6}&{6}&{3}\cr{6}&{-6}&{0}\cr{3}&{0}&{-2}\cr\endpmatrix$.}

\vbox{\noindent
$N=8\ $ $d=7,\,\eta=0,\,h=1$:
$\langle 7 \rangle \oplus \langle -1 \rangle \oplus \langle -1\rangle $
\nobreak
\newline
$\pmatrix
{0}&{-1}&{1}\cr
{0}&{1}&{0}\cr
{1}&{0}&{-3}\cr
{1}&{-2}&{-2}\cr
\endpmatrix
\hskip20pt
\pmatrix
{-2}&{1}&{3}&{0}\cr
{1}&{-1}&{0}&{2}\cr
{3}&{0}&{-2}&{1}\cr
{0}&{2}&{1}&{-1}\cr
\endpmatrix$.}

\vbox{\noindent
$N=9\ $ $d=7,\,\eta=1,\,h=0$: $U\oplus \langle -7 \rangle$
\nobreak
\newline
$\pmatrix{7}&{0}&{-1}\cr{0}&{0}&{1}\cr{-1}&{1}&{0}\cr
{7}&{7}&{-4}\cr{3}&{1}&{-1}\cr\endpmatrix
\hskip20pt
\pmatrix{-7}&{7}&{7}&{21}&{0}\cr{7}&{-7}&{0}&{28}&{7}\cr
{7}&{0}&{-2}&{0}&{2}\cr{21}&{28}&{0}&{-14}&{0}\cr
{0}&{7}&{2}&{0}&{-1}\cr\endpmatrix$.}

\vbox{\noindent
$N=10\ $ $d=10=2\cdot 5,\,\eta=1,h=0$: $U\oplus \langle -10 \rangle$
\nobreak
\newline
$\pmatrix{5}&{0}&{-1}\cr{0}&{0}&{1}\cr{-1}&{1}&{0}\cr{2}&{2}&{-1}\cr
\endpmatrix
\hskip20pt
\pmatrix{-10}&{10}&{5}&{0}\cr{10}&{-10}&{0}&{10}\cr
{5}&{0}&{-2}&{0}\cr{0}&{10}&{0}&{-2}\cr\endpmatrix$.}

\vbox{\noindent
$N=11\ $ $d=11,\,\eta=0,\,h=0$:
$\langle 11 \rangle\oplus \langle -1\rangle \oplus \langle -1\rangle$
\nobreak
\newline
$\pmatrix
{0}&{-1}&{1}\cr{0}&{1}&{0}\cr{3}&{0}&{-11}\cr{1}&{-2}&{-3}\cr
\endpmatrix
\hskip20pt
\pmatrix
{-2}&{1}&{11}&{1}\cr
{1}&{-1}&{0}&{2}\cr
{11}&{0}&{-22}&{0}
\cr{1}&{2}&{0}&{-2}\cr
\endpmatrix$.}

\vbox{\noindent
$N=12\ $ $d=11,\,\eta=1,\,h=1$: $U\oplus \langle -11 \rangle$
\nobreak
\newline
$\pmatrix{11}&{0}&{-1}\cr{0}&{0}&{1}\cr{-1}&{1}&{0}\cr{7}&{7}&{-3}\cr
{55}&{44}&{-21}\cr{66}&{44}&{-23}\cr
{7}&{3}&{-2}\cr{5}&{1}&{-1}\cr\endpmatrix
\hskip20pt
\pmatrix{-11}&{11}&{11}&{44}&{253}&{231}&{11}&{0}\cr
{11}&{-11}&{0}&{33}&{231}&{253}&{22}&{11}\cr
{11}&{0}&{-2}&{0}&{11}&{22}&{4}&{4}\cr
{44}&{33}&{0}&{-1}&{0}&{11}&{4}&{9}\cr
{253}&{231}&{11}&{0}&{-11}&{11}&{11}&{44}\cr
{231}&{253}&{22}&{11}&{11}&{-11}&{0}&{33}\cr
{11}&{22}&{4}&{4}&{11}&{0}&{-2}&{0}\cr
{0}&{11}&{4}&{9}&{44}&{33}&{0}&{-1}\cr\endpmatrix$.}

\vbox{\noindent
$N=13\ $ $d=13,\,\eta =0,\,h=1$: $U\oplus \langle -13 \rangle$
\nobreak
\newline
$\pmatrix{13}&{0}&{-1}\cr{0}&{0}&{1}\cr{-1}&{1}&{0}\cr
{3}&{2}&{-1}\cr{104}&{39}&{-25}\cr{234}&{78}&{-53}\cr
{55}&{17}&{-12}\cr{6}&{1}&{-1}\cr\endpmatrix
\hskip20pt
\pmatrix{-13}&{13}&{13}&{13}&{182}&{325}&{65}&{0}\cr
{13}&{-13}&{0}&{13}&{325}&{689}&{156}&{13}\cr
{13}&{0}&{-2}&{1}&{65}&{156}&{38}&{5}\cr
{13}&{13}&{1}&{-1}&{0}&{13}&{5}&{2}\cr
{182}&{325}&{65}&{0}&{-13}&{13}&{13}&{13}\cr
{325}&{689}&{156}&{13}&{13}&{-13}&{0}&{13}\cr
{65}&{156}&{38}&{5}&{13}&{0}&{-2}&{1}\cr
{0}&{13}&{5}&{2}&{13}&{13}&{1}&{-1}\cr\endpmatrix$.}

\vbox{\noindent
$N=14\ $ $d=13,\,\eta=1,\,h=1$:
$\langle 26 \rangle\oplus \langle -2\rangle \oplus\langle -1\rangle
(1/2,1/2,0)$
\nobreak
\newline
$\pmatrix{0}&{0}&{1}\cr{0}&{1}&{0}\cr{5}&{0}&{-26}\cr
{1}&{-1}&{-5}\cr{{1}\over{2}}&{{-3}\over{2}}&{-2}\cr
{{3}\over{2}}&{{-13}\over{2}}&{0}\cr\endpmatrix
\hskip20pt
\pmatrix{-1}&{0}&{26}&{5}&{2}&{0}\cr
{0}&{-2}&{0}&{2}&{3}&{13}\cr
{26}&{0}&{-26}&{0}&{13}&{195}\cr
{5}&{2}&{0}&{-1}&{0}&{26}\cr
{2}&{3}&{13}&{0}&{-2}&{0}\cr
{0}&{13}&{195}&{26}&{0}&{-26}\cr\endpmatrix$.}

\vbox{\noindent
$N=15\ $ $d=14=2\cdot 7,\,\eta=1,\,h=0$: $U\oplus \langle -14 \rangle$
\nobreak
\newline
$\pmatrix{7}&{0}&{-1}\cr{0}&{0}&{1}\cr{-1}&{1}&{0}\cr{3}&{2}&{-1}\cr
\endpmatrix
\hskip20pt
\pmatrix{-14}&{14}&{7}&{0}\cr{14}&{-14}&{0}&{14}\cr{7}&{0}&{-2}&{1}\cr
{0}&{14}&{1}&{-2}\cr\endpmatrix$.}

\vbox{\noindent
$N=16\ $ $d=15=3\cdot 5,\,\eta=0,\,h=0$:
$\langle 3\rangle \oplus \langle -5\rangle\oplus \langle -1\rangle$
\nobreak
\newline
$\pmatrix{0}&{0}&{1}\cr{0}&{1}&{0}\cr{1}&{0}&{-3}\cr
{1}&{-1}&{0}\cr\endpmatrix
\hskip20pt
\pmatrix
{-1}&{0}&{3}&{0}\cr
{0}&{-5}&{0}&{5}\cr
{3}&{0}&{-6}&{3}\cr
{0}&{5}&{3}&{-2}\cr\endpmatrix$.}

\vbox{\noindent
$N=17\ $ $d=15=3\cdot 5,\,\eta=1,\,h=0$: $\langle 5\rangle \oplus
\langle -3 \rangle \oplus \langle -1 \rangle $
\nobreak
\newline
$\pmatrix{0}&{0}&{1}\cr{0}&{1}&{0}\cr{2}&{0}&{-5}
\cr{1}&{-1}&{-2}\cr{3}&{-5}&{0}\cr\endpmatrix
\hskip20pt
\pmatrix
{-1}&{0}&{5}&{2}&{0}\cr
{0}&{-3}&{0}&{3}&{15}\cr
{5}&{0}&{-5}&{0}&{30}\cr
{2}&{3}&{0}&{-2}&{0}\cr
{0}&{15}&{30}&{0}&{-30}\cr\endpmatrix$.}

\vbox{\noindent
$N=18\ $ $d=15=3\cdot 5,\,\eta=2,\,h=0$: $U\oplus \langle -15 \rangle$
\nobreak
\newline
$\pmatrix{15}&{0}&{-1}\cr{0}&{0}&{1}\cr{-1}&{1}&{0}\cr
{5}&{5}&{-2}\cr{9}&{3}&{-2}\cr{7}&{1}&{-1}\cr\endpmatrix
\hskip20pt
\pmatrix{-15}&{15}&{15}&{45}&{15}&{0}\cr
{15}&{-15}&{0}&{30}&{30}&{15}\cr
{15}&{0}&{-2}&{0}&{6}&{6}\cr
{45}&{30}&{0}&{-10}&{0}&{10}\cr
{15}&{30}&{6}&{0}&{-6}&{0}\cr
{0}&{15}&{6}&{10}&{0}&{-1}\cr\endpmatrix$.}

\vbox{\noindent
$N=19\ $ $d=15=3\cdot 5,\,\eta=3,\,h=0$: $\langle 15 \rangle \oplus
\langle -1 \rangle \oplus \langle -1\rangle$
\nobreak
\newline
$\pmatrix
{0}&{-1}&{1}\cr{0}&{1}&{0}\cr{1}&{0}&{-5}\cr
{1}&{-3}&{-3}\cr\endpmatrix
\hskip20pt
\pmatrix
{-2}&{1}&{5}&{0}\cr{1}&{-1}&{0}&{3}\cr
{5}&{0}&{-10}&{0}\cr{0}&{3}&{0}&{-3}\cr\endpmatrix$.}

\vbox{\noindent
$N=20\ $ $d=17,\,\eta=0,\,h=0$: $U\oplus \langle -17 \rangle$
\nobreak
\newline
$\pmatrix{17}&{0}&{-1}\cr{0}&{0}&{1}\cr{-1}&{1}&{0}\cr
{17}&{17}&{-6}\cr{4}&{2}&{-1}\cr{11}&{3}&{-2}\cr
{8}&{1}&{-1}\cr\endpmatrix
\hskip20pt
\pmatrix{-17}&{17}&{17}&{187}&{17}&{17}&{0}\cr
{17}&{-17}&{0}&{102}&{17}&{34}&{17}\cr
{17}&{0}&{-2}&{0}&{2}&{8}&{7}\cr
{187}&{102}&{0}&{-34}&{0}&{34}&{51}\cr
{17}&{17}&{2}&{0}&{-1}&{0}&{3}\cr
{17}&{34}&{8}&{34}&{0}&{-2}&{1}\cr
{0}&{17}&{7}&{51}&{3}&{1}&{-1}\cr\endpmatrix$.}

\vbox{\noindent
$N=21\ $ $d=19,\,\eta=0,\,h=1$:
$\langle 19 \rangle \oplus \langle -1\rangle \oplus
\langle -1\rangle $
\nobreak
\newline
$\pmatrix{0}&{-1}&{1}\cr{0}&{1}&{0}\cr{13}&{0}&{-57}\cr
{3}&{-2}&{-13}\cr{1}&{-2}&{-4}\cr{6}&{-19}&{-19}\cr\endpmatrix
\hskip20pt
\pmatrix
{-2}&{1}&{57}&{11}&{2}&{0}\cr
{1}&{-1}&{0}&{2}&{2}&{19}\cr
{57}&{0}&{-38}&{0}&{19}&{399}\cr
{11}&{2}&{0}&{-2}&{1}&{57}\cr
{2}&{2}&{19}&{1}&{-1}&{0}\cr
{0}&{19}&{399}&{57}&{0}&{-38}\cr\endpmatrix$.}

\vbox{\noindent
$N=22\ $ $d=21=3\cdot 7,\,\eta=0,\,h=0$:
$\langle 14\rangle \oplus \langle -6\rangle \oplus
\langle -1\rangle(1/2,1/2,0)$
\nobreak
\newline
$\pmatrix{0}&{0}&{1}\cr{0}&{1}&{0}\cr{{1}\over{2}}&
{{-1}\over{2}}&{-2}\cr{{3}\over{2}}&{{-7}\over{2}}&{0}\cr\endpmatrix
\hskip20pt
\pmatrix{-1}&{0}&{2}&{0}\cr{0}&{-6}&{3}&{21}\cr
{2}&{3}&{-2}&{0}\cr{0}&{21}&{0}&{-42}\cr\endpmatrix$.}

\vbox{\noindent
$N=23\ $ $d=21=3\cdot 7,\,\eta=1,\,h=1$: $U\oplus \langle -21 \rangle$
\nobreak
\newline
$\left(\smallmatrix{21}&{0}&{-1}\cr{0}&{0}&{1}\cr{-1}&{1}&{0}\cr
{3}&{3}&{-1}\cr{5}&{2}&{-1}\cr{84}&{21}&{-13}\cr
{210}&{42}&{-29}\cr{61}&{11}&{-8}\cr{57}&{9}&{-7}\cr
{10}&{1}&{-1}\cr\endsmallmatrix\right)
\left(\smallmatrix{-21}&{21}&{21}&{42}&{21}&{168}&{273}&{63}&{42}&{0}\cr
{21}&{-21}&{0}&{21}&{21}&{273}&{609}&{168}&{147}&{21}\cr
{21}&{0}&{-2}&{0}&{3}&{63}&{168}&{50}&{48}&{9}\cr
{42}&{21}&{0}&{-3}&{0}&{42}&{147}&{48}&{51}&{12}\cr
{21}&{21}&{3}&{0}&{-1}&{0}&{21}&{9}&{12}&{4}\cr
{168}&{273}&{63}&{42}&{0}&{-21}&{21}&{21}&{42}&{21}\cr
{273}&{609}&{168}&{147}&{21}&{21}&{-21}&{0}&{21}&{21}\cr
{63}&{168}&{50}&{48}&{9}&{21}&{0}&{-2}&{0}&{3}\cr
{42}&{147}&{48}&{51}&{12}&{42}&{21}&{0}&{-3}&{0}\cr
{0}&{21}&{9}&{12}&{4}&{21}&{21}&{3}&{0}&{-1}\cr\endsmallmatrix\right)$.}

\vbox{\noindent
$N=24\ $ $d=21=3 \cdot 7,\,\eta=2,\,h=0$:
$\langle 3\rangle \oplus \langle -7\rangle \oplus \langle-1\rangle$
\nobreak
\newline
$\pmatrix
{0}&{0}&{1}\cr{0}&{1}&{0}\cr{1}&{0}&{-3}\cr{7}&{-4}&{-7}\cr
{3}&{-2}&{-1}\cr\endpmatrix
\hskip20pt
\pmatrix{-1}&{0}&{3}&{7}&{1}\cr{0}&{-7}&{0}&{28}&{14}\cr
{3}&{0}&{-6}&{0}&{6}\cr{7}&{28}&{0}&{-14}&{0}\cr
{1}&{14}&{6}&{0}&{-2}\cr\endpmatrix$.}

\vbox{\noindent
$N=25\ $ $d=21=3\cdot 7,\,\eta=3,\,h=0$:
$\langle 1 \rangle \oplus \langle -7\rangle \oplus
\langle -3\rangle$
\nobreak
\newline
$\pmatrix
{0}&{0}&{1}\cr{0}&{1}&{0}\cr{1}&{0}&{-1}\cr
{3}&{-1}&{-1}\cr{7}&{-3}&{0}\cr\endpmatrix
\hskip20pt
\pmatrix{-3}&{0}&{3}&{3}&{0}\cr{0}&{-7}&{0}&{7}&{21}\cr
{3}&{0}&{-2}&{0}&{7}\cr{3}&{7}&{0}&{-1}&{0}\cr
{0}&{21}&{7}&{0}&{-14}\cr\endpmatrix$.}

\vbox{\noindent
$N=26\ $ $d=22=2\cdot 11,\,\eta=0,\,h=1$: $U\oplus \langle -22 \rangle$
\nobreak
\newline
$\pmatrix{11}&{0}&{-1}\cr{0}&{0}&{1}\cr{-1}&{1}&{0}\cr{33}&{33}&{-10}\cr
{88}&{66}&{-23}\cr{5}&{2}&{-1}\cr\endpmatrix
\hskip20pt
\pmatrix{-22}&{22}&{11}&{143}&{220}&{0}\cr{22}&{-22}&{0}&{220}&{506}&{22}\cr
{11}&{0}&{-2}&{0}&{22}&{3}\cr{143}&{220}&{0}&{-22}&{22}&{11}\cr
{220}&{506}&{22}&{22}&{-22}&{0}\cr{0}&{22}&{3}&{11}&{0}&{-2}\cr
\endpmatrix$.}

\vbox{\noindent
$N=27\ $ $d=23,\,\eta=0,\,h=1$:
$\langle 23 \rangle \oplus \langle -1 \rangle \oplus \langle -1 \rangle$
\nobreak
\newline
$\pmatrix{0}&{-1}&{1}\cr{0}&{1}&{0}\cr{1}&{0}&{-5}\cr
{12}&{-17}&{-55}\cr{6}&{-10}&{-27}\cr{1}&{-3}&{-4}\cr\endpmatrix
\hskip20pt
\pmatrix
{-2}&{1}&{5}&{38}&{17}&{1}\cr
{1}&{-1}&{0}&{17}&{10}&{3}\cr
{5}&{0}&{-2}&{1}&{3}&{3}\cr
{38}&{17}&{1}&{-2}&{1}&{5}\cr
{17}&{10}&{3}&{1}&{-1}&{0}\cr
{1}&{3}&{3}&{5}&{0}&{-2}\cr\endpmatrix$.}

\vbox{\noindent
$N=28\ $ $d=26=2\cdot 13,\,\eta=1,\,h=0$: $U\oplus \langle -26 \rangle$
\nobreak
\newline
$\pmatrix{13}&{0}&{-1}\cr{0}&{0}&{1}\cr{-1}&{1}&{0}\cr{4}&{3}&{-1}\cr
{6}&{2}&{-1}\cr\endpmatrix
\hskip20pt
\pmatrix{-26}&{26}&{13}&{13}&{0}\cr
{26}&{-26}&{0}&{26}&{26}\cr{13}&{0}&{-2}&{1}&{4}\cr
{13}&{26}&{1}&{-2}&{0}\cr{0}&{26}&{4}&{0}&{-2}\cr\endpmatrix$.}

\vbox{\noindent
$N=29\ $ $d=29,\,\eta=1,\,h=1$:
$\langle 1 \rangle \oplus \langle -58 \rangle \oplus
\langle -2 \rangle (0,1/2,1/2)$
\nobreak
\newline
$\pmatrix{0}&{0}&{1}\cr{0}&{1}&{0}\cr{1}&{0}&{-1}\cr
{5}&{{-1}\over{2}}&{{-5}\over{2}}\cr{8}&{-1}&{-2}\cr
{232}&{-30}&{-29}\cr{23}&{-3}&{-2}\cr
{19}&{{-5}\over{2}}&{{-1}\over{2}}\cr\endpmatrix
\hskip20pt
\pmatrix{-2}&{0}&{2}&{5}&{4}&{58}&{4}&{1}\cr
{0}&{-58}&{0}&{29}&{58}&{1740}&{174}&{145}\cr
{2}&{0}&{-1}&{0}&{4}&{174}&{19}&{18}\cr
{5}&{29}&{0}&{-2}&{1}&{145}&{18}&{20}\cr
{4}&{58}&{4}&{1}&{-2}&{0}&{2}&{5}\cr
{58}&{1740}&{174}&{145}&{0}&{-58}&{0}&{29}\cr
{4}&{174}&{19}&{18}&{2}&{0}&{-1}&{0}\cr
{1}&{145}&{18}&{20}&{5}&{29}&{0}&{-2}\cr\endpmatrix$.}

\vbox{\noindent
$N=30\ $ $d=30=2\cdot 3 \cdot 5,\,\eta=1,\,h=1$:
$U\oplus \langle -30 \rangle$
\nobreak
\newline
$\pmatrix{15}&{0}&{-1}\cr{0}&{0}&{1}\cr{-1}&{1}&{0}\cr
{15}&{15}&{-4}\cr{60}&{30}&{-11}\cr{7}&{2}&{-1}\cr\endpmatrix
\hskip20pt
\pmatrix{-30}&{30}&{15}&{105}&{120}&{0}\cr
{30}&{-30}&{0}&{120}&{330}&{30}\cr{15}&{0}&{-2}&{0}&{30}&{5}\cr
{105}&{120}&{0}&{-30}&{30}&{15}\cr{120}&{330}&{30}&{30}&{-30}&{0}\cr
{0}&{30}&{5}&{15}&{0}&{-2}\cr\endpmatrix$.}

\vbox{\noindent
$N=31\ $ $d=30=2\cdot 3 \cdot 5,\,\eta=2,\,h=0$:
$\langle 10\rangle \oplus \langle -6\rangle \oplus
\langle -2\rangle(0,1/2,1/2)$
\nobreak
\newline
$\pmatrix{0}&{{-1}\over{2}}&{{1}\over{2}}\cr{0}&{1}&{0}\cr
{2}&{0}&{-5}\cr{3}&{-2}&{-6}\cr\endpmatrix
\hskip20pt
\pmatrix{-2}&{3}&{5}&{0}\cr{3}&{-6}&{0}&{12}\cr
{5}&{0}&{-10}&{0}\cr{0}&{12}&{0}&{-6}\cr\endpmatrix$.}

\vbox{\noindent
$N=32\ $ $d=33=3\cdot 11,\,\eta=0,\,h=1$:
$\langle 3\rangle \oplus \langle -22\rangle
\oplus \langle -2 \rangle (0,1/2,1/2)$
\nobreak
\newline
$\pmatrix{0}&{0}&{1}\cr{0}&{1}&{0}\cr{2}&{0}&{-3}\cr
{6}&{{-3}\over{2}}&{{-11}\over{2}}\cr
{22}&{{-13}\over{2}}&{{-33}\over{2}}\cr
{4}&{{-3}\over{2}}&{{-3}\over{2}}\cr\endpmatrix
\hskip20pt
\pmatrix{-2}&{0}&{6}&{11}&{33}&{3}\cr{0}&{-22}&{0}&{33}&{143}&{33}\cr
{6}&{0}&{-6}&{3}&{33}&{15}\cr{11}&{33}&{3}&{-2}&{0}&{6}\cr
{33}&{143}&{33}&{0}&{-22}&{0}\cr{3}&{33}&{15}&{6}&{0}&{-6}\cr
\endpmatrix$.}

\vbox{\noindent
$N=33\ $ $d=33=3\cdot 11,\,\eta=1,\,h=0$:
$\langle 11\rangle \oplus \langle -3\rangle \oplus \langle -1\rangle$
\nobreak
\newline
$\pmatrix{0}&{0}&{1}\cr{0}&{1}&{0}\cr{3}&{0}&{-11}\cr
{1}&{-1}&{-3}\cr{15}&{-22}&{-33}\cr{1}&{-2}&{-1}\cr\endpmatrix
\hskip20pt
\pmatrix{-1}&{0}&{11}&{3}&{33}&{1}\cr
{0}&{-3}&{0}&{3}&{66}&{6}\cr{11}&{0}&{-22}&{0}&{132}&{22}\cr
{3}&{3}&{0}&{-1}&{0}&{2}\cr{33}&{66}&{132}&{0}&{-66}&{0}\cr
{1}&{6}&{22}&{2}&{0}&{-2}\cr\endpmatrix $.}

\vbox{\noindent
$N=34\ $ $d=33=3\cdot 11,\,\eta=2,\,h=1$: $U\oplus \langle -33 \rangle$
\nobreak
\newline
$\left(\smallmatrix{33}&{0}&{-1}\cr{0}&{0}&{1}\cr{-1}&{1}&{0}\cr
{4}&{4}&{-1}\cr{39}&{27}&{-8}\cr{407}&{253}&{-79}\cr
{68}&{41}&{-13}\cr{231}&{132}&{-43}\cr{132}&{66}&{-23}\cr
{13}&{5}&{-2}\cr{8}&{2}&{-1}\cr{21}&{3}&{-2}\cr
{121}&{11}&{-9}\cr{16}&{1}&{-1}\cr\endsmallmatrix\right)
\left(\smallmatrix{-33}&{33}&{33}&{99}&{627}&{5742}&{924}&{2937}&{1419}&
{99}&{33}&{33}&{66}&{0}\cr
{33}&{-33}&{0}&{33}&{264}&{2607}&{429}&{1419}&{759}&{66}&{33}&
{66}&{297}&{33}\cr
{33}&{0}&{-2}&{0}&{12}&{154}&{27}&{99}&{66}&{8}&{6}&{18}&{110}&
{15}\cr
{99}&{33}&{0}&{-1}&{0}&{33}&{7}&{33}&{33}&{6}&{7}&{30}&{231}&{35}\cr
{627}&{264}&{12}&{0}&{-6}&{0}&{3}&{33}&{66}&{18}&{30}&{156}&{1320}&{207}\cr
{5742}&{2607}&{154}&{33}&{0}&{-11}&{0}&{66}&{297}&{110}&{231}&{1320}&
{11627}&{1848}\cr
{924}&{429}&{27}&{7}&{3}&{0}&{-1}&{0}&{33}&{15}&{35}&{207}&{1848}&{295}\cr
{2937}&{1419}&{99}&{33}&{33}&{66}&{0}&{-33}&{33}&{33}&{99}&
{627}&{5742}&{924}\cr{1419}&{759}&{66}&{33}&{66}&{297}&{33}&{33}&{-33}&{0}
&{33}&{264}&{2607}&{429}\cr
{99}&{66}&{8}&{6}&{18}&{110}&{15}&{33}&{0}&{-2}&{0}&{12}&{154}&{27}\cr
{33}&{33}&{6}&{7}&{30}&{231}&{35}&{99}&{33}&{0}&
{-1}&{0}&{33}&{7}\cr
{33}&{66}&{18}&{30}&{156}&{1320}&{207}&{627}&{264}
&{12}&{0}&{-6}&{0}&{3}\cr
{66}&{297}&{110}&{231}&{1320}&{11627}&{1848}&{5742}&{2607}&
{154}&{33}&{0}&{-11}&{0}\cr
{0}&{33}&{15}&{35}&{207}&{1848}&{295}&{924}&{429}&{27}&{7}&
{3}&{0}&{-1}\cr\endsmallmatrix\right)$.}

\vbox{\noindent
$N=35\ $ $d=33=3\cdot 11,\,\eta=3,\,h=1$:
$\langle 1\rangle \oplus \langle -11\rangle
\oplus \langle -3\rangle$
\nobreak
\newline
$\pmatrix{0}&{0}&{1}\cr{0}&{1}&{0}\cr{1}&{0}&{-1}\cr
{12}&{-3}&{-4}\cr{44}&{-12}&{-11}\cr{3}&{-1}&{0}\cr\endpmatrix
\hskip20pt
\pmatrix
{-3}&{0}&{3}&{12}&{33}&{0}\cr{0}&{-11}&{0}&{33}&{132}&{11}\cr
{3}&{0}&{-2}&{0}&{11}&{3}\cr{12}&{33}&{0}&{-3}&{0}&{3}\cr
{33}&{132}&{11}&{0}&{-11}&{0}\cr{0}&{11}&{3}&{3}&{0}&{-2}\cr
\endpmatrix$.}

\vbox{\noindent
$N=36\ $ $d=34=2\cdot 17,\,\eta=0,\,h=1$: $U\oplus \langle -34 \rangle$
\nobreak
\newline
$\pmatrix{17}&{0}&{-1}\cr{0}&{0}&{1}\cr{-1}&{1}&{0}\cr{4}&{4}&{-1}\cr
{85}&{51}&{-16}\cr{204}&{102}&{-35}\cr{19}&{8}&{-3}\cr{8}&{2}&{-1}\cr
\endpmatrix
\hskip20pt
\pmatrix{-34}&{34}&{17}&{34}&{323}&{544}&{34}&{0}\cr
{34}&{-34}&{0}&{34}&{544}&{1190}&{102}&{34}\cr
{17}&{0}&{-2}&{0}&{34}&{102}&{11}&{6}\cr
{34}&{34}&{0}&{-2}&{0}&{34}&{6}&{6}\cr
{323}&{544}&{34}&{0}&{-34}&{34}&{17}&{34}\cr
{544}&{1190}&{102}&{34}&{34}&{-34}&{0}&{34}\cr
{34}&{102}&{11}&{6}&{17}&{0}&{-2}&{0}\cr
{0}&{34}&{6}&{6}&{34}&{34}&{0}&{-2}\cr\endpmatrix$.}

\vbox{\noindent
$N=37\ $ $d=35=5\cdot 7,\,\eta=0,\,h=0$:
$\langle 7 \rangle\oplus \langle -5 \rangle \oplus \langle -1 \rangle$
\nobreak
\newline
$\pmatrix{0}&{0}&{1}\cr{0}&{1}&{0}\cr{1}&{0}&{-3}\cr
{1}&{-1}&{-2}\cr{5}&{-7}&{0}\cr\endpmatrix
\hskip20pt
\pmatrix{-1}&{0}&{3}&{2}&{0}\cr{0}&{-5}&{0}&{5}&{35}\cr
{3}&{0}&{-2}&{1}&{35}\cr{2}&{5}&{1}&{-2}&{0}\cr
{0}&{35}&{35}&{0}&{-70}\cr\endpmatrix$.}

\vbox{\noindent
$N=38\ $ $d=35=5\cdot 7,\,\eta=2,\,h=1$:
$\langle 1\rangle\oplus \langle -7 \rangle \oplus \langle -5 \rangle$
\nobreak
\newline
$\pmatrix{0}&{0}&{1}\cr{0}&{1}&{0}\cr
{2}&{0}&{-1}\cr{5}&{-1}&{-2}\cr{91}&{-25}&{-28}\cr
{140}&{-40}&{-41}\cr{196}&{-57}&{-56}\cr
{54}&{-16}&{-15}\cr{9}&{-3}&{-2}\cr{7}&{-3}&{0}\cr\endpmatrix
\pmatrix{-5}&{0}&{5}&{10}&{140}&{205}&{280}&{75}&{10}&{0}\cr
{0}&{-7}&{0}&{7}&{175}&{280}&{399}&{112}&{21}&{21}\cr
{5}&{0}&{-1}&{0}&{42}&{75}&{112}&{33}&{8}&{14}\cr
{10}&{7}&{0}&{-2}&{0}&{10}&{21}&{8}&{4}&{14}\cr
{140}&{175}&{42}&{0}&{-14}&{0}&{21}&{14}&{14}&{112}\cr
{205}&{280}&{75}&{10}&{0}&{-5}&{0}&{5}&{10}&{140}\cr
{280}&{399}&{112}&{21}&{21}&{0}&{-7}&{0}&{7}&{175}\cr
{75}&{112}&{33}&{8}&{14}&{5}&{0}&{-1}&{0}&{42}\cr
{10}&{21}&{8}&{4}&{14}&{10}&{7}&{0}&{-2}&{0}\cr
{0}&{21}&{14}&{14}&{112}&{140}&{175}&{42}&{0}&{-14}\cr
\endpmatrix$.}

\vbox{\noindent
$N=39\ $ $d=35=5\cdot 7,\,\eta=3,\,h=0$:
$\langle 35 \rangle \oplus \langle -1 \rangle \oplus \langle -1 \rangle$
\nobreak
\newline
$\pmatrix{0}&{-1}&{1}\cr{0}&{1}&{0}\cr{1}&{0}&{-7}\cr
{3}&{-10}&{-15}\cr{5}&{-21}&{-21}\cr\endpmatrix
\hskip20pt
\pmatrix{-2}&{1}&{7}&{5}&{0}\cr{1}&{-1}&{0}&{10}&{21}\cr
{7}&{0}&{-14}&{0}&{28}\cr{5}&{10}&{0}&{-10}&{0}\cr
{0}&{21}&{28}&{0}&{-7}\cr\endpmatrix$.}

\vbox{\noindent
$N=40\ $ $d=38=2\cdot 19,\,\eta=0,\,h=1$: $U\oplus \langle -38 \rangle$
\nobreak
\newline
$\pmatrix{19}&{0}&{-1}\cr{0}&{0}&{1}\cr{-1}&{1}&{0}\cr{31}&{30}&{-7}\cr
{171}&{152}&{-37}\cr{190}&{152}&{-39}\cr{6}&{3}&{-1}\cr{9}&{2}&{-1}\cr
\endpmatrix
\hskip20pt
\pmatrix{-38}&{38}&{19}&{304}&{1482}&{1406}&{19}&{0}\cr
{38}&{-38}&{0}&{266}&{1406}&{1482}&{38}&{38}\cr
{19}&{0}&{-2}&{1}&{19}&{38}&{3}&{7}\cr
{304}&{266}&{1}&{-2}&{0}&{38}&{7}&{66}\cr
{1482}&{1406}&{19}&{0}&{-38}&{38}&{19}&{304}\cr
{1406}&{1482}&{38}&{38}&{38}&{-38}&{0}&{266}\cr
{19}&{38}&{3}&{7}&{19}&{0}&{-2}&{1}\cr
{0}&{38}&{7}&{66}&{304}&{266}&{1}&{-2}\cr\endpmatrix$.}

\vbox{\noindent
$N=41\ $ $d=39=3\cdot 13,\,\eta=0,\,h=0$:
$\langle 3 \rangle \oplus \langle -13 \rangle \oplus \langle -1 \rangle$
\nobreak
\newline
$\pmatrix{0}&{0}&{1}\cr{0}&{1}&{0}\cr
{1}&{0}&{-3}\cr{2}&{-1}&{-1}\cr\endpmatrix
\hskip20pt
\pmatrix{-1}&{0}&{3}&{1}\cr{0}&{-13}&{0}&{13}\cr
{3}&{0}&{-6}&{3}\cr{1}&{13}&{3}&{-2}\cr\endpmatrix$.}

\vbox{\noindent
$N=42\ $ $d=39=3\cdot 13,\,\eta=2,\,h=1$:
$\langle 3 \rangle \oplus \langle -26 \rangle
\oplus \langle-2 \rangle (0,1/2,1/2)$
\nobreak
\newline
$\left(\smallmatrix{0}&{0}&{1}\cr{0}&{1}&{0}\cr{2}&{0}&{-3}\cr
{5}&{-1}&{-5}\cr{130}&{-30}&{-117}\cr{12}&{-3}&{-10}\cr
{52}&{-14}&{-39}\cr{10}&{-3}&{-6}\cr{3}&{-1}&{-1}\cr
{26}&{-9}&{0}\cr\endsmallmatrix\right)
\hskip20pt
\left(\smallmatrix{-2}&{0}&{6}&{10}&{234}&{20}&{78}&{12}&{2}&{0}\cr
{0}&{-26}&{0}&{26}&{780}&{78}&{364}&{78}&{26}&{234}\cr
{6}&{0}&{-6}&{0}&{78}&{12}&{78}&{24}&{12}&{156}\cr
{10}&{26}&{0}&{-1}&{0}&{2}&{26}&{12}&{9}&{156}\cr
{234}&{780}&{78}&{0}&{-78}&{0}&{234}&{156}&{156}&{3120}\cr
{20}&{78}&{12}&{2}&{0}&{-2}&{0}&{6}&{10}&{234}\cr
{78}&{364}&{78}&{26}&{234}&{0}&{-26}&{0}&{26}&{780}\cr
{12}&{78}&{24}&{12}&{156}&{6}&{0}&{-6}&{0}&{78}\cr
{2}&{26}&{12}&{9}&{156}&{10}&{26}&{0}&{-1}&{0}\cr
{0}&{234}&{156}&{156}&{3120}&{234}&{780}&{78}&{0}&{-78}\cr
\endsmallmatrix\right)$.}

\vbox{\noindent
$N=43\ $ $d=39=3\cdot 13,\,\eta=3,\, h=1$:
$\langle 2 \rangle\oplus \langle -26\rangle \oplus \langle -3 \rangle
(1/2,1/2,0)$
\nobreak
\newline
$\pmatrix{0}&{0}&{1}\cr{0}&{1}&{0}\cr{1}&{0}&{-1}\cr
{6}&{-1}&{-4}\cr{30}&{-6}&{-17}\cr{{195}\over{2}}&{{-41}\over{2}}&{-52}\cr
{{13}\over{2}}&{{-3}\over{2}}&{-3}\cr{{3}\over{2}}&{{-1}\over{2}}&{0}\cr
\endpmatrix
\hskip20pt
\pmatrix{-3}&{0}&{3}&{12}&{51}&{156}&{9}&{0}\cr
{0}&{-26}&{0}&{26}&{156}&{533}&{39}&{13}\cr
{3}&{0}&{-1}&{0}&{9}&{39}&{4}&{3}\cr
{12}&{26}&{0}&{-2}&{0}&{13}&{3}&{5}\cr
{51}&{156}&{9}&{0}&{-3}&{0}&{3}&{12}\cr
{156}&{533}&{39}&{13}&{0}&{-26}&{0}&{26}\cr
{9}&{39}&{4}&{3}&{3}&{0}&{-1}&{0}\cr
{0}&{13}&{3}&{5}&{12}&{26}&{0}&{-2}\cr\endpmatrix$.}

\vbox{\noindent
$N=44\ $ $d=42=2\cdot 3\cdot 7,\,\eta=0,\,h=0$: $U\oplus \langle -42 \rangle$
\nobreak
\newline
$\pmatrix{21}&{0}&{-1}\cr{0}&{0}&{1}\cr{-1}&{1}&{0}\cr{6}&{3}&{-1}\cr
{10}&{2}&{-1}\cr\endpmatrix
\hskip20pt
\pmatrix{-42}&{42}&{21}&{21}&{0}\cr{42}&{-42}&{0}&{42}&{42}\cr
{21}&{0}&{-2}&{3}&{8}\cr{21}&{42}&{3}&{-6}&{0}\cr
{0}&{42}&{8}&{0}&{-2}\cr\endpmatrix$.}

\vbox{\noindent
$N=45\ $ $d=42=2\cdot 3\cdot 7,\,\eta=3,\,h=1$:
$\langle 6 \rangle\oplus \langle -14\rangle \oplus \langle -2 \rangle
(1/2,1/2,0)$
\nobreak
\newline
$\pmatrix{0}&{0}&{1}\cr{0}&{1}&{0}\cr{1}&{0}&{-2}\cr
{6}&{-2}&{-9}\cr{{21}\over{2}}&{{-9}\over{2}}&{-14}\cr
{{1}\over{2}}&{{-1}\over{2}}&{0}\cr\endpmatrix
\hskip20pt
\pmatrix{-2}&{0}&{4}&{18}&{28}&{0}\cr
{0}&{-14}&{0}&{28}&{63}&{7}\cr{4}&{0}&{-2}&{0}&{7}&{3}\cr
{18}&{28}&{0}&{-2}&{0}&{4}\cr{28}&{63}&{7}&{0}&{-14}&{0}\cr
{0}&{7}&{3}&{4}&{0}&{-2}\cr\endpmatrix$.}

\vbox{\noindent
$N=46\ $ $d=51=3\cdot 17,\,\eta=0,\,h=0$:
$\langle 3 \rangle \oplus \langle -17 \rangle \oplus \langle -1 \rangle$
\nobreak
\newline
$\pmatrix{0}&{0}&{1}\cr{0}&{1}&{0}\cr{1}&{0}&{-3}\cr
{17}&{-6}&{-17}\cr{5}&{-2}&{-3}\cr{7}&{-3}&{0}\cr\endpmatrix
\hskip20pt
\pmatrix{-1}&{0}&{3}&{17}&{3}&{0}\cr
{0}&{-17}&{0}&{102}&{34}&{51}\cr{3}&{0}&{-6}&{0}&{6}&{21}\cr
{17}&{102}&{0}&{-34}&{0}&{51}\cr{3}&{34}&{6}&{0}&{-2}&{3}\cr
{0}&{51}&{21}&{51}&{3}&{-6}\cr\endpmatrix$.}

\vbox{\noindent
$N=47\ $ $d=51=3\cdot 17,\,\eta=1,\,h=1$:
$\langle 17 \rangle \oplus \langle -3\rangle \oplus \langle -1 \rangle$
\nobreak
\newline
$\left(\smallmatrix{0}&{0}&{1}\cr{0}&{1}&{0}\cr{4}&{0}&{-17}\cr
{1}&{-1}&{-4}\cr{3}&{-5}&{-9}\cr{27}&{-51}&{-68}\cr
{4}&{-8}&{-9}\cr{12}&{-25}&{-24}\cr{64}&{-136}&{-119}\cr
{5}&{-11}&{-8}\cr{3}&{-7}&{-3}\cr{7}&{-17}&{0}\cr
\endsmallmatrix\right)
\hskip20pt
\left(\smallmatrix{-1}&{0}&{17}&{4}&{9}&{68}&{9}&
{24}&{119}&{8}&{3}&{0}\cr
{0}&{-3}&{0}&{3}&{15}&{153}&{24}&{75}&{408}&{33}&{21}&{51}\cr
{17}&{0}&{-17}&{0}&{51}&{680}&{119}&{408}&{2329}&{204}&{153}&{476}\cr
{4}&{3}&{0}&{-2}&{0}&{34}&{8}&{33}&{204}&{20}&{18}&{68}\cr
{9}&{15}&{51}&{0}&{-3}&{0}&{3}&{21}&{153}&{18}&{21}&{102}\cr
{68}&{153}&{680}&{34}&{0}&{-34}&{0}&{51}&{476}&{68}&{102}&{612}\cr
{9}&{24}&{119}&{8}&{3}&{0}&{-1}&{0}&{17}&{4}&{9}&{68}\cr
{24}&{75}&{408}&{33}&{21}&{51}&{0}&{-3}&{0}&{3}&{15}&{153}\cr
{119}&{408}&{2329}&{204}&{153}&{476}&{17}&{0}&{-17}&{0}&{51}&{680}\cr
{8}&{33}&{204}&{20}&{18}&{68}&{4}&{3}&{0}&{-2}&{0}&{34}\cr
{3}&{21}&{153}&{18}&{21}&{102}&{9}&{15}&{51}&{0}&{-3}&{0}\cr
{0}&{51}&{476}&{68}&{102}&{612}&{68}&{153}&{680}&{34}&{0}&{-34}\cr
\endsmallmatrix\right)$.}

\vbox{\noindent
$N=48\ $ $d=51=3\cdot 17,\,\eta=3,\,h=1$:
$\langle 51 \rangle \oplus \langle -1\rangle \oplus \langle -1 \rangle$
\nobreak
\newline
$\pmatrix{0}&{-1}&{1}\cr{0}&{1}&{0}\cr
{7}&{0}&{-51}\cr{1}&{-2}&{-7}\cr
{1}&{-4}&{-6}\cr{10}&{-51}&{-51}\cr\endpmatrix
\hskip20pt
\pmatrix{-2}&{1}&{51}&{5}&{2}&{0}\cr
{1}&{-1}&{0}&{2}&{4}&{51}\cr{51}&{0}&{-102}&{0}&{51}&{969}\cr
{5}&{2}&{0}&{-2}&{1}&{51}\cr{2}&{4}&{51}&{1}&{-1}&{0}\cr
{0}&{51}&{969}&{51}&{0}&{-102}\cr\endpmatrix$.}

\vbox{\noindent
$N=49\ $ $d=55=5\cdot 11,\,\eta=0,\,h=1$:
$\langle 11 \rangle \oplus \langle -5\rangle \oplus \langle -1 \rangle$
\nobreak
\newline
$\pmatrix{0}&{0}&{1}\cr{0}&{1}&{0}\cr{3}&{0}&{-11}\cr
{3}&{-2}&{-9}\cr{2}&{-2}&{-5}\cr{5}&{-6}&{-10}\cr
{8}&{-11}&{-11}\cr{2}&{-3}&{-1}\cr\endpmatrix
\hskip20pt
\pmatrix{-1}&{0}&{11}&{9}&{5}&{10}&{11}&{1}\cr
{0}&{-5}&{0}&{10}&{10}&{30}&{55}&{15}\cr
{11}&{0}&{-22}&{0}&{11}&{55}&{143}&{55}\cr
{9}&{10}&{0}&{-2}&{1}&{15}&{55}&{27}\cr
{5}&{10}&{11}&{1}&{-1}&{0}&{11}&{9}\cr
{10}&{30}&{55}&{15}&{0}&{-5}&{0}&{10}\cr
{11}&{55}&{143}&{55}&{11}&{0}&{-22}&{0}\cr
{1}&{15}&{55}&{27}&{9}&{10}&{0}&{-2}\cr\endpmatrix$.}

\vbox{\noindent
$N=50\ $ $d=55=5\cdot 11,\,\eta=3,\,h=1$:
$\langle 2 \rangle \oplus \langle -11\rangle \oplus \langle -10 \rangle
(1/2,0,1/2)$
\nobreak
\newline
$\pmatrix{66}&{-28}&{-3}\cr{7}&{-3}&{0}\cr{0}&{0}&{1}\cr
{0}&{1}&{0}\cr{{1}\over{2}}&{0}&{{-1}\over{2}}\cr
{{5}\over{2}}&{-1}&{{-1}\over{2}}\cr{{95}\over{2}}&{-20}&{{-7}\over{2}}\cr
{418}&{-177}&{-22}\cr\endpmatrix
\hskip20pt
\pmatrix{-2}&{0}&{30}&{308}&{51}&{7}&{5}&{0}\cr
{0}&{-1}&{0}&{33}&{7}&{2}&{5}&{11}\cr
{30}&{0}&{-10}&{0}&{5}&{5}&{35}&{220}\cr
{308}&{33}&{0}&{-11}&{0}&{11}&{220}&{1947}\cr
{51}&{7}&{5}&{0}&{-2}&{0}&{30}&{308}\cr
{7}&{2}&{5}&{11}&{0}&{-1}&{0}&{33}\cr
{5}&{5}&{35}&{220}&{30}&{0}&{-10}&{0}\cr
{0}&{11}&{220}&{1947}&{308}&{33}&{0}&{-11}\cr\endpmatrix$.}

\vbox{\noindent
$N=51\ $ $d=57=3\cdot 19,\,\eta=0,\,h=1$:
$\langle 1 \rangle \oplus \langle -38 \rangle \oplus \langle -6 \rangle
(0,1/2,1/2)$
\nobreak
\newline
$\pmatrix{0}&{0}&{1}\cr{0}&{1}&{0}\cr{2}&{0}&{-1}\cr
{45}&{{-9}\over{2}}&{{-29}\over{2}}\cr
{95}&{{-21}\over{2}}&{{-57}\over{2}}\cr
{3}&{{-1}\over{2}}&{{-1}\over{2}}\cr\endpmatrix
\hskip20pt
\pmatrix{-6}&{0}&{6}&{87}&{171}&{3}\cr{0}&{-38}&{0}&{171}&{399}&{19}\cr
{6}&{0}&{-2}&{3}&{19}&{3}\cr{87}&{171}&{3}&{-6}&{0}&{6}\cr
{171}&{399}&{19}&{0}&{-38}&{0}\cr{3}&{19}&{3}&{6}&{0}&{-2}\cr
\endpmatrix$.}

\vbox{\noindent
$N=52\ $ $d=57=3\cdot 19,\,\eta=2,\,h=1$:
$\langle 3 \rangle \oplus \langle -19 \rangle \oplus \langle -1 \rangle$
\nobreak
\newline
$\left(\smallmatrix{223}&{-84}&{-123}\cr{29}&{-11}&{-15}\cr
{361}&{-138}&{-171}\cr{5}&{-2}&{-1}\cr{0}&{0}&{1}\cr
{0}&{1}&{0}\cr{1}&{0}&{-3}\cr{3}&{-1}&{-3}\cr
{247}&{-90}&{-171}\cr{27}&{-10}&{-17}\cr
{96}&{-36}&{-55}\cr{1216}&{-457}&{-684}\cr\endsmallmatrix\right)
\hskip20pt
\left(\smallmatrix{-6}&{0}&{228}&{30}&{123}&{1596}
&{300}&{42}&{570}&{12}&{3}&{0}\cr
{0}&{-1}&{0}&{2}&{15}&{209}&{42}&{7}&{114}&{4}&{3}&{19}\cr
{228}&{0}&{-114}&{0}&{171}&{2622}&{570}&{114}&{2280}&{114}&{171}&{1710}\cr
{30}&{2}&{0}&{-2}&{1}&{38}&{12}&{4}&{114}&{8}&{17}&{190}\cr
{123}&{15}&{171}&{1}&{-1}&{0}&{3}&{3}&{171}&{17}&{55}&{684}\cr
{1596}&{209}&{2622}&{38}&{0}&{-19}&{0}&{19}&{1710}&{190}&{684}&{8683}\cr
{300}&{42}&{570}&{12}&{3}&{0}&{-6}&{0}&{228}&{30}&{123}&{1596}\cr
{42}&{7}&{114}&{4}&{3}&{19}&{0}&{-1}&{0}&{2}&{15}&{209}\cr
{570}&{114}&{2280}&{114}&{171}&{1710}&{228}&{0}&{-114}&{0}&{171}&{2622}\cr
{12}&{4}&{114}&{8}&{17}&{190}&{30}&{2}&{0}&{-2}&{1}&{38}\cr
{3}&{3}&{171}&{17}&{55}&{684}&{123}&{15}&{171}&{1}&{-1}&{0}\cr
{0}&{19}&{1710}&{190}&{684}&{8683}&{1596}&{209}&{2622}&{38}&{0}&{-19}\cr
\endsmallmatrix\right)$.}

\vbox{\noindent
$N=53\ $ $d=65=5\cdot 13,\,\eta =2,\,h=1$:
$\langle 5 \rangle \oplus \langle -26 \rangle \oplus \langle -2 \rangle
(0,1/2,1/2)$
\nobreak
\newline
$\pmatrix
{0}&{0}&{1}\cr{0}&{1}&{0}\cr{3}&{0}&{-5}\cr
{3}&{{-1}\over{2}}&{{-9}\over{2}}\cr{20}&{-5}&{-26}\cr
{52}&{-14}&{-65}\cr{17}&{-5}&{-20}\cr
{1}&{{-1}\over{2}}&{{-1}\over{2}}\cr\endpmatrix
\hskip20pt
\pmatrix
{-2}&{0}&{10}&{9}&{52}&{130}&{40}&{1}\cr
{0}&{-26}&{0}&{13}&{130}&{364}&{130}&{13}\cr
{10}&{0}&{-5}&{0}&{40}&{130}&{55}&{10}\cr
{9}&{13}&{0}&{-2}&{1}&{13}&{10}&{4}\cr
{52}&{130}&{40}&{1}&{-2}&{0}&{10}&{9}\cr
{130}&{364}&{130}&{13}&{0}&{-26}&{0}&{13}\cr
{40}&{130}&{55}&{10}&{10}&{0}&{-5}&{0}\cr
{1}&{13}&{10}&{4}&{9}&{13}&{0}&{-2}\cr
\endpmatrix$.}

\vbox{\noindent
$N=54\ $ $d=66=2\cdot 3\cdot 11,\,\eta=1,\,h=1$:
$U\oplus \langle -66 \rangle$
\nobreak
\newline
$\left(\smallmatrix{33}&{0}&{-1}\cr{0}&{0}&{1}\cr{-1}&{1}&{0}\cr
{11}&{11}&{-2}\cr{8}&{4}&{-1}\cr{99}&{33}&{-10}\cr
{264}&{66}&{-23}\cr{37}&{8}&{-3}\cr{121}&{22}&{-9}\cr{16}&{2}&{-1}\cr
\endsmallmatrix\right)
\hskip20pt
\left(\smallmatrix{-66}&{66}&{33}&{231}&{66}&{429}&{660}&{66}&{132}&{0}\cr
{66}&{-66}&{0}&{132}&{66}&{660}&{1518}&{198}&{594}&{66}\cr
{33}&{0}&{-2}&{0}&{4}&{66}&{198}&{29}&{99}&{14}\cr
{231}&{132}&{0}&{-22}&{0}&{132}&{594}&{99}&{385}&{66}\cr
{66}&{66}&{4}&{0}&{-2}&{0}&{66}&{14}&{66}&{14}\cr
{429}&{660}&{66}&{132}&{0}&{-66}&{66}&{33}&{231}&{66}\cr
{660}&{1518}&{198}&{594}&{66}&{66}&{-66}&{0}&{132}&{66}\cr
{66}&{198}&{29}&{99}&{14}&{33}&{0}&{-2}&{0}&{4}\cr
{132}&{594}&{99}&{385}&{66}&{231}&{132}&{0}&{-22}&{0}\cr
{0}&{66}&{14}&{66}&{14}&{66}&{66}&{4}&{0}&{-2}\cr
\endsmallmatrix\right)$.}

\vbox{\noindent
$N=55\ $ $d=66=2\cdot 3\cdot 11,\,\eta=2,\,h=1$:
$\langle 22 \rangle \oplus \langle -6 \rangle \oplus \langle -2 \rangle
(0,1/2,1/2)$
\nobreak
\newline
$\pmatrix{0}&{{-1}\over{2}}&{{1}\over{2}}\cr{0}&{1}&{0}\cr
{3}&{0}&{-10}\cr{10}&{-2}&{-33}\cr
{6}&{{-5}\over{2}}&{{-39}\over{2}}\cr
{1}&{-1}&{-3}\cr\endpmatrix
\hskip20pt
\pmatrix{-2}&{3}&{10}&{27}&{12}&{0}\cr
{3}&{-6}&{0}&{12}&{15}&{6}\cr{10}&{0}&{-2}&{0}&{6}&{6}\cr
{27}&{12}&{0}&{-2}&{3}&{10}\cr{12}&{15}&{6}&{3}&{-6}&{0}\cr
{0}&{6}&{6}&{10}&{0}&{-2}\cr\endpmatrix$.}

\vbox{\noindent
$N=56\ $ $d=69=3\cdot 23,\,\eta=3,\,h=0$:
$\langle 1 \rangle \oplus \langle -23 \rangle \oplus \langle -3 \rangle$
\nobreak
\newline
$\pmatrix{0}&{0}&{1}\cr{0}&{1}&{0}\cr{1}&{0}&{-1}\cr
{69}&{-12}&{-23}\cr{5}&{-1}&{-1}\cr{23}&{-5}&{0}\cr\endpmatrix
\hskip20pt
\pmatrix
{-3}&{0}&{3}&{69}&{3}&{0}\cr
{0}&{-23}&{0}&{276}&{23}&{115}\cr
{3}&{0}&{-2}&{0}&{2}&{23}\cr{69}&{276}&{0}&{-138}&{0}&{207}\cr
{3}&{23}&{2}&{0}&{-1}&{0}\cr{0}&{115}&{23}&{207}&{0}&{-46}\cr
\endpmatrix$.}

\vbox{\noindent
$N=57\ $ $d=70=2\cdot 5\cdot 7,\,\eta=3,\,h=1$:
$\langle 2 \rangle \oplus \langle -14\rangle \oplus \langle -10\rangle
(0,1/2,1/2)$
\nobreak
\newline
$\pmatrix{0}&{0}&{1}\cr{0}&{1}&{0}\cr
{2}&{0}&{-1}\cr{35}&{-6}&{-14}\cr
{10}&{{-5}\over{2}}&{{-7}\over{2}}\cr
{14}&{{-9}\over{2}}&{{-7}\over{2}}\cr
{4}&{{-3}\over{2}}&{{-1}\over{2}}\cr{21}&{-8}&{0}\cr\endpmatrix
\hskip20pt
\pmatrix{-10}&{0}&{10}&{140}&{35}&{35}&{5}&{0}\cr
{0}&{-14}&{0}&{84}&{35}&{63}&{21}&{112}\cr
{10}&{0}&{-2}&{0}&{5}&{21}&{11}&{84}\cr
{140}&{84}&{0}&{-14}&{0}&{112}&{84}&{798}\cr
{35}&{35}&{5}&{0}&{-10}&{0}&{10}&{140}\cr
{35}&{63}&{21}&{112}&{0}&{-14}&{0}&{84}\cr
{5}&{21}&{11}&{84}&{10}&{0}&{-2}&{0}\cr
{0}&{112}&{84}&{798}&{140}&{84}&{0}&{-14}\cr\endpmatrix$.}

\vbox{\noindent
$N=58\ $ $d=77=7\cdot 11,\,\eta=3,\,h=0$:
$\langle 1 \rangle \oplus \langle -11\rangle \oplus \langle -7\rangle$
\nobreak
\newline
$\pmatrix{0}&{0}&{1}\cr{0}&{1}&{0}\cr{7}&{0}&{-3}\cr
{21}&{-3}&{-7}\cr{308}&{-49}&{-99}\cr
{4}&{-1}&{-1}\cr{3}&{-1}&{0}\cr\endpmatrix
\hskip20pt
\pmatrix{-7}&{0}&{21}&{49}&{693}&{7}&{0}\cr
{0}&{-11}&{0}&{33}&{539}&{11}&{11}\cr
{21}&{0}&{-14}&{0}&{77}&{7}&{21}\cr
{49}&{33}&{0}&{-1}&{0}&{2}&{30}\cr
{693}&{539}&{77}&{0}&{-154}&{0}&{385}\cr
{7}&{11}&{7}&{2}&{0}&{-2}&{1}\cr
{0}&{11}&{21}&{30}&{385}&{1}&{-2}\cr\endpmatrix$.}

\vbox{\noindent
$N=59\ $ $d=78=2\cdot 3\cdot 13,\,\eta=2,h=0$: $U\oplus \langle -78 \rangle$
\nobreak
\newline
$\pmatrix{39}&{0}&{-1}\cr{0}&{0}&{1}\cr{-1}&{1}&{0}\cr
{6}&{6}&{-1}\cr{26}&{13}&{-3}\cr{12}&{3}&{-1}\cr{19}&{2}&{-1}\cr
\endpmatrix
\hskip20pt
\pmatrix{-78}&{78}&{39}&{156}&{273}&{39}&{0}\cr
{78}&{-78}&{0}&{78}&{234}&{78}&{78}\cr
{39}&{0}&{-2}&{0}&{13}&{9}&{17}\cr{156}&{78}&{0}&{-6}&{0}&{12}&{48}\cr
{273}&{234}&{13}&{0}&{-26}&{0}&{65}\cr{39}&{78}&{9}&{12}&{0}&{-6}&{3}\cr
{0}&{78}&{17}&{48}&{65}&{3}&{-2}\cr\endpmatrix$.}

\vbox{\noindent
$N=60\ $ $d=85=5\cdot 17,\,\eta=1,\,h=1$:
$\langle 1 \rangle \oplus \langle -34 \rangle \oplus \langle -10 \rangle
(0,1/2,1/2)$
\nobreak
\newline
$\left(\smallmatrix{476}&{-35}&{-136}\cr{25}&{-2}&{-7}\cr
{25}&{{-5}\over{2}}&{{-13}\over{2}}\cr
{3}&{{-1}\over{2}}&{{-1}\over{2}}\cr
{17}&{-3}&{0}\cr{0}&{0}&{1}\cr{0}&{1}&{0}\cr
{3}&{0}&{-1}\cr{45}&{{-5}\over{2}}&{{-27}\over{2}}\cr
{67}&{{-9}\over{2}}&{{-39}\over{2}}\cr{1887}&{-133}&{-544}\cr
{560}&{-40}&{-161}\cr\endsmallmatrix\right)
\hskip20pt
\left(\smallmatrix
{-34}&{0}&{85}&{153}&{4522}&{1360}&{1190}&{68}&{85}&{17}&{102}&{0}\cr
{0}&{-1}&{0}&{6}&{221}&{70}&{68}&{5}&{10}&{4}&{51}&{10}\cr
{85}&{0}&{-10}&{0}&{170}&{65}&{85}&{10}&{35}&{25}&{510}&{135}\cr
{153}&{6}&{0}&{-2}&{0}&{5}&{17}&{4}&{25}&{27}&{680}&{195}\cr
{4522}&{221}&{170}&{0}&{-17}&{0}&{102}&{51}&{510}&{680}&{18513}&{5440}\cr
{1360}&{70}&{65}&{5}&{0}&{-10}&{0}&{10}&{135}&{195}&{5440}&{1610}\cr
{1190}&{68}&{85}&{17}&{102}&{0}&{-34}&{0}&{85}&{153}&{4522}&{1360}\cr
{68}&{5}&{10}&{4}&{51}&{10}&{0}&{-1}&{0}&{6}&{221}&{70}\cr
{85}&{10}&{35}&{25}&{510}&{135}&{85}&{0}&{-10}&{0}&{170}&{65}\cr
{17}&{4}&{25}&{27}&{680}&{195}&{153}&{6}&{0}&{-2}&{0}&{5}\cr
{102}&{51}&{510}&{680}&{18513}&{5440}&{4522}&{221}&{170}&{0}&{-17}&{0}\cr
{0}&{10}&{135}&{195}&{5440}&{1610}&{1360}&{70}&{65}&{5}&{0}&{-10}\cr
\endsmallmatrix\right)$.}

\vbox{\noindent
$N=61\ $ $d=87=3\cdot 29,\,\eta=3,\,h=1$:
$\langle 2 \rangle \oplus \langle -58 \rangle \oplus \langle -3 \rangle
(1/2,1/2,0)$
\nobreak
\newline
$\pmatrix{0}&{0}&{1}\cr{0}&{1}&{0}\cr{1}&{0}&{-1}\cr
{{7}\over{2}}&{{-1}\over{2}}&{-2}\cr{36}&{-6}&{-13}\cr
{348}&{-59}&{-116}\cr{35}&{-6}&{-11}\cr
{{5}\over{2}}&{{-1}\over{2}}&{0}\cr\endpmatrix
\hskip20pt
\pmatrix{-3}&{0}&{3}&{6}&{39}&{348}&{33}&{0}\cr
{0}&{-58}&{0}&{29}&{348}&{3422}&{348}&{29}\cr
{3}&{0}&{-1}&{1}&{33}&{348}&{37}&{5}\cr
{6}&{29}&{1}&{-2}&{0}&{29}&{5}&{3}\cr
{39}&{348}&{33}&{0}&{-3}&{0}&{3}&{6}\cr
{348}&{3422}&{348}&{29}&{0}&{-58}&{0}&{29}\cr
{33}&{348}&{37}&{5}&{3}&{0}&{-1}&{1}\cr
{0}&{29}&{5}&{3}&{6}&{29}&{1}&{-2}\cr\endpmatrix$.}

\vbox{\noindent
$N=62\ $ $d=91=7\cdot 13,\,\eta=3,\,h=1$:
$\langle 1 \rangle \oplus \langle -26 \rangle \oplus \langle -14 \rangle
(0,1/2,1/2)$
\nobreak
\newline
$\left(\smallmatrix{138}&{{-53}\over{2}}&{{-15}\over{2}}\cr
{5}&{-1}&{0}\cr{0}&{0}&{1}\cr{0}&{1}&{0}\cr
{7}&{0}&{-2}\cr{6}&{{-1}\over{2}}&{{-3}\over{2}}\cr
{3}&{{-1}\over{2}}&{{-1}\over{2}}\cr
{56}&{{-21}\over{2}}&{{-9}\over{2}}\cr
{312}&{{-119}\over{2}}&{{-39}\over{2}}\cr
{329}&{-63}&{-19}\cr\endsmallmatrix\right)
\hskip20pt 
\left(\smallmatrix{-2}&{1}&{105}&{689}&{756}&{326}&{17}&{21}&{13}&{0}\cr
{1}&{-1}&{0}&{26}&{35}&{17}&{2}&{7}&{13}&{7}\cr
{105}&{0}&{-14}&{0}&{28}&{21}&{7}&{63}&{273}&{266}\cr
{689}&{26}&{0}&{-26}&{0}&{13}&{13}&{273}&{1547}&{1638}\cr
{756}&{35}&{28}&{0}&{-7}&{0}&{7}&{266}&{1638}&{1771}\cr
{326}&{17}&{21}&{13}&{0}&{-2}&{1}&{105}&{689}&{756}\cr
{17}&{2}&{7}&{13}&{7}&{1}&{-1}&{0}&{26}&{35}\cr
{21}&{7}&{63}&{273}&{266}&{105}&{0}&{-14}&{0}&{28}\cr
{13}&{13}&{273}&{1547}&{1638}&{689}&{26}&{0}&{-26}&{0}\cr
{0}&{7}&{266}&{1638}&{1771}&{756}&{35}&{28}&{0}&{-7}\cr
\endsmallmatrix\right)$.}

\vbox{\noindent
$N=63\ $ $d=93=3\cdot 31,\,\eta=0,\,h=1$:
$\langle 31 \rangle \oplus \langle -6 \rangle \oplus \langle -2 \rangle
(0,1/2,1/2)$
\nobreak
\newline
$\pmatrix{0}&{{-1}\over{2}}&{{1}\over{2}}\cr{0}&{1}&{0}\cr
{1}&{0}&{-4}\cr{36}&{{-31}\over{2}}&{{-279}\over{2}}\cr
{2}&{{-3}\over{2}}&{{-15}\over{2}}\cr{24}&{-25}&{-84}\cr
{11}&{-12}&{-38}\cr{150}&{{-341}\over{2}}&{{-1023}\over{2}}\cr\endpmatrix
\hskip20pt
\pmatrix{-2}&{3}&{4}&{93}&{3}&{9}&{2}&{0}\cr
{3}&{-6}&{0}&{93}&{9}&{150}&{72}&{1023}\cr
{4}&{0}&{-1}&{0}&{2}&{72}&{37}&{558}\cr
{93}&{93}&{0}&{-186}&{0}&{1023}&{558}&{8835}\cr
{3}&{9}&{2}&{0}&{-2}&{3}&{4}&{93}\cr
{9}&{150}&{72}&{1023}&{3}&{-6}&{0}&{93}\cr
{2}&{72}&{37}&{558}&{4}&{0}&{-1}&{0}\cr
{0}&{1023}&{558}&{8835}&{93}&{93}&{0}&{-186}\cr\endpmatrix$.}

\vbox{\noindent
$N=64\ $ $d=95=5\cdot 19,\,\eta=3,\,h=1$:
$\langle 2 \rangle \oplus \langle -19 \rangle \oplus \langle-10\rangle
(1/2,0,1/2)$
\nobreak
\newline
$\left(\smallmatrix{{129}\over{2}}&{-20}&{{-17}\over{2}}\cr
{19}&{-6}&{-2}\cr{3}&{-1}&{0}\cr{0}&{0}&{1}\cr
{0}&{1}&{0}\cr{{1}\over{2}}&{0}&{{-1}\over{2}}\cr
{{27}\over{2}}&{-4}&{{-5}\over{2}}\cr
{{59}\over{2}}&{-9}&{{-9}\over{2}}\cr
{{585}\over{2}}&{-90}&{{-83}\over{2}}\cr
{1235}&{-381}&{-171}\cr\endsmallmatrix\right)
\hskip20pt
\left(\smallmatrix{-2}&{1}&{7}&{85}&{380}&{22}&{9}&{3}&{5}&{0}\cr
{1}&{-2}&{0}&{20}&{114}&{9}&{7}&{5}&{25}&{76}\cr
{7}&{0}&{-1}&{0}&{19}&{3}&{5}&{6}&{45}&{171}\cr
{85}&{20}&{0}&{-10}&{0}&{5}&{25}&{45}&{415}&{1710}\cr
{380}&{114}&{19}&{0}&{-19}&{0}&{76}&{171}&{1710}&{7239}\cr
{22}&{9}&{3}&{5}&{0}&{-2}&{1}&{7}&{85}&{380}\cr
{9}&{7}&{5}&{25}&{76}&{1}&{-2}&{0}&{20}&{114}\cr
{3}&{5}&{6}&{45}&{171}&{7}&{0}&{-1}&{0}&{19}\cr
{5}&{25}&{45}&{415}&{1710}&{85}&{20}&{0}&{-10}&{0}\cr
{0}&{76}&{171}&{1710}&{7239}&{380}&{114}&{19}&{0}&{-19}\cr
\endsmallmatrix\right)$.}

\vbox{\noindent
$N=65\ $ $d=102=2\cdot 3\cdot 17,\,\eta=0,\,h=1$:
$\langle 2 \rangle \oplus \langle-34\rangle \oplus \langle-6\rangle
(0,1/2,1/2)$
\nobreak
\newline
$\pmatrix
{0}&{0}&{1}\cr{0}&{1}&{0}\cr{3}&{0}&{-2}\cr{8}&{-1}&{-4}\cr
{270}&{{-81}\over{2}}&{{-245}\over{2}}\cr
{510}&{{-155}\over{2}}&{{-459}\over{2}}\cr{57}&{-9}&{-25}\cr
{2}&{{-1}\over{2}}&{{-1}\over{2}}\cr\endpmatrix
\hskip20pt
\pmatrix
{-6}&{0}&{12}&{24}&{735}&{1377}&{150}&{3}\cr
{0}&{-34}&{0}&{34}&{1377}&{2635}&{306}&{17}\cr
{12}&{0}&{-6}&{0}&{150}&{306}&{42}&{6}\cr
{24}&{34}&{0}&{-2}&{3}&{17}&{6}&{3}\cr
{735}&{1377}&{150}&{3}&{-6}&{0}&{12}&{24}\cr
{1377}&{2635}&{306}&{17}&{0}&{-34}&{0}&{34}\cr
{150}&{306}&{42}&{6}&{12}&{0}&{-6}&{0}\cr
{3}&{17}&{6}&{3}&{24}&{34}&{0}&{-2}\cr
\endpmatrix$.}

\vbox{\noindent
$N=66\ $ $d=105=3\cdot 5\cdot 7,\,\eta=0,h=1$:
$\langle 2 \rangle \oplus \langle-42\rangle \oplus \langle-5\rangle
(1/2,1/2,0)$
\nobreak
\newline
$\left(\smallmatrix{0}&{0}&{1}\cr{0}&{1}&{0}\cr{3}&{0}&{-2}\cr
{70}&{-5}&{-42}\cr{21}&{-2}&{-12}\cr{90}&{-10}&{-49}\cr
{{315}\over{2}}&{{-37}\over{2}}&{-84}\cr{{39}\over{2}}&
{{-5}\over{2}}&{-10}\cr{{175}\over{2}}&{{-25}\over{2}}
&{-42}\cr{{3}\over{2}}&{{-1}\over{2}}&{0}\cr\endsmallmatrix\right)
\hskip20pt
\left(\smallmatrix{-5}&{0}&{10}&{210}&{60}&{245}&{420}&{50}&{210}&{0}\cr
{0}&{-42}&{0}&{210}&{84}&{420}&{777}&{105}&{525}&{21}\cr
{10}&{0}&{-2}&{0}&{6}&{50}&{105}&{17}&{105}&{9}\cr
{210}&{210}&{0}&{-70}&{0}&{210}&{525}&{105}&{805}&{105}\cr
{60}&{84}&{6}&{0}&{-6}&{0}&{21}&{9}&{105}&{21}\cr
{245}&{420}&{50}&{210}&{0}&{-5}&{0}&{10}&{210}&{60}\cr
{420}&{777}&{105}&{525}&{21}&{0}&{-42}&{0}&{210}&{84}\cr
{50}&{105}&{17}&{105}&{9}&{10}&{0}&{-2}&{0}&{6}\cr
{210}&{525}&{105}&{805}&{105}&{210}&{210}&{0}&{-70}&{0}\cr
{0}&{21}&{9}&{105}&{21}&{60}&{84}&{6}&{0}&{-6}\cr\endsmallmatrix
\right)$.}

\vbox{\noindent
$N=67\ $ $d=105=3\cdot 5\cdot 7,\,\eta=1,\,h=0$:
$\langle 1 \rangle \oplus \langle-21\rangle \oplus \langle-5\rangle$
\nobreak
\newline
$\pmatrix{0}&{0}&{1}\cr{0}&{1}&{0}\cr{2}&{0}&{-1}\cr
{8}&{-1}&{-3}\cr{30}&{-5}&{-9}\cr{5}&{-1}&{-1}\cr
{70}&{-15}&{-7}\cr{9}&{-2}&{0}\cr\endpmatrix
\hskip20pt
\pmatrix{-5}&{0}&{5}&{15}&{45}&{5}&{35}&{0}\cr
{0}&{-21}&{0}&{21}&{105}&{21}&{315}&{42}\cr
{5}&{0}&{-1}&{1}&{15}&{5}&{105}&{18}\cr
{15}&{21}&{1}&{-2}&{0}&{4}&{140}&{30}\cr
{45}&{105}&{15}&{0}&{-30}&{0}&{210}&{60}\cr
{5}&{21}&{5}&{4}&{0}&{-1}&{0}&{3}\cr
{35}&{315}&{105}&{140}&{210}&{0}&{-70}&{0}\cr
{0}&{42}&{18}&{30}&{60}&{3}&{0}&{-3}\cr\endpmatrix$.}

\vbox{\noindent
$N=68\ $ $d=105=3\cdot 5 \cdot 7,\,\eta=2,\,h=1$:
$\langle 7 \rangle \oplus \langle-15\rangle \oplus \langle -1\rangle$
\nobreak
\newline
$\left(\smallmatrix
{108}&{360}&{53}&{15}&{57}&{15}&{3}&{10}&{0}&{0}
&{1}&{3}&{69}&{75}&{87}&{1250}\cr
{-48}&{-161}&{-24}&{-7}&{-28}&{-8}&{-2}&{-7}&{0}
&{1}&{0}&{-1}&{-28}&{-32}&{-38}&{-553}\cr
{-217}&{-720}&{-105}&{-29}&{-105}&{-25}&{-3}&{0}
&{1}&{0}&{-3}&{-7}&{-147}&{-155}&{-177}&{-2520}\cr
\endsmallmatrix\right)$
\nobreak
\newline
$\left(\smallmatrix
{-1}&{0}&{3}&{7}&{147}&{155}&{177}&{2520}&{217}&{720}&
{105}&{29}&{105}&{25}&{3}&{0}\cr
{0}&{-15}&{0}&{15}&{420}&{480}&{570}&{8295}&{720}
&{2415}&{360}&{105}&{420}&{120}&{30}&{105}\cr
{3}&{0}&{-2}&{0}&{42}&{60}&{78}&{1190}&{105}&{360}
&{56}&{18}&{84}&{30}&{12}&{70}\cr
{7}&{15}&{0}&{-1}&{0}&{10}&{18}&{315}&{29}&{105}&{18}
&{7}&{42}&{20}&{12}&{105}\cr
{147}&{420}&{42}&{0}&{-42}&{0}&{42}&{1050}&{105}&{420}
&{84}&{42}&{336}&{210}&{168}&{1890}\cr
{155}&{480}&{60}&{10}&{0}&{-10}&{0}&{210}&{25}&{120}&{30}
&{20}&{210}&{160}&{150}&{1890}\cr
{177}&{570}&{78}&{18}&{42}&{0}&{-6}&{0}&{3}&{30}&{12}&{12}
&{168}&{150}&{156}&{2100}\cr
{2520}&{8295}&{1190}&{315}&{1050}&{210}&{0}&{-35}&{0}&{105}
&{70}&{105}&{1890}&{1890}&{2100}&{29435}\cr
{217}&{720}&{105}&{29}&{105}&{25}&{3}&{0}&{-1}&{0}&{3}&{7}
&{147}&{155}&{177}&{2520}\cr
{720}&{2415}&{360}&{105}&{420}&{120}&{30}&{105}&{0}&{-15}&
{0}&{15}&{420}&{480}&{570}&{8295}\cr
{105}&{360}&{56}&{18}&{84}&{30}&{12}&{70}&{3}&{0}&{-2}&{0}
&{42}&{60}&{78}&{1190}\cr
{29}&{105}&{18}&{7}&{42}&{20}&{12}&{105}&{7}&{15}&{0}&{-1}
&{0}&{10}&{18}&{315}\cr
{105}&{420}&{84}&{42}&{336}&{210}&{168}&{1890}&{147}&{420}
&{42}&{0}&{-42}&{0}&{42}&{1050}\cr
{25}&{120}&{30}&{20}&{210}&{160}&{150}&{1890}&{155}&{480}
&{60}&{10}&{0}&{-10}&{0}&{210}\cr
{3}&{30}&{12}&{12}&{168}&{150}&{156}&{2100}&{177}&{570}
&{78}&{18}&{42}&{0}&{-6}&{0}\cr
{0}&{105}&{70}&{105}&{1890}&{1890}&{2100}&{29435}&{2520}
&{8295}&{1190}&{315}&{1050}&{210}&{0}&{-35}\cr
\endsmallmatrix\right)$.}

\vbox{\noindent
$N=69\ $ $d=105=3\cdot 5\cdot 7,\,\eta=3,\,h=0$:
$\langle 14 \rangle \oplus \langle -10 \rangle \oplus \langle -3 \rangle
(1/2,1/2,0)$
\nobreak
\newline
$\pmatrix{0}&{0}&{1}\cr{0}&{1}&{0}\cr{3}&{0}&{-7}\cr
{{1}\over{2}}&{{-1}\over{2}}&{-1}\cr
{{5}\over{2}}&{{-7}\over{2}}&{0}\cr\endpmatrix
\hskip20pt
\pmatrix{-3}&{0}&{21}&{3}&{0}\cr{0}&{-10}&{0}&{5}&{35}\cr
{21}&{0}&{-21}&{0}&{105}\cr{3}&{5}&{0}&{-2}&{0}\cr
{0}&{35}&{105}&{0}&{-35}\cr\endpmatrix$.}

\vbox{\noindent
$N=70\ $ $d=105=3\cdot 5\cdot 7,\,\eta=5,\,h=0$:
$\langle 5 \rangle \oplus \langle -7 \rangle \oplus \langle -3 \rangle$
\nobreak
\newline
$\pmatrix{0}&{0}&{1}\cr{0}&{1}&{0}\cr{3}&{0}&{-5}\cr
{7}&{-4}&{-7}\cr{7}&{-5}&{-5}\cr{1}&{-1}&{0}\cr\endpmatrix
\hskip20pt
\pmatrix{-3}&{0}&{15}&{21}&{15}&{0}\cr
{0}&{-7}&{0}&{28}&{35}&{7}\cr{15}&{0}&{-30}&{0}&{30}&{15}\cr
{21}&{28}&{0}&{-14}&{0}&{7}\cr{15}&{35}&{30}&{0}&{-5}&{0}\cr
{0}&{7}&{15}&{7}&{0}&{-2}\cr\endpmatrix$.}

\vbox{\noindent
$N=71\ $ $d=105=3\cdot 5\cdot 7,\,\eta=6,\,h=0$:
$\langle 1 \rangle \oplus \langle -15 \rangle \oplus \langle -7 \rangle$
\nobreak
\newline
$\pmatrix
{0}&{0}&{1}\cr{0}&{1}&{0}\cr{7}&{0}&{-3}\cr{15}&{-2}&{-5}\cr
{11}&{-2}&{-3}\cr{3}&{-1}&{0}\cr\endpmatrix
\hskip20pt
\pmatrix
{-7}&{0}&{21}&{35}&{21}&{0}\cr
{0}&{-15}&{0}&{30}&{30}&{15}\cr
{21}&{0}&{-14}&{0}&{14}&{21}\cr
{35}&{30}&{0}&{-10}&{0}&{15}\cr
{21}&{30}&{14}&{0}&{-2}&{3}\cr
{0}&{15}&{21}&{15}&{3}&{-6}\cr
\endpmatrix$.}

\vbox{\noindent
$N=72\ $ $d=105=3\cdot 5\cdot 7,\,\eta =7,\,h=0$:
$\langle 10 \rangle \oplus \langle -14 \rangle \oplus \langle -3 \rangle
(1/2,1/2,0)$
\nobreak
\newline
$\pmatrix{0}&{0}&{1}\cr{0}&{1}&{0}\cr{1}&{0}&{-2}\cr
{42}&{-15}&{-70}\cr{2}&{-1}&{-3}\cr{8}&{-5}&{-10}\cr
{7}&{-5}&{-7}\cr{{1}\over{2}}&{{-1}\over{2}}&{0}\cr\endpmatrix
\hskip20pt
\pmatrix{-3}&{0}&{6}&{210}&{9}&{30}&{21}&{0}\cr
{0}&{-14}&{0}&{210}&{14}&{70}&{70}&{7}\cr
{6}&{0}&{-2}&{0}&{2}&{20}&{28}&{5}\cr
{210}&{210}&{0}&{-210}&{0}&{210}&{420}&{105}\cr
{9}&{14}&{2}&{0}&{-1}&{0}&{7}&{3}\cr
{30}&{70}&{20}&{210}&{0}&{-10}&{0}&{5}\cr
{21}&{70}&{28}&{420}&{7}&{0}&{-7}&{0}\cr
{0}&{7}&{5}&{105}&{3}&{5}&{0}&{-1}\cr\endpmatrix$.}

\vbox{\noindent
$N=73\ $ $d=110=2\cdot 5\cdot 11,\,\eta=1,\,h=1$:
$U\oplus \langle -110 \rangle$
\nobreak
\newline
$\left(\smallmatrix
{55}&{0}&{-1}\cr{0}&{0}&{1}\cr{-1}&{1}&{0}\cr{22}&{22}&{-3}\cr
{10}&{5}&{-1}\cr{44}&{11}&{-3}\cr{18}&{3}&{-1}\cr{440}&{55}&{-21}\cr
{1980}&{220}&{-89}\cr{361}&{39}&{-16}\cr{3938}&{418}&{-173}\cr
{485}&{50}&{-21}\cr{451}&{44}&{-19}\cr{27}&{2}&{-1}\cr
\endsmallmatrix\right) \ \
\left(\smallmatrix
{-110}&{110}&{55}&{880}&{165}&{275}&{55}&{715}&{2310}
&{385}&{3960}&{440}&{330}&{0}\cr
{110}&{-110}&{0}&{330}&{110}&{330}&{110}&{2310}&{9790}&
{1760}&{19030}&{2310}&{2090}&{110}\cr
{55}&{0}&{-2}&{0}&{5}&{33}&{15}&{385}&{1760}&{322}&{3520}&{435}&{407}&{25}\cr
{880}&{330}&{0}&{-22}&{0}&{220}&{132}&{3960}&{19030}&
{3520}&{38742}&{4840}&{4620}&{308}\cr
{165}&{110}&{5}&{0}&{-10}&{0}&{10}&{440}&{2310}&{435}&{4840}
&{615}&{605}&{45}\cr
{275}&{330}&{33}&{220}&{0}&{-22}&{0}&{330}&{2090}&{407}
&{4620}&{605}&{627}&{55}\cr
{55}&{110}&{15}&{132}&{10}&{0}&{-2}&{0}&{110}&{25}&{308}&{45}&{55}&{7}\cr
{715}&{2310}&{385}&{3960}&{440}&{330}&{0}&{-110}&
{110}&{55}&{880}&{165}&{275}&{55}\cr
{2310}&{9790}&{1760}&{19030}&{2310}&{2090}&{110}&{110}&{-110}&
{0}&{330}&{110}&{330}&{110}\cr
{385}&{1760}&{322}&{3520}&{435}&{407}&{25}&{55}&{0}&{-2}&{0}&{5}&{33}&{15}\cr
{3960}&{19030}&{3520}&{38742}&{4840}&{4620}&{308}&
{880}&{330}&{0}&{-22}&{0}&{220}&{132}\cr
{440}&{2310}&{435}&{4840}&{615}&{605}&{45}&{165}
&{110}&{5}&{0}&{-10}&{0}&{10}\cr
{330}&{2090}&{407}&{4620}&{605}&{627}&{55}&{275}&{330}&{33}&
{220}&{0}&{-22}&{0}\cr
{0}&{110}&{25}&{308}&{45}&{55}&{7}&{55}&{110}&{15}&{132}&{10}&{0}&{-2}\cr
\endsmallmatrix\right)$.}

\vbox{\noindent
$N=74\ $ $d=111=3\cdot 37,\,\eta=0,\,h=1$:
$\langle 3\rangle \oplus \langle-37 \rangle \oplus \langle-1\rangle$
\nobreak
\newline
$\pmatrix{0}&{0}&{1}\cr{0}&{1}&{0}\cr{1}&{0}&{-3}\cr
{21}&{-5}&{-20}\cr{12}&{-3}&{-10}\cr{148}&{-38}&{-111}\cr
{11}&{-3}&{-6}\cr{7}&{-2}&{-1}\cr\endpmatrix
\hskip20pt
\pmatrix{-1}&{0}&{3}&{20}&{10}&{111}&{6}&{1}\cr
{0}&{-37}&{0}&{185}&{111}&{1406}&{111}&{74}\cr
{3}&{0}&{-6}&{3}&{6}&{111}&{15}&{18}\cr
{20}&{185}&{3}&{-2}&{1}&{74}&{18}&{51}\cr
{10}&{111}&{6}&{1}&{-1}&{0}&{3}&{20}\cr
{111}&{1406}&{111}&{74}&{0}&{-37}&{0}&{185}\cr
{6}&{111}&{15}&{18}&{3}&{0}&{-6}&{3}\cr
{1}&{74}&{18}&{51}&{20}&{185}&{3}&{-2}\cr\endpmatrix$.}

\vbox{\noindent
$N=75\ $ $d=130=2\cdot 5\cdot 13,\,\eta=3,\,h=1$:
$\langle 10 \rangle \oplus \langle-26 \rangle \oplus \langle -2 \rangle
(1/2,1/2,0)$
\nobreak
\newline
$\left(\smallmatrix{0}&{0}&{1}\cr{0}&{1}&{0}\cr{2}&{0}&{-5}\cr
{{3}\over{2}}&{{-1}\over{2}}&{-3}\cr
{26}&{-12}&{-39}\cr{{9}\over{2}}&{{-5}\over{2}}&{-5}\cr
{5}&{-3}&{-3}\cr{{195}\over{2}}&{{-119}\over{2}}&{-39}\cr
{98}&{-60}&{-35}\cr{{57}\over{2}}&{{-35}\over{2}}&{-9}\cr
{234}&{-144}&{-65}\cr{8}&{-5}&{0}\cr\endsmallmatrix\right)
\hskip20pt
\left(\smallmatrix{-2}&{0}&{10}&{6}&{78}&{10}&
{6}&{78}&{70}&{18}&{130}&{0}\cr
{0}&{-26}&{0}&{13}&{312}&{65}&{78}&{1547}&{1560}&{455}&{3744}&{130}\cr
{10}&{0}&{-10}&{0}&{130}&{40}&{70}&{1560}&{1610}&{480}&{4030}&{160}\cr
{6}&{13}&{0}&{-2}&{0}&{5}&{18}&{455}&{480}&{146}&{1248}&{55}\cr
{78}&{312}&{130}&{0}&{-26}&{0}&{130}&{3744}&{4030}&{1248}&{10842}&{520}\cr
{10}&{65}&{40}&{5}&{0}&{-10}&{0}&{130}&{160}&{55}&{520}&{35}\cr
{6}&{78}&{70}&{18}&{130}&{0}&{-2}&{0}&{10}&{6}&{78}&{10}\cr
{78}&{1547}&{1560}&{455}&{3744}&{130}&{0}&{-26}&{0}&{13}&{312}&{65}\cr
{70}&{1560}&{1610}&{480}&{4030}&{160}&{10}&{0}&{-10}&{0}&{130}&{40}\cr
{18}&{455}&{480}&{146}&{1248}&{55}&{6}&{13}&{0}&{-2}&{0}&{5}\cr
{130}&{3744}&{4030}&{1248}&{10842}&{520}&{78}&{312}&{130}&{0}&{-26}&{0}\cr
{0}&{130}&{160}&{55}&{520}&{35}&{10}&{65}&{40}&{5}&{0}&{-10}\cr
\endsmallmatrix\right)$.}

\vbox{\noindent
$N=76\ $ $d=141=3\cdot 47,\,\eta=3,\,h=1$:
$\langle 1 \rangle \oplus \langle-47 \rangle \oplus \langle -3 \rangle$
\nobreak
\newline
$\left(\smallmatrix{0}&{0}&{1}\cr{0}&{1}&{0}\cr{1}&{0}&{-1}\cr
{705}&{-84}&{-235}\cr{41}&{-5}&{-13}\cr{1081}&{-134}&{-329}\cr
{168}&{-21}&{-50}\cr{1128}&{-142}&{-329}\cr{23}&{-3}&{-6}\cr
{423}&{-57}&{-94}\cr{7}&{-1}&{-1}\cr{47}&{-7}&{0}\cr
\endsmallmatrix\right)
\hskip20pt
\left(\smallmatrix{-3}&{0}&{3}&{705}&{39}&{987}&{150}&
{987}&{18}&{282}&{3}&{0}\cr
{0}&{-47}&{0}&{3948}&{235}&{6298}&{987}&{6674}
&{141}&{2679}&{47}&{329}\cr
{3}&{0}&{-2}&{0}&{2}&{94}&{18}&{141}&{5}&{141}&{4}&{47}\cr
{705}&{3948}&{0}&{-282}&{0}&{1128}&{282}
&{2679}&{141}&{6909}&{282}&{5499}\cr
{39}&{235}&{2}&{0}&{-1}&{0}&{3}&{47}&
{4}&{282}&{13}&{282}\cr
{987}&{6298}&{94}&{1128}&{0}&{-94}&{0}&{329}
&{47}&{5499}&{282}&{6721}\cr
{150}&{987}&{18}&{282}&{3}&{0}&{-3}&{0}
&{3}&{705}&{39}&{987}\cr
{987}&{6674}&{141}&{2679}&{47}&{329}&{0}&
{-47}&{0}&{3948}&{235}&{6298}\cr
{18}&{141}&{5}&{141}&{4}&{47}&{3}&{0}&
{-2}&{0}&{2}&{94}\cr
{282}&{2679}&{141}&{6909}&{282}&{5499}
&{705}&{3948}&{0}&{-282}&{0}&{1128}\cr
{3}&{47}&{4}&{282}&{13}&{282}&{39}&{235}
&{2}&{0}&{-1}&{0}\cr
{0}&{329}&{47}&{5499}&{282}&{6721}&{987}
&{6298}&{94}&{1128}&{0}&{-94}\cr\endsmallmatrix\right)$.}

\vbox{\noindent
$N=77\ $ $d=155=5\cdot 31,\,\eta=3,\,h=1$:
$\langle 2 \rangle \oplus \langle-31 \rangle \oplus \langle -10 \rangle
(1/2,0,1/2)$
\nobreak
\newline
$\left(\smallmatrix{114}&{1705}&{12}&{31}&{0}&{0}&
{{1}\over{2}}&{{17}\over{2}}&{{4185}\over{2}}&{{221}\over{2}}&
{{7533}\over{2}}&{{7105}\over{2}}&{22785}&{{2449}\over{2}}\cr
{-28}&{-420}&{-3}&{-8}&{0}&{1}&{0}&{-2}&{-510}&{-27}&{-922}&
{-870}&{-5581}&{-300}\cr
{-13}&{-186}&{-1}&{0}&{1}&{0}&{{-1}\over{2}}&{{-3}\over{2}}&
{{-527}\over{2}}&{{-27}\over{2}}&{{-899}\over{2}}&
{{-843}\over{2}}&{-2697}&{{-289}\over{2}}\cr\endsmallmatrix\right)$
\nobreak
\newline
$\left(\smallmatrix
{-2}&{0}&{2}&{124}&{130}&{868}&{49}&{7}&{155}&{3}&{31}&{15}&{62}&{1}\cr
{0}&{-310}&{0}&{1550}&{1860}&{13020}&{775}&{155}&{5115}&{155}&
{3255}&{2635}&{15810}&{775}\cr
{2}&{0}&{-1}&{0}&{10}&{93}&{7}&{3}&{155}&{6}&{155}&{135}&{837}&{43}\cr
{124}&{1550}&{0}&{-62}&{0}&{248}&{31}&{31}&{3255}&{155}&
{4867}&{4495}&{28582}&{1519}\cr
{130}&{1860}&{10}&{0}&{-10}&{0}&{5}&{15}&{2635}&{135}
&{4495}&{4215}&{26970}&{1445}\cr
{868}&{13020}&{93}&{248}&{0}&{-31}&{0}&{62}&{15810}
&{837}&{28582}&{26970}&{173011}&{9300}\cr
{49}&{775}&{7}&{31}&{5}&{0}&{-2}&{1}&{775}&{43}
&{1519}&{1445}&{9300}&{502}\cr
{7}&{155}&{3}&{31}&{15}&{62}&{1}&{-2}&{0}&{2}&{124}
&{130}&{868}&{49}\cr
{155}&{5115}&{155}&{3255}&{2635}&{15810}&{775}&{0}
&{-310}&{0}&{1550}&{1860}&{13020}&{775}\cr
{3}&{155}&{6}&{155}&{135}&{837}&{43}&{2}&{0}&{-1}&{0}
&{10}&{93}&{7}\cr
{31}&{3255}&{155}&{4867}&{4495}&{28582}&{1519}&{124}
&{1550}&{0}&{-62}&{0}&{248}&{31}\cr
{15}&{2635}&{135}&{4495}&{4215}&{26970}&{1445}&{130}&{1860}
&{10}&{0}&{-10}&{0}&{5}\cr
{62}&{15810}&{837}&{28582}&{26970}&{173011}&{9300}&{868}
&{13020}&{93}&{248}&{0}&{-31}&{0}\cr
{1}&{775}&{43}&{1519}&{1445}&{9300}&{502}&{49}&{775}
&{7}&{31}&{5}&{0}&{-2}\cr\endsmallmatrix\right)$.}

\vbox{\noindent
$N=78\ $ $d=165=3\cdot 5\cdot 11,\,\eta=0,\,h=0$:
$\langle 5 \rangle \oplus \langle -6\rangle \oplus \langle -22 \rangle
(0,1/2,1/2)$
\nobreak
\newline
$\pmatrix{0}&{0}&{1}\cr{0}&{1}&{0}\cr
{1}&{{-1}\over{2}}&{{-1}\over{2}}\cr
{33}&{{-55}\over{2}}&{{-15}\over{2}}\cr{1}&{-1}&{0}\cr\endpmatrix
\hskip20pt
\pmatrix{-22}&{0}&{11}&{165}&{0}\cr
{0}&{-6}&{3}&{165}&{6}\cr
{11}&{3}&{-2}&{0}&{2}\cr
{165}&{165}&{0}&{-330}&{0}\cr
{0}&{6}&{2}&{0}&{-1}\cr\endpmatrix$.}

\vbox{\noindent
$N=79\ $ $d=165=3\cdot 5\cdot 11,\,\eta=3,h=1$:
$\langle 15 \rangle \oplus \langle -22 \rangle \oplus \langle -2 \rangle
(0,1/2,1/2)$
\nobreak
\newline
$\left(\smallmatrix{0}&{0}&{1}\cr{0}&{1}&{0}\cr{1}&{0}&{-3}\cr
{44}&{-15}&{-110}\cr{5}&{-2}&{-12}\cr
{22}&{{-21}\over{2}}&{{-99}\over{2}}\cr
{4}&{{-5}\over{2}}&{{-15}\over{2}}\cr
{2}&{{-3}\over{2}}&{{-5}\over{2}}\cr
{22}&{{-35}\over{2}}&{{-33}\over{2}}\cr
{13}&{{-21}\over{2}}&{{-15}\over{2}}\cr
{286}&{{-465}\over{2}}&{{-275}\over{2}}\cr
{27}&{-22}&{-12}\cr{88}&{-72}&{-33}\cr{6}&{-5}&{0}\cr
\endsmallmatrix\right)
\left(\smallmatrix
{-2}&{0}&{6}&{220}&{24}&{99}&{15}&{5}&{33}&{15}&{275}&{24}&{66}&{0}\cr
{0}&{-22}&{0}&{330}&{44}&{231}&{55}&{33}&{385}&{231}&{5115}&
{484}&{1584}&{110}\cr{6}&{0}&{-3}&{0}&{3}&{33}&{15}&{15}&
{231}&{150}&{3465}&{333}&{1122}&{90}\cr{220}&{330}&{0}&{-110}
&{0}&{165}&{165}&{275}&{5115}&{3465}&{81785}&{7920}&{27060}&{2310}\cr
{24}&{44}&{3}&{0}&{-1}&{0}&{10}&{24}&{484}&{333}&{7920}&{769}
&{2640}&{230}\cr{99}&{231}&{33}&{165}&{0}&{-66}&{0}&{66}&{1584}
&{1122}&{27060}&{2640}&{9141}&{825}\cr
{15}&{55}&{15}&{165}&{10}&{0}&{-10}&{0}&{110}&{90}&{2310}
&{230}&{825}&{85}\cr
{5}&{33}&{15}&{275}&{24}&{66}&{0}&{-2}&{0}&{6}&{220}
&{24}&{99}&{15}\cr
{33}&{385}&{231}&{5115}&{484}&{1584}&{110}&{0}&{-22}&
{0}&{330}&{44}&{231}&{55}\cr
{15}&{231}&{150}&{3465}&{333}&{1122}&{90}&{6}&{0}&{-3}
&{0}&{3}&{33}&{15}\cr
{275}&{5115}&{3465}&{81785}&{7920}&{27060}&{2310}&{220}&
{330}&{0}&{-110}&{0}&{165}&{165}\cr
{24}&{484}&{333}&{7920}&{769}&{2640}&{230}&{24}&{44}&{3}
&{0}&{-1}&{0}&{10}\cr
{66}&{1584}&{1122}&{27060}&{2640}&{9141}&{825}&{99}&{231}
&{33}&{165}&{0}&{-66}&{0}\cr
{0}&{110}&{90}&{2310}&{230}&{825}&{85}&{15}&{55}&{15}&{165}
&{10}&{0}&{-10}\cr\endsmallmatrix\right)$.}

\vbox{\noindent
$N=80\ $ $d=165=3\cdot 5\cdot 11,\,\eta=4,\,h=1$:
$\langle 3 \rangle \oplus \langle -5 \rangle \oplus \langle -11 \rangle$
\nobreak
\newline
$\left(\smallmatrix{0}&{0}&{1}\cr{0}&{1}&{0}\cr{11}&{0}&{-6}\cr
{4}&{-1}&{-2}\cr{440}&{-154}&{-205}\cr{91}&{-33}&{-42}\cr
{55}&{-21}&{-25}\cr{27}&{-11}&{-12}\cr{154}&{-66}&{-67}\cr
{70}&{-31}&{-30}\cr{143}&{-66}&{-60}\cr{10}&{-5}&{-4}\cr
{330}&{-176}&{-125}\cr{49}&{-27}&{-18}\cr{15}&{-9}&{-5}\cr
{1}&{-1}&{0}\cr\endsmallmatrix\right)
\left(\smallmatrix{-11}&{0}&{66}&{22}&{2255}&{462}&{275}&{132}&{737}&
{330}&{660}&{44}&{1375}&{198}&{55}&{0}\cr
{0}&{-5}&{0}&{5}&{770}&{165}&{105}&{55}&
{330}&{155}&{330}&{25}&{880}&{135}&{45}&{5}\cr
{66}&{0}&{-33}&{0}&{990}&{231}&{165}&{99}&{660}
&{330}&{759}&{66}&{2640}&{429}&{165}&{33}\cr
{22}&{5}&{0}&{-1}&{0}&{3}&{5}&{5}&{44}&{25}&
{66}&{7}&{330}&{57}&{25}&{7}\cr
{2255}&{770}&{990}&{0}&{-55}&{0}&{55}&{110}&
{1375}&{880}&{2640}&{330}&{18205}&{3300}&{1595}&{550}\cr
{462}&{165}&{231}&{3}&{0}&{-6}&{0}&{12}&{198}&{135}&{429}
&{57}&{3300}&{606}&{300}&{108}\cr
{275}&{105}&{165}&{5}&{55}&{0}&{-5}&{0}&{55}&{45}&{165}
&{25}&{1595}&{300}&{155}&{60}\cr
{132}&{55}&{99}&{5}&{110}&{12}&{0}&{-2}&{0}&{5}&{33}&{7}
&{550}&{108}&{60}&{26}\cr
{737}&{330}&{660}&{44}&{1375}&{198}&{55}&{0}&{-11}&{0}
&{66}&{22}&{2255}&{462}&{275}&{132}\cr
{330}&{155}&{330}&{25}&{880}&{135}&{45}&{5}&{0}&{-5}&{0}
&{5}&{770}&{165}&{105}&{55}\cr
{660}&{330}&{759}&{66}&{2640}&{429}&{165}&{33}&{66}&{0}
&{-33}&{0}&{990}&{231}&{165}&{99}\cr
{44}&{25}&{66}&{7}&{330}&{57}&{25}&{7}&{22}&{5}&{0}&{-1}&
{0}&{3}&{5}&{5}\cr
{1375}&{880}&{2640}&{330}&{18205}&{3300}&{1595}&
{550}&{2255}&{770}&{990}&{0}&{-55}&{0}&{55}&{110}\cr
{198}&{135}&{429}&{57}&{3300}&{606}&{300}&{108}&
{462}&{165}&{231}&{3}&{0}&{-6}&{0}&{12}\cr
{55}&{45}&{165}&{25}&{1595}&{300}&{155}&{60}&{275}
&{105}&{165}&{5}&{55}&{0}&{-5}&{0}\cr
{0}&{5}&{33}&{7}&{550}&{108}&{60}&{26}&{132}&{55}
&{99}&{5}&{110}&{12}&{0}&{-2}\cr\endsmallmatrix\right)$.}

\vbox{\noindent
$N=81\ $ $d=165=3\cdot 5\cdot 11,\,\eta=5,\,h=1$:
$\langle 1 \rangle \oplus \langle -55 \rangle \oplus \langle -3 \rangle$
\nobreak
\newline
$\left(\smallmatrix{0}&{0}&{1}\cr{0}&{1}&{0}\cr{1}&{0}&{-1}\cr
{9}&{-1}&{-3}\cr{75}&{-9}&{-20}\cr
{242}&{-30}&{-55}\cr{120}&{-15}&{-26}\cr
{1320}&{-166}&{-275}\cr{71}&{-9}&{-14}\cr
{39}&{-5}&{-7}\cr{45}&{-6}&{-5}\cr{22}&{-3}&{0}\cr
\endsmallmatrix\right)
\hskip20pt
\left(\smallmatrix
{-3}&{0}&{3}&{9}&{60}&{165}&{78}&{825}&{42}&{21}&{15}&{0}\cr
{0}&{-55}&{0}&{55}&{495}&{1650}&{825}&{9130}&{495}&{275}&{330}&{165}\cr
{3}&{0}&{-2}&{0}&{15}&{77}&{42}&{495}&{29}&{18}&{30}&{22}\cr
{9}&{55}&{0}&{-1}&{0}&{33}&{21}&{275}&{18}&{13}&{30}&{33}\cr
{60}&{495}&{15}&{0}&{-30}&{0}&{15}&{330}&{30}&{30}&{105}&{165}\cr
{165}&{1650}&{77}&{33}&{0}&{-11}&{0}&{165}&{22}&{33}&{165}&{374}\cr
{78}&{825}&{42}&{21}&{15}&{0}&{-3}&{0}&{3}&{9}&{60}&{165}\cr
{825}&{9130}&{495}&{275}&{330}&{165}&{0}&{-55}&{0}&{55}&{495}&{1650}\cr
{42}&{495}&{29}&{18}&{30}&{22}&{3}&{0}&{-2}&{0}&{15}&{77}\cr
{21}&{275}&{18}&{13}&{30}&{33}&{9}&{55}&{0}&{-1}&{0}&{33}\cr
{15}&{330}&{30}&{30}&{105}&{165}&{60}&{495}&{15}&{0}&{-30}&{0}\cr
{0}&{165}&{22}&{33}&{165}&{374}&{165}&{1650}&{77}&{33}&{0}&{-11}\cr
\endsmallmatrix\right)$.}

\vbox{\noindent
$N=82\ $ $d=165=3\cdot 5\cdot 11,\,\eta=6,\,h=0$
$\langle 1 \rangle \oplus \langle -15 \rangle \oplus \langle -11 \rangle$
\nobreak
\newline
$\pmatrix{0}&{0}&{1}\cr{0}&{1}&{0}\cr{3}&{0}&{-1}\cr
{99}&{-11}&{-27}\cr{20}&{-3}&{-5}\cr
{5}&{-1}&{-1}\cr{3}&{-1}&{0}\cr\endpmatrix
\hskip20pt
\pmatrix{-11}&{0}&{11}&{297}&{55}&{11}&{0}\cr
{0}&{-15}&{0}&{165}&{45}&{15}&{15}\cr
{11}&{0}&{-2}&{0}&{5}&{4}&{9}\cr
{297}&{165}&{0}&{-33}&{0}&{33}&{132}\cr
{55}&{45}&{5}&{0}&{-10}&{0}&{15}\cr
{11}&{15}&{4}&{33}&{0}&{-1}&{0}\cr
{0}&{15}&{9}&{132}&{15}&{0}&{-6}\cr\endpmatrix$.}

\vbox{\noindent
$N=83\ $ $d=170=2\cdot 5\cdot 17,\,\eta=1,\,h=1$:
$\langle 2 \rangle \oplus \langle-10 \rangle \oplus \langle-34\rangle
(1/2,1/2,0)$
\nobreak
\newline
$\left(\smallmatrix{205}&{{147}\over{2}}&{{1}\over{2}}&{0}&{0}&{4}&{{153}
\over{2}}&{65}&{{123}\over{2}}&{{89}\over{2}}&{3060}&{720}&{176}&{{1377}
\over{2}}\cr{-39}&{{-29}\over{2}}&{{-1}\over{2}}&{0}&{1}&
{0}&{{-17}\over{2}}&{-9}&{{-19}\over{2}}&{{-15}\over{2}}&{-544}&{-129}&
{-32}&{{-255}\over{2}}\cr{-45}&{-16}&{0}&{1}&{0}&{-1}&{-18}&{-15}&
{-14}&{-10}&{-681}&{-160}&{-39}&{-152}\cr\endsmallmatrix\right)$
\nobreak
\newline
$\left(\smallmatrix
{-10}&{0}&{10}&{1530}&{390}&{110}&{510}&{190}&{90}&{20}&
{510}&{90}&{10}&{0}\cr{0}
&{-2}&{1}&{544}&{145}&{44}&{221}&{90}&{47}
&{14}&{476}&{95}&{16}&{34}\cr
{10}&{1}&{-2}&{0}&{5}&{4}&{34}&{20}&{14}&{7}
&{340}&{75}&{16}&{51}\cr
{1530}&{544}&{0}&{-34}&{0}&{34}&{612}&{510}&
{476}&{340}&{23154}&{5440}&{1326}&{5168}\cr
{390}&{145}&{5}&{0}&{-10}&{0}&{85}&{90}&{95}&{75}&
{5440}&{1290}&{320}&{1275}\cr
{110}&{44}&{4}&{34}&{0}&{-2}&{0}&{10}&{16}&{16}&
{1326}&{320}&{82}&{340}\cr
{510}&{221}&{34}&{612}&{85}&{0}&{-34}&{0}&{34}&
{51}&{5168}&{1275}&{340}&{1479}\cr
{190}&{90}&{20}&{510}&{90}&{10}&{0}&{-10}&{0}
&{10}&{1530}&{390}&{110}&{510}\cr
{90}&{47}&{14}&{476}&{95}&{16}&{34}&{0}&{-2}
&{1}&{544}&{145}&{44}&{221}\cr
{20}&{14}&{7}&{340}&{75}&{16}&{51}&{10}&{1}&
{-2}&{0}&{5}&{4}&{34}\cr
{510}&{476}&{340}&{23154}&{5440}&{1326}&{5168}&{1530}&
{544}&{0}&{-34}&{0}&{34}&{612}\cr
{90}&{95}&{75}&{5440}&{1290}&{320}&{1275}&{390}&{145}&
{5}&{0}&{-10}&{0}&{85}\cr
{10}&{16}&{16}&{1326}&{320}&{82}&{340}&{110}&{44}&{4}&
{34}&{0}&{-2}&{0}\cr
{0}&{34}&{51}&{5168}&{1275}&{340}&{1479}&{510}&{221}&
{34}&{612}&{85}&{0}&{-34}\cr\endsmallmatrix\right)$.}

\vbox{\noindent
$N=84\ $ $d=195=3\cdot 5\cdot 13,\,\eta=3,\,h=0$:
$\langle 15 \rangle \oplus \langle-13 \rangle \oplus \langle-1\rangle$
\nobreak
\newline
$\pmatrix{0}&{0}&{1}\cr{0}&{1}&{0}\cr{1}&{0}&{-5}\cr
{13}&{-9}&{-39}\cr{1}&{-1}&{-2}\cr{13}&{-15}&{0}\cr\endpmatrix
\hskip20pt
\pmatrix{-1}&{0}&{5}&{39}&{2}&{0}\cr
{0}&{-13}&{0}&{117}&{13}&{195}\cr
{5}&{0}&{-10}&{0}&{5}&{195}\cr
{39}&{117}&{0}&{-39}&{0}&{780}\cr
{2}&{13}&{5}&{0}&{-2}&{0}\cr
{0}&{195}&{195}&{780}&{0}&{-390}\cr\endpmatrix$.}

\vbox{\noindent
$N=85\ $ $d=195=3\cdot 5\cdot 13,\,\eta=5,h=1$:
$\langle 6 \rangle \oplus \langle-26 \rangle \oplus \langle -5 \rangle
(1/2,1/2,0)$
\nobreak
\newline
$\left(\smallmatrix{0}&{0}&{1}\cr{0}&{1}&{0}\cr{5}&{0}&{-6}\cr
{5}&{-1}&{-5}\cr{{11}\over{2}}&{{-3}\over{2}}&{-5}\cr
{{15}\over{2}}&{{-5}\over{2}}&{-6}\cr
{{39}\over{2}}&{{-15}\over{2}}&{-13}\cr
{{35}\over{2}}&{{-15}\over{2}}&{-9}\cr
{{11}\over{2}}&{{-5}\over{2}}&{-2}\cr{2}&{-1}&{0}\cr
\endsmallmatrix\right)
\hskip20pt
\left(\smallmatrix{-5}&{0}&{30}&{25}&{25}&{30}&{65}&{45}&{10}&{0}\cr
{0}&{-26}&{0}&{26}&{39}&{65}&{195}&{195}&{65}&{26}\cr
{30}&{0}&{-30}&{0}&{15}&{45}&{195}&{255}&{105}&{60}\cr
{25}&{26}&{0}&{-1}&{1}&{10}&{65}&{105}&{50}&{34}\cr
{25}&{39}&{15}&{1}&{-2}&{0}&{26}&{60}&{34}&{27}\cr
{30}&{65}&{45}&{10}&{0}&{-5}&{0}&{30}&{25}&{25}\cr
{65}&{195}&{195}&{65}&{26}&{0}&{-26}&{0}&{26}&{39}\cr
{45}&{195}&{255}&{105}&{60}&{30}&{0}&{-30}&{0}&{15}\cr
{10}&{65}&{105}&{50}&{34}&{25}&{26}&{0}&{-1}&{1}\cr
{0}&{26}&{60}&{34}&{27}&{25}&{39}&{15}&{1}&{-2}\cr
\endsmallmatrix\right)$.}

\vbox{\noindent
$N=86\ $ $d=195=3\cdot 5\cdot 13,\,\eta =6,\,h=0$:
$\langle 3 \rangle \oplus \langle-65 \rangle \oplus \langle -1 \rangle$
\nobreak
\newline
$\pmatrix{0}&{0}&{1}\cr{0}&{1}&{0}\cr{1}&{0}&{-3}\cr
{28}&{-5}&{-27}\cr{115}&{-21}&{-105}\cr{91}&{-17}&{-78}\cr
{15}&{-3}&{-10}\cr{13}&{-3}&{0}\cr\endpmatrix
\hskip20pt
\pmatrix{-1}&{0}&{3}&{27}&{105}&{78}&{10}&{0}\cr
{0}&{-65}&{0}&{325}&{1365}&{1105}&{195}&{195}\cr
{3}&{0}&{-6}&{3}&{30}&{39}&{15}&{39}\cr
{27}&{325}&{3}&{-2}&{0}&{13}&{15}&{117}\cr
{105}&{1365}&{30}&{0}&{-15}&{0}&{30}&{390}\cr
{78}&{1105}&{39}&{13}&{0}&{-26}&{0}&{234}\cr
{10}&{195}&{15}&{15}&{30}&{0}&{-10}&{0}\cr
{0}&{195}&{39}&{117}&{390}&{234}&{0}&{-78}\cr\endpmatrix$.}

\vbox{\noindent
$N=87\ $ $d=205=5\cdot 41,\,\eta =1,\,h=1$:
$\langle 2 \rangle \oplus \langle-41 \rangle \oplus \langle -10 \rangle
(1/2,0,1/2)$
\nobreak
\newline
$\left(\smallmatrix{107}&{{185}\over{2}}&{{533}\over{2}}
&{32}&{820}&{9}&{0}&{0}
&{{1}\over{2}}&{5}&{{95}\over{2}}&{{1763}\over{2}}&{248}
&{10660}&{327}&{840}&{6888}&{{727}\over{2}}\cr{-23}&
{-20}&{-58}&{-7}&{-180}&{-2}&{0}&{1}&{0}&{-1}&{-10}&
{-188}&{-53}&{-2280}&{-70}&{-180}&{-1477}&{-78}\cr
{-11}&{{-17}\over{2}}&{{-41}\over{2}}&{-2}&{-41}&{0}&{1}&
{0}&{{-1}\over{2}}&{-1}&{{-13}\over{2}}&{{-205}\over{2}}&
{-28}&{-1189}&{-36}&{-91}&{-738}&{{-77}\over{2}}\cr
\endsmallmatrix\right)$
\nobreak
\newline
$\pmatrix
A&B\\
B&A\\
\endpmatrix\ \
\text{where\ \ }
\pmatrix
A\\B\\
\endpmatrix =
\left(\smallmatrix
{-1}&{0}&{82}&{27}&{1230}&{40}&{110}&{943}&{52}\cr
{0}&{-10}&{0}&{10}&{615}&{25}&{85}&{820}&{50}\cr
{82}&{0}&{-82}&{0}&{615}&{41}&{205}&{2378}&{164}\cr
{27}&{10}&{0}&{-1}&{0}&{2}&{20}&{287}&{22}\cr
{1230}&{615}&{615}&{0}&{-410}&{0}&{410}&{7380}&{615}\cr
{40}&{25}&{41}&{2}&{0}&{-2}&{0}&{82}&{9}\cr
{110}&{85}&{205}&{20}&{410}&{0}&{-10}&{0}&{5}\cr
{943}&{820}&{2378}&{287}&{7380}&{82}&{0}&{-41}&{0}\cr
{52}&{50}&{164}&{22}&{615}&{9}&{5}&{0}&{-2}\cr
{17}&{20}&{82}&{13}&{410}&{8}&{10}&{41}&{0}\cr
{20}&{35}&{205}&{40}&{1435}&{35}&{65}&{410}&{15}\cr
{82}&{205}&{1763}&{410}&{16195}&{451}&{1025}&{7708}&{369}\cr
{13}&{40}&{410}&{101}&{4100}&{118}&{280}&{2173}&{108}\cr
{410}&{1435}&{16195}&{4100}&{168510}&{4920}&{11890}&{93480}&{4715}\cr
{8}&{35}&{451}&{118}&{4920}&{146}&{360}&{2870}&{147}\cr
{10}&{65}&{1025}&{280}&{11890}&{360}&{910}&{7380}&{385}\cr
{41}&{410}&{7708}&{2173}&{93480}&{2870}&{7380}&{60557}&{3198}\cr
{0}&{15}&{369}&{108}&{4715}&{147}&{385}&{3198}&{171}\cr
\endsmallmatrix\right)$.}

\vbox{\noindent
$N=88\ $ $d=210=2\cdot 3\cdot 5\cdot 7,\,\eta=2,\,h=0$:
$\langle 2 \rangle \oplus \langle -10\rangle \oplus \langle -42 \rangle
(1/2,1/2,0)$
\nobreak
\newline
$\pmatrix{0}&{0}&{1}\cr{0}&{1}&{0}\cr{9}&{0}&{-2}\cr
{70}&{-6}&{-15}\cr{{189}\over{2}}&{{-21}\over{2}}&{-20}\cr
{5}&{-1}&{-1}\cr{{1}\over{2}}&{{-1}\over{2}}&{0}\cr\endpmatrix
\hskip20pt
\pmatrix{-42}&{0}&{84}&{630}&{840}&{42}&{0}\cr
{0}&{-10}&{0}&{60}&{105}&{10}&{5}\cr
{84}&{0}&{-6}&{0}&{21}&{6}&{9}\cr
{630}&{60}&{0}&{-10}&{0}&{10}&{40}\cr
{840}&{105}&{21}&{0}&{-42}&{0}&{42}\cr
{42}&{10}&{6}&{10}&{0}&{-2}&{0}\cr
{0}&{5}&{9}&{40}&{42}&{0}&{-2}\cr\endpmatrix$.}

\vbox{\noindent
$N=89\ $ $d=210=2\cdot 3\cdot 5\cdot 7,\,\eta=4,\,h=1$:
$\langle 30 \rangle \oplus \langle -14\rangle \oplus \langle -2 \rangle
(0,1/2,1/2)$
\nobreak
\newline
$\pmatrix{0}&{0}&{1}\cr{0}&{1}&{0}\cr
{1}&{0}&{-4}\cr{8}&{-3}&{-30}\cr
{7}&{{-7}\over{2}}&{{-51}\over{2}}\cr
{7}&{{-9}\over{2}}&{{-49}\over{2}}\cr
{1}&{-1}&{-3}\cr{1}&{{-3}\over{2}}&{{-3}\over{2}}\cr\endpmatrix
\hskip20pt
\pmatrix
{-2}&{0}&{8}&{60}&{51}&{49}&{6}&{3}\cr
{0}&{-14}&{0}&{42}&{49}&{63}&{14}&{21}\cr
{8}&{0}&{-2}&{0}&{6}&{14}&{6}&{18}\cr
{60}&{42}&{0}&{-6}&{3}&{21}&{18}&{87}\cr
{51}&{49}&{6}&{3}&{-2}&{0}&{8}&{60}\cr
{49}&{63}&{14}&{21}&{0}&{-14}&{0}&{42}\cr
{6}&{14}&{6}&{18}&{8}&{0}&{-2}&{0}\cr
{3}&{21}&{18}&{87}&{60}&{42}&{0}&{-6}\cr\endpmatrix$.}

\vbox{\noindent
$N=90\ $ $d=219=3\cdot 73,\,\eta=0,\,h=0$:
$\langle 3 \rangle \oplus \langle -73\rangle \oplus \langle-1\rangle$
\nobreak
\newline
$\pmatrix{0}&{0}&{1}\cr{0}&{1}&{0}\cr{1}&{0}&{-3}\cr
{6}&{-1}&{-6}\cr{28}&{-5}&{-23}\cr{49}&{-9}&{-36}\cr
{219}&{-41}&{-146}\cr{5}&{-1}&{-2}\cr{73}&{-15}&{0}\cr\endpmatrix
\pmatrix{-1}&{0}&{3}&{6}&{23}&{36}&{146}&{2}&{0}\cr
{0}&{-73}&{0}&{73}&{365}&{657}&{2993}&{73}&{1095}\cr
{3}&{0}&{-6}&{0}&{15}&{39}&{219}&{9}&{219}\cr
{6}&{73}&{0}&{-1}&{1}&{9}&{73}&{5}&{219}\cr
{23}&{365}&{15}&{1}&{-2}&{3}&{73}&{9}&{657}\cr
{36}&{657}&{39}&{9}&{3}&{-6}&{0}&{6}&{876}\cr
{146}&{2993}&{219}&{73}&{73}&{0}&{-146}&{0}&{3066}\cr
{2}&{73}&{9}&{5}&{9}&{6}&{0}&{-2}&{0}\cr
{0}&{1095}&{219}&{219}&{657}&{876}&{3066}&{0}&{-438}\cr\endpmatrix$.}

\vbox{\noindent
$N=91\ $ $d=231=3\cdot 7\cdot 11,\,\eta=1,\,h=1$:
$\langle 14 \rangle \oplus \langle -22\rangle \oplus \langle-3\rangle
(1/2,1/2,0)$
\nobreak
\newline
$\left(\smallmatrix{0}&{0}&{1}\cr{0}&{1}&{0}\cr{3}&{0}&{-7}\cr
{{3}\over{2}}&{{-1}\over{2}}&{-3}\cr
{{99}\over{2}}&{{-45}\over{2}}&{-88}\cr
{{35}\over{2}}&{{-17}\over{2}}&{-30}\cr{60}&{-30}&{-101}\cr
{132}&{-67}&{-220}\cr{81}&{-42}&{-133}\cr
{{9}\over{2}}&{{-5}\over{2}}&{-7}\cr
{{33}\over{2}}&{{-21}\over{2}}&{-22}\cr
{{1}\over{2}}&{{-1}\over{2}}&{0}\cr\endsmallmatrix\right)
\hskip20pt
\left(\smallmatrix
{-3}&{0}&{21}&{9}&{264}&{90}&{303}&{660}&{399}&{21}&{66}&{0}\cr
{0}&{-22}&{0}&{11}&{495}&{187}&{660}&{1474}&{924}&{55}&{231}&{11}\cr
{21}&{0}&{-21}&{0}&{231}&{105}&{399}&{924}&{609}&{42}&{231}&{21}\cr
{9}&{11}&{0}&{-1}&{0}&{4}&{21}&{55}&{42}&{4}&{33}&{5}\cr
{264}&{495}&{231}&{0}&{-66}&{0}&{66}&{231}&{231}&{33}&{429}&{99}\cr
{90}&{187}&{105}&{4}&{0}&{-2}&{0}&{11}&{21}&{5}&{99}&{29}\cr
{303}&{660}&{399}&{21}&{66}&{0}&{-3}&{0}&{21}&{9}&{264}&{90}\cr
{660}&{1474}&{924}&{55}&{231}&{11}&{0}&{-22}&{0}&{11}&{495}&{187}\cr
{399}&{924}&{609}&{42}&{231}&{21}&{21}&{0}&{-21}&{0}&{231}&{105}\cr
{21}&{55}&{42}&{4}&{33}&{5}&{9}&{11}&{0}&{-1}&{0}&{4}\cr
{66}&{231}&{231}&{33}&{429}&{99}&{264}&{495}&{231}&{0}&{-66}&{0}\cr
{0}&{11}&{21}&{5}&{99}&{29}&{90}&{187}&{105}&{4}&{0}&{-2}\cr
\endsmallmatrix\right)$.}

\vbox{\noindent
$N=92\ $ $d=231=3\cdot 7\cdot 11,\,\eta=7,\,h=0$:
$\langle 6 \rangle \oplus \langle -14\rangle \oplus \langle-11\rangle
(1/2,1/2,0)$
\nobreak
\newline
$\pmatrix{0}&{0}&{1}\cr{0}&{1}&{0}\cr
{4}&{0}&{-3}\cr{{77}\over{2}}&{{-11}\over{2}}&{-28}\cr
{{3}\over{2}}&{{-1}\over{2}}&{-1}\cr
{{1}\over{2}}&{{-1}\over{2}}&{0}\cr\endpmatrix
\hskip20pt
\pmatrix{-11}&{0}&{33}&{308}&{11}&{0}\cr
{0}&{-14}&{0}&{77}&{7}&{7}\cr{33}&{0}&{-3}&{0}&{3}&{12}\cr
{308}&{77}&{0}&{-154}&{0}&{77}\cr{11}&{7}&{3}&{0}&{-1}&{1}\cr
{0}&{7}&{12}&{77}&{1}&{-2}\cr\endpmatrix$.}

\vbox{\noindent
$N=93\ $ $d=255=3\cdot 5\cdot 17,\,\eta=0,\,h=1$:
$\langle 3 \rangle \oplus \langle -17\rangle \oplus \langle-5\rangle$
\nobreak
\newline
$\left(\smallmatrix{0}&{0}&{1}\cr{0}&{1}&{0}\cr{1}&{0}&{-1}\cr
{7}&{-2}&{-4}\cr{34}&{-11}&{-17}\cr{8}&{-3}&{-3}\cr
{25}&{-10}&{-6}\cr{170}&{-69}&{-34}\cr{49}&{-20}&{-9}\cr
{78}&{-32}&{-13}\cr{221}&{-91}&{-34}\cr{7}&{-3}&{0}\cr
\endsmallmatrix\right)
\hskip20pt
\left(\smallmatrix
{-5}&{0}&{5}&{20}&{85}&{15}&{30}&{170}&{45}&{65}&{170}&{0}\cr
{0}&{-17}&{0}&{34}&{187}&{51}&{170}&{1173}&{340}&{544}&{1547}&{51}\cr
{5}&{0}&{-2}&{1}&{17}&{9}&{45}&{340}&{102}&{169}&{493}&{21}\cr
{20}&{34}&{1}&{-1}&{0}&{6}&{65}&{544}&{169}&{290}&{867}&{45}\cr
{85}&{187}&{17}&{0}&{-34}&{0}&{170}&{1547}&{493}&{867}&{2635}&{153}\cr
{15}&{51}&{9}&{6}&{0}&{-6}&{0}&{51}&{21}&{45}&{153}&{15}\cr
{30}&{170}&{45}&{65}&{170}&{0}&{-5}&{0}&{5}&{20}&{85}&{15}\cr
{170}&{1173}&{340}&{544}&{1547}&{51}&{0}&{-17}&{0}&{34}&{187}&{51}\cr
{45}&{340}&{102}&{169}&{493}&{21}&{5}&{0}&{-2}&{1}&{17}&{9}\cr
{65}&{544}&{169}&{290}&{867}&{45}&{20}&{34}&{1}&{-1}&{0}&{6}\cr
{170}&{1547}&{493}&{867}&{2635}&{153}&{85}&{187}&{17}&{0}&{-34}&{0}\cr
{0}&{51}&{21}&{45}&{153}&{15}&{15}&{51}&{9}&{6}&{0}&{-6}\cr
\endsmallmatrix\right)$.}

\vbox{\noindent
$N=94\ $ $d=255=3\cdot 5\cdot 17,\,\eta=6,h=0$:
$\langle 3 \rangle \oplus \langle -85\rangle \oplus \langle-1\rangle$
\nobreak
\newline
$\pmatrix
{0}&{0}&{1}\cr{0}&{1}&{0}\cr{1}&{0}&{-3}\cr
{6}&{-1}&{-5}\cr{85}&{-15}&{-51}\cr{5}&{-1}&{0}\cr\endpmatrix
\hskip20pt
\pmatrix{-1}&{0}&{3}&{5}&{51}&{0}\cr{0}&{-85}&{0}&{85}&{1275}&{85}\cr
{3}&{0}&{-6}&{3}&{102}&{15}\cr{5}&{85}&{3}&{-2}&{0}&{5}\cr
{51}&{1275}&{102}&{0}&{-51}&{0}\cr{0}&{85}&{15}&{5}&{0}&{-10}\cr
\endpmatrix$.}

\vbox{\noindent
$N=95\ $ $d=273=3\cdot 7\cdot 13,\,\eta=0,\,h=1$:
$\langle 14 \rangle \oplus \langle -13\rangle \oplus \langle-6\rangle
(1/2,0,1/2)$
\nobreak
\newline
$\left(\smallmatrix{93}&{-84}&{-70}\cr{2}&{-2}&{-1}\cr{13}&{-14}&{0}\cr
{0}&{0}&{1}\cr{0}&{1}&{0}\cr{{3}\over{2}}&{0}&{{-7}\over{2}}\cr
{{5}\over{2}}&{-2}&{{-5}\over{2}}\cr
{{793}\over{2}}&{-350}&{{-637}\over{2}}\cr
{{189}\over{2}}&{-84}&{{-149}\over{2}}\cr{234}&{-209}&{-182}\cr
\endsmallmatrix\right)
\hskip20pt
\left(\smallmatrix
{-42}&{0}&{1638}&{420}&{1092}&{483}&{21}&{273}&{21}&{0}\cr
{0}&{-2}&{0}&{6}&{26}&{21}&{3}&{91}&{15}&{26}\cr
{1638}&{0}&{-182}&{0}&{182}&{273}&{91}&{8463}&{1911}&{4550}\cr
{420}&{6}&{0}&{-6}&{0}&{21}&{15}&{1911}&{447}&{1092}\cr
{1092}&{26}&{182}&{0}&{-13}&{0}&{26}&{4550}&{1092}&{2717}\cr
{483}&{21}&{273}&{21}&{0}&{-42}&{0}&{1638}&{420}&{1092}\cr
{21}&{3}&{91}&{15}&{26}&{0}&{-2}&{0}&{6}&{26}\cr
{273}&{91}&{8463}&{1911}&{4550}&{1638}&{0}&{-182}&{0}&{182}\cr
{21}&{15}&{1911}&{447}&{1092}&{420}&{6}&{0}&{-6}&{0}\cr
{0}&{26}&{4550}&{1092}&{2717}&{1092}&{26}&{182}&{0}&{-13}\cr
\endsmallmatrix\right)$.}

\vbox{\noindent
$N=96\ $ $d=273=3\cdot 7\cdot 13,\,\eta=6,\,h=1$:
$\langle 26\rangle \oplus \langle -7 \rangle \oplus \langle -6\rangle
(1/2,0,1/2)$
\nobreak
\newline
$\left(\smallmatrix{{147}\over{2}}&{{261}\over{2}}&{{35}\over{2}}&
{{23}\over{2}}&{{1}\over{2}}&{0}&{0}&{6}&{7}&{34}&{17}&{{273}\over{2}}\cr
{-64}&{-117}&{-17}&{-13}&{-1}&{0}&{1}&{0}&{-4}&{-26}&{-14}&{-117}\cr
{{-273}\over{2}}&{{-481}\over{2}}&{{-63}\over{2}}&
{{-39}\over{2}}&{{-1}\over{2}}&{1}&{0}&{-13}&{-14}&{-65}&
{-32}&{{-509}\over{2}}\cr\endsmallmatrix\right)$
\nobreak
\newline
$\left(\smallmatrix
{-7}&{0}&{28}&{182}&{98}&{819}&{448}&{819}&{119}&{91}&{7}&{0}\cr
{0}&{-78}&{0}&{234}&{156}&{1443}&{819}&{1599}&{273}&{273}&{39}&{78}\cr
{28}&{0}&{-14}&{0}&{14}&{189}&{119}&{273}&{63}&{91}&{21}&{84}\cr
{182}&{234}&{0}&{-26}&{0}&{117}&{91}&{273}&{91}&{195}&{65}&{390}\cr
{98}&{156}&{14}&{0}&{-2}&{3}&{7}&{39}&{21}&{65}&{27}&{192}\cr
{819}&{1443}&{189}&{117}&{3}&{-6}&{0}&{78}&{84}&{390}&{192}&{1527}\cr
{448}&{819}&{119}&{91}&{7}&{0}&{-7}&{0}&{28}&{182}&{98}&{819}\cr
{819}&{1599}&{273}&{273}&{39}&{78}&{0}&{-78}&{0}&{234}&{156}&{1443}\cr
{119}&{273}&{63}&{91}&{21}&{84}&{28}&{0}&{-14}&{0}&{14}&{189}\cr
{91}&{273}&{91}&{195}&{65}&{390}&{182}&{234}&{0}&{-26}&{0}&{117}\cr
{7}&{39}&{21}&{65}&{27}&{192}&{98}&{156}&{14}&{0}&{-2}&{3}\cr
{0}&{78}&{84}&{390}&{192}&{1527}&{819}&{1443}&{189}&{117}&{3}&{-6}\cr
\endsmallmatrix\right)$.}

\vbox{\noindent
$N=97\ $ $d=285=3\cdot 5\cdot 19,\,\eta=3,\,h=1$:
$\langle 10\rangle \oplus \langle -38 \rangle \oplus \langle -3\rangle
(1/2,1/2,0)$
\nobreak
\newline
$\left(\smallmatrix
{0}&{0}&{1}\cr{0}&{1}&{0}\cr{1}&{0}&{-2}\cr{57}&{-12}&{-95}\cr
{{19}\over{2}}&{{-5}\over{2}}&{-15}\cr{13}&{-4}&{-19}\cr
{266}&{-85}&{-380}\cr{18}&{-6}&{-25}\cr{57}&{-20}&{-76}\cr
{8}&{-3}&{-10}\cr{114}&{-45}&{-133}\cr
{{11}\over{2}}&{{-5}\over{2}}&{-5}\cr{2}&{-1}&{-1}\cr
{19}&{-10}&{0}\cr\endsmallmatrix\right)
\left(\smallmatrix
{-3}&{0}&{6}&{285}&{45}&{57}&{1140}&{75}&{228}&{30}&{399}&{15}&{3}&{0}\cr
{0}&{-38}&{0}&{456}&{95}&{152}&{3230}&{228}&{760}&
{114}&{1710}&{95}&{38}&{380}\cr
{6}&{0}&{-2}&{0}&{5}&{16}&{380}&{30}&{114}&{20}&{342}
&{25}&{14}&{190}\cr
{285}&{456}&{0}&{-57}&{0}&{171}&{4560}&{399}&{1710}&
{342}&{6555}&{570}&{399}&{6270}\cr
{45}&{95}&{5}&{0}&{-10}&{0}&{95}&{15}&{95}&{25}&{570}&
{60}&{50}&{855}\cr
{57}&{152}&{16}&{171}&{0}&{-1}&{0}&{3}&{38}&{14}&{399}&
{50}&{51}&{950}\cr
{1140}&{3230}&{380}&{4560}&{95}&{0}&{-190}&{0}&{380}&
{190}&{6270}&{855}&{950}&{18240}\cr
{75}&{228}&{30}&{399}&{15}&{3}&{0}&{-3}&{0}&{6}&{285}&
{45}&{57}&{1140}\cr
{228}&{760}&{114}&{1710}&{95}&{38}&{380}&{0}&{-38}&{0}&
{456}&{95}&{152}&{3230}\cr
{30}&{114}&{20}&{342}&{25}&{14}&{190}&{6}&{0}&{-2}&{0}&
{5}&{16}&{380}\cr
{399}&{1710}&{342}&{6555}&{570}&{399}&{6270}&{285}&{456}&
{0}&{-57}&{0}&{171}&{4560}\cr
{15}&{95}&{25}&{570}&{60}&{50}&{855}&{45}&{95}&{5}&{0}&
{-10}&{0}&{95}\cr
{3}&{38}&{14}&{399}&{50}&{51}&{950}&{57}&{152}&{16}&{171}
&{0}&{-1}&{0}\cr
{0}&{380}&{190}&{6270}&{855}&{950}&{18240}&{1140}&{3230}&
{380}&{4560}&{95}&{0}&{-190}\cr\endsmallmatrix\right)$.}

\vbox{\noindent
$N=98\ $ $d=285=3\cdot 5\cdot 19,\,\eta=5,\,h=0$:
$\langle 1 \rangle \oplus \langle -95\rangle \oplus \langle -3\rangle$
\nobreak
\newline
$\pmatrix{0}&{0}&{1}\cr{0}&{1}&{0}\cr{1}&{0}&{-1}\cr
{285}&{-24}&{-95}\cr{11}&{-1}&{-3}\cr
{30}&{-3}&{-5}\cr{19}&{-2}&{0}\cr\endpmatrix
\hskip20pt
\pmatrix
{-3}&{0}&{3}&{285}&{9}&{15}&{0}\cr
{0}&{-95}&{0}&{2280}&{95}&{285}&{190}\cr
{3}&{0}&{-2}&{0}&{2}&{15}&{19}\cr
{285}&{2280}&{0}&{-570}&{0}&{285}&{855}\cr
{9}&{95}&{2}&{0}&{-1}&{0}&{19}\cr
{15}&{285}&{15}&{285}&{0}&{-30}&{0}\cr
{0}&{190}&{19}&{855}&{19}&{0}&{-19}\cr\endpmatrix$.}

\vbox{\noindent
$N=99\ $ $d=291=3\cdot 97,\,\eta=0,\,h=1$:
$\langle 3 \rangle \oplus \langle -97 \rangle \oplus \langle -1\rangle$
\nobreak
\newline
$\left(\smallmatrix
{0}&{0}&{1}&{97}&{13}&{25}&{2134}&{80}&{73}&{18}&{582}&
{53}&{485}&{29}&{23}&{1358}&{40}&{17}\cr
{0}&{1}&{0}&{-14}&{-2}&{-4}&{-345}&{-13}&{-12}&{-3}&{-98}&{-9}
&{-83}&{-5}&{-4}&{-237}&{-7}&{-3}\cr
{1}&{0}&{-3}&{-97}&{-11}&{-18}&{-1455}&{-53}&{-45}&{-10}&{-291}
&{-24}&{-194}&{-10}&{-6}&{-291}&{-7}&{0}\cr
\endsmallmatrix\right)$
\nobreak
\newline
$\pmatrix
A&B\\
B&A\\
\endpmatrix
\ \text{where\ \ }
\pmatrix
A\\
B\\
\endpmatrix =
\left(\smallmatrix
{-1}&{0}&{3}&{97}&{11}&{18}&{1455}&{53}&{45}\cr
{0}&{-97}&{0}&{1358}&{194}&{388}&{33465}&{1261}&{1164}\cr
{3}&{0}&{-6}&{0}&{6}&{21}&{2037}&{81}&{84}\cr
{97}&{1358}&{0}&{-194}&{0}&{97}&{11349}&{485}&{582}\cr
{11}&{194}&{6}&{0}&{-2}&{1}&{291}&{15}&{24}\cr
{18}&{388}&{21}&{97}&{1}&{-1}&{0}&{2}&{9}\cr
{1455}&{33465}&{2037}&{11349}&{291}&{0}&{-582}&{0}&{291}\cr
{53}&{1261}&{81}&{485}&{15}&{2}&{0}&{-2}&{3}\cr
{45}&{1164}&{84}&{582}&{24}&{9}&{291}&{3}&{-6}\cr
{10}&{291}&{24}&{194}&{10}&{6}&{291}&{7}&{0}\cr
{291}&{9506}&{873}&{8051}&{485}&{388}&{22989}&{679}&{291}\cr
{24}&{873}&{87}&{873}&{57}&{51}&{3201}&{99}&{51}\cr
{194}&{8051}&{873}&{9603}&{679}&{679}&{45105}&{1455}&{873}\cr
{10}&{485}&{57}&{679}&{51}&{55}&{3783}&{125}&{81}\cr
{6}&{388}&{51}&{679}&{55}&{65}&{4656}&{158}&{111}\cr
{291}&{22989}&{3201}&{45105}&{3783}&{4656}&{339306}&{11640}&{8439}\cr
{7}&{679}&{99}&{1455}&{125}&{158}&{11640}&{402}&{297}\cr
{0}&{291}&{51}&{873}&{81}&{111}&{8439}&{297}&{231}\cr
\endsmallmatrix\right)$.}

\vbox{\noindent
$N=100\ $ $d=330=2\cdot 3\cdot 5\cdot 11,\,\eta=3,\,h=0$:
$\langle 6 \rangle \oplus \langle -110 \rangle \oplus \langle -2\rangle
(1/2,1/2,0)$
\nobreak
\newline
$\pmatrix{0}&{0}&{1}\cr{0}&{1}&{0}\cr{1}&{0}&{-2}\cr
{154}&{-18}&{-231}\cr{{55}\over{2}}&{{-7}\over{2}}&{-40}\cr
{{33}\over{2}}&{{-5}\over{2}}&{-22}\cr{5}&{-1}&{-5}\cr
{{11}\over{2}}&{{-3}\over{2}}&{0}\cr\endpmatrix
\hskip20pt
\pmatrix{-2}&{0}&{4}&{462}&{80}&{44}&{10}&{0}\cr
{0}&{-110}&{0}&{1980}&{385}&{275}&{110}&{165}\cr
{4}&{0}&{-2}&{0}&{5}&{11}&{10}&{33}\cr
{462}&{1980}&{0}&{-66}&{0}&{132}&{330}&{2112}\cr
{80}&{385}&{5}&{0}&{-10}&{0}&{40}&{330}\cr
{44}&{275}&{11}&{132}&{0}&{-22}&{0}&{132}\cr
{10}&{110}&{10}&{330}&{40}&{0}&{-10}&{0}\cr
{0}&{165}&{33}&{2112}&{330}&{132}&{0}&{-66}\cr\endpmatrix$.}

\vbox{\noindent
$N=101\ $ $d=345=3\cdot 5\cdot 23,\,\eta=6,\,h=0$:
$\langle 1 \rangle \oplus \langle -23 \rangle \oplus \langle -15 \rangle$
\nobreak
\newline
$\pmatrix{0}&{0}&{1}\cr{0}&{1}&{0}\cr{3}&{0}&{-1}\cr
{6}&{-1}&{-1}\cr{207}&{-39}&{-23}\cr
{25}&{-5}&{-2}\cr{23}&{-5}&{0}\cr\endpmatrix
\hskip20pt
\pmatrix{-15}&{0}&{15}&{15}&{345}&{30}&{0}\cr
{0}&{-23}&{0}&{23}&{897}&{115}&{115}\cr
{15}&{0}&{-6}&{3}&{276}&{45}&{69}\cr
{15}&{23}&{3}&{-2}&{0}&{5}&{23}\cr
{345}&{897}&{276}&{0}&{-69}&{0}&{276}\cr
{30}&{115}&{45}&{5}&{0}&{-10}&{0}\cr
{0}&{115}&{69}&{23}&{276}&{0}&{-46}\cr\endpmatrix$.}

\vbox{\noindent
$N=102\ $ $d=357=3\cdot 7\cdot 17,\,\eta=3,\,h=1$:
$\langle 17 \rangle \oplus \langle -7 \rangle \oplus \langle -3 \rangle$
\nobreak
\newline
$\left(\smallmatrix
{13}&{55}&{7}&{0}&{0}&{7}&{3}&{105}&{1}&{13}&{7}&{6}&
{42}&{367}&{117}&{3465}\cr
{-20}&{-85}&{-11}&{0}&{1}&{0}&{-1}&{-51}&{-1}&{-17}&{-10}
&{-9}&{-64}&{-561}&{-179}&{-5304}\cr
{-5}&{-17}&{0}&{1}&{0}&{-17}&{-7}&{-238}&{-2}&{-17}&{-7}&
{-4}&{-21}&{-170}&{-53}&{-1547}\cr
\endsmallmatrix\right)$
\nobreak
\newline
$\pmatrix
A&B\\
B&A\\
\endpmatrix
\ \text{where\ \ }
\pmatrix
A\\
B\\
\endpmatrix =
\left(\smallmatrix
{-2}&{0}&{7}&{15}&{140}&{1292}&{418}&{12495}\cr
{0}&{-17}&{0}&{51}&{595}&{5678}&{1853}&{55692}\cr
{7}&{0}&{-14}&{0}&{77}&{833}&{280}&{8568}\cr
{15}&{51}&{0}&{-3}&{0}&{51}&{21}&{714}\cr
{140}&{595}&{77}&{0}&{-7}&{0}&{7}&{357}\cr
{1292}&{5678}&{833}&{51}&{0}&{-34}&{0}&{357}\cr
{418}&{1853}&{280}&{21}&{7}&{0}&{-1}&{0}\cr
{12495}&{55692}&{8568}&{714}&{357}&{357}&{0}&{-714}\cr
{51}&{238}&{42}&{6}&{7}&{17}&{2}&{0}\cr
{238}&{1173}&{238}&{51}&{119}&{680}&{187}&{4998}\cr
{42}&{238}&{63}&{21}&{70}&{476}&{140}&{3927}\cr
{6}&{51}&{21}&{12}&{63}&{510}&{159}&{4641}\cr
{7}&{119}&{70}&{63}&{448}&{3927}&{1253}&{37128}\cr
{17}&{680}&{476}&{510}&{3927}&{35003}&{11220}&{333438}\cr
{2}&{187}&{140}&{159}&{1253}&{11220}&{3601}&{107100}\cr
{0}&{4998}&{3927}&{4641}&{37128}&{333438}&{107100}&{3186939}\cr
\endsmallmatrix\right)$.}

\vbox{\noindent
$N=103\ $ $d=357=3\cdot 7\cdot 17,\,\eta=5,\,h=1$:
$\langle 1 \rangle \oplus \langle -119 \rangle \oplus \langle -3 \rangle$
\nobreak
\newline
$\left(\smallmatrix
{0}&{0}&{1}\cr{0}&{1}&{0}\cr{1}&{0}&{-1}\cr
{561}&{-42}&{-187}\cr{145}&{-11}&{-47}\cr{2975}&{-227}&{-952}\cr
{156}&{-12}&{-49}\cr{1547}&{-120}&{-476}\cr{12}&{-1}&{-3}\cr
{102}&{-9}&{-17}\cr{11}&{-1}&{-1}\cr{119}&{-11}&{0}\cr
\endsmallmatrix\right)
\left(\smallmatrix
{-3}&{0}&{3}&{561}&{141}&{2856}&{147}&{1428}&{9}&{51}&{3}&{0}\cr
{0}&{-119}&{0}&{4998}&{1309}&{27013}&{1428}&{14280}
&{119}&{1071}&{119}&{1309}\cr
{3}&{0}&{-2}&{0}&{4}&{119}&{9}&{119}&{3}&{51}&{8}&{119}\cr
{561}&{4998}&{0}&{-102}&{0}&{357}&
{51}&{1071}&{51}&{2703}&{612}&{11781}\cr
{141}&{1309}&{4}&{0}&{-1}&{0}&{3}&{119}&{8}&{612}
&{145}&{2856}\cr
{2856}&{27013}&{119}&{357}&{0}&{-238}&{0}&{1309}&
{119}&{11781}&{2856}&{56882}\cr
{147}&{1428}&{9}&{51}&{3}&{0}&{-3}&{0}&{3}&{561}&{141}
&{2856}\cr
{1428}&{14280}&{119}&{1071}&{119}&{1309}&{0}&{-119}&
{0}&{4998}&{1309}&{27013}\cr
{9}&{119}&{3}&{51}&{8}&{119}&{3}&{0}&{-2}&{0}&{4}&{119}\cr
{51}&{1071}&{51}&{2703}&{612}&{11781}&{561}&{4998}&{0}&{-102}&{0}&{357}\cr
{3}&{119}&{8}&{612}&{145}&{2856}&{141}&{1309}&{4}&{0}&{-1}&{0}\cr
{0}&{1309}&{119}&{11781}&{2856}&{56882}&{2856}
&{27013}&{119}&{357}&{0}&{-238}\cr
\endsmallmatrix\right)$.}

\vbox{\noindent
$N=104\ $ $d=385=5\cdot 7\cdot 11,\,\eta=6,\,h=1$:
$\langle 5 \rangle \oplus \langle -11 \rangle \oplus \langle -7 \rangle$
\nobreak
\newline
$\left(\smallmatrix{0}&{0}&{1}\cr{0}&{1}&{0}\cr{1}&{0}&{-1}\cr
{77}&{-28}&{-55}\cr{37}&{-15}&{-25}\cr{49}&{-21}&{-32}\cr
{88}&{-40}&{-55}\cr{14}&{-7}&{-8}\cr{22}&{-12}&{-11}\cr
{5}&{-3}&{-2}\cr{77}&{-49}&{-22}\cr{23}&{-15}&{-5}\cr
{21}&{-14}&{-3}\cr{22}&{-15}&{0}\cr\endsmallmatrix\right)
\hskip20pt
\left(\smallmatrix
{-7}&{0}&{7}&{385}&{175}&{224}&{385}&{56}&{77}&{14}&{154}&{35}&{21}&{0}\cr
{0}&{-11}&{0}&{308}&{165}&{231}&{440}&{77}&{132}
&{33}&{539}&{165}&{154}&{165}\cr
{7}&{0}&{-2}&{0}&{10}&{21}&{55}&{14}&{33}&{11}&
{231}&{80}&{84}&{110}\cr
{385}&{308}&{0}&{-154}&{0}&{77}&{385}&{154}&{539}
&{231}&{6083}&{2310}&{2618}&{3850}\cr
{175}&{165}&{10}&{0}&{-5}&{0}&{55}&{35}&{165}&{80}
&{2310}&{905}&{1050}&{1595}\cr
{224}&{231}&{21}&{77}&{0}&{-14}&{0}&{21}&{154}&{84}
&{2618}&{1050}&{1239}&{1925}\cr
{385}&{440}&{55}&{385}&{55}&{0}&{-55}&{0}&{165}&{110}
&{3850}&{1595}&{1925}&{3080}\cr
{56}&{77}&{14}&{154}&{35}&{21}&{0}&{-7}&{0}&{7}&{385}
&{175}&{224}&{385}\cr
{77}&{132}&{33}&{539}&{165}&{154}&{165}&{0}&{-11}&{0}
&{308}&{165}&{231}&{440}\cr
{14}&{33}&{11}&{231}&{80}&{84}&{110}&{7}&{0}&{-2}&{0}
&{10}&{21}&{55}\cr
{154}&{539}&{231}&{6083}&{2310}&{2618}&{3850}&{385}&
{308}&{0}&{-154}&{0}&{77}&{385}\cr
{35}&{165}&{80}&{2310}&{905}&{1050}&{1595}&{175}&{165}&
{10}&{0}&{-5}&{0}&{55}\cr
{21}&{154}&{84}&{2618}&{1050}&{1239}&{1925}&{224}&{231}
&{21}&{77}&{0}&{-14}&{0}\cr
{0}&{165}&{110}&{3850}&{1595}&{1925}&{3080}&{385}&{440}
&{55}&{385}&{55}&{0}&{-55}\cr\endsmallmatrix\right)$.}

\vbox{\noindent
$N=105\ $ $d=390=2\cdot 3\cdot 5\cdot 13,\,\eta=6,h=1$:
$\langle 10 \rangle \oplus \langle -26 \rangle \oplus \langle -6 \rangle
(1/2,1/2,0)$
\nobreak
\newline
$\left(\smallmatrix{45}&{{65}\over{2}}&{{9}\over{2}}&{{11}\over{2}}&
{{39}\over{2}}&{{5}\over{2}}&{78}&{8}&{0}&{0}&{3}&{32}&{273}&
{60}&{2652}&{{659}\over{2}}\cr
{-9}&{{-15}\over{2}}&{{-3}\over{2}}&
{{-5}\over{2}}&{{-21}\over{2}}&{{-3}\over{2}}&{-48}&{-5}&{0}&
{1}&{0}&{-5}&{-48}&{-11}&{-498}&{{-125}\over{2}}\cr
{-55}&{-39}&{-5}&{-5}&{-13}&{-1}&{-13}&{0}&{1}&{0}&{-4}&{-40}&
{-338}&{-74}&{-3263}&{-405}\cr\endsmallmatrix\right)$
\nobreak
\newline
$\pmatrix
A&B\\
B&A\\
\endpmatrix
\ \text{where\ \ }
\pmatrix
A\\
B\\
\endpmatrix =
\left(\smallmatrix
{-6}&{0}&{24}&{240}&{2028}&{444}&{19578}&{2430}\cr
{0}&{-26}&{0}&{130}&{1248}&{286}&{12948}&{1625}\cr
{24}&{0}&{-6}&{0}&{78}&{24}&{1248}&{165}\cr
{240}&{130}&{0}&{-10}&{0}&{10}&{780}&{115}\cr
{2028}&{1248}&{78}&{0}&{-78}&{0}&{1092}&{195}\cr
{444}&{286}&{24}&{10}&{0}&{-2}&{0}&{5}\cr
{19578}&{12948}&{1248}&{780}&{1092}&{0}&{-78}&{0}\cr
{2430}&{1625}&{165}&{115}&{195}&{5}&{0}&{-10}\cr
{330}&{234}&{30}&{30}&{78}&{6}&{78}&{0}\cr
{234}&{195}&{39}&{65}&{273}&{39}&{1248}&{130}\cr
{30}&{39}&{15}&{45}&{273}&{51}&{2028}&{240}\cr
{30}&{65}&{45}&{235}&{1755}&{365}&{15600}&{1910}\cr
{78}&{273}&{273}&{1755}&{13767}&{2925}&{126672}&{15600}\cr
{6}&{39}&{51}&{365}&{2925}&{627}&{27300}&{3370}\cr
{78}&{1248}&{2028}&{15600}&{126672}&{27300}&{1192542}&{147420}\cr
{0}&{130}&{240}&{1910}&{15600}&{3370}&{147420}&{18235}\cr
\endsmallmatrix\right)$.}

\vbox{\noindent
$N=106\ $ $d=399=3\cdot 7\cdot 19,\,\eta=4,\,h=0$:
$\langle 14 \rangle \oplus \langle -19 \rangle \oplus \langle -6 \rangle
(1/2,0,1/2)$
\nobreak
\newline
$\pmatrix{0}&{0}&{1}\cr{0}&{1}&{0}\cr
{{3}\over{2}}&{0}&{{-7}\over{2}}\cr
{{3}\over{2}}&{-1}&{{-3}\over{2}}\cr
{{5}\over{2}}&{-2}&{{-3}\over{2}}\cr
{{57}\over{2}}&{-24}&{{-19}\over{2}}\cr
{{7}\over{2}}&{-3}&{{-1}\over{2}}\cr
{57}&{-49}&{0}\cr\endpmatrix
\hskip20pt
\pmatrix{-6}&{0}&{21}&{9}&{9}&{57}&{3}&{0}\cr
{0}&{-19}&{0}&{19}&{38}&{456}&{57}&{931}\cr
{21}&{0}&{-42}&{0}&{21}&{399}&{63}&{1197}\cr
{9}&{19}&{0}&{-1}&{1}&{57}&{12}&{266}\cr
{9}&{38}&{21}&{1}&{-2}&{0}&{4}&{133}\cr
{57}&{456}&{399}&{57}&{0}&{-114}&{0}&{399}\cr
{3}&{57}&{63}&{12}&{4}&{0}&{-1}&{0}\cr
{0}&{931}&{1197}&{266}&{133}&{399}&{0}&{-133}\cr
\endpmatrix$.}

\vbox{\noindent
$N=107\ $ $d=429=3\cdot 11 \cdot 13,\,\eta=3,\,h=1$:
$\langle 1 \rangle \oplus \langle -39 \rangle \oplus \langle -11\rangle$
\nobreak
\newline
$\left(\smallmatrix
{0}&{0}&{1}\cr{0}&{1}&{0}\cr{3}&{0}&{-1}\cr{1287}&{-88}&{-351}\cr
{41}&{-3}&{-11}\cr{1573}&{-121}&{-416}\cr{126}&{-10}&{-33}\cr
{132}&{-11}&{-34}\cr{156}&{-14}&{-39}\cr{9}&{-1}&{-2}\cr
{429}&{-55}&{-78}\cr{7}&{-1}&{-1}\cr{143}&{-22}&{-13}\cr
{6}&{-1}&{0}\cr\endsmallmatrix\right)
\left(\smallmatrix
{-11}&{0}&{11}&{3861}&{121}&{4576}&{363}&{374}
&{429}&{22}&{858}&{11}&{143}&{0}\cr
{0}&{-39}&{0}&{3432}&{117}&{4719}&{390}
&{429}&{546}&{39}&{2145}&{39}&{858}&{39}\cr
{11}&{0}&{-2}&{0}&{2}&{143}&{15}&{22}&{39}&{5}
&{429}&{10}&{286}&{18}\cr
{3861}&{3432}&{0}&{-858}&{0}&{3003}&{429}&{858}
&{2145}&{429}&{62205}&{1716}&{58344}&{4290}\cr
{121}&{117}&{2}&{0}&{-1}&{0}&{3}&{11}&{39}&{10}
&{1716}&{49}&{1716}&{129}\cr
{4576}&{4719}&{143}&{3003}&{0}&{-286}&{0}&{143}
&{858}&{286}&{58344}&{1716}&{61633}&{4719}\cr
{363}&{390}&{15}&{429}&{3}&{0}&{-3}&{0}&{39}&{18}
&{4290}&{129}&{4719}&{366}\cr
{374}&{429}&{22}&{858}&{11}&{143}&{0}&{-11}&{0}&{11}
&{3861}&{121}&{4576}&{363}\cr
{429}&{546}&{39}&{2145}&{39}&{858}&{39}&{0}&{-39}&{0}
&{3432}&{117}&{4719}&{390}\cr
{22}&{39}&{5}&{429}&{10}&{286}&{18}&{11}&{0}&{-2}&{0}
&{2}&{143}&{15}\cr
{858}&{2145}&{429}&{62205}&{1716}&{58344}&{4290}&{3861}
&{3432}&{0}&{-858}&{0}&{3003}&{429}\cr
{11}&{39}&{10}&{1716}&{49}&{1716}&{129}&{121}&{117}&{2}
&{0}&{-1}&{0}&{3}\cr
{143}&{858}&{286}&{58344}&{1716}&{61633}&{4719}&{4576}
&{4719}&{143}&{3003}&{0}&{-286}&{0}\cr
{0}&{39}&{18}&{4290}&{129}&{4719}&{366}&{363}&{390}
&{15}&{429}&{3}&{0}&{-3}\cr\endsmallmatrix\right)$.}

\vbox{\noindent
$N=108\ $ $d=435=3\cdot 5\cdot 29,\,\eta=0,\,h=0$:
$\langle 3 \rangle \oplus \langle -29 \rangle \oplus \langle -5\rangle$
\nobreak
\newline
$\pmatrix{0}&{0}&{1}\cr{0}&{1}&{0}\cr{1}&{0}&{-1}\cr
{145}&{-30}&{-87}\cr{4}&{-1}&{-2}\cr{145}&{-40}&{-58}\cr
{10}&{-3}&{-3}\cr{3}&{-1}&{0}\cr\endpmatrix
\hskip20pt
\pmatrix
{-5}&{0}&{5}&{435}&{10}&{290}&{15}&{0}\cr
{0}&{-29}&{0}&{870}&{29}&{1160}&{87}&{29}\cr
{5}&{0}&{-2}&{0}&{2}&{145}&{15}&{9}\cr
{435}&{870}&{0}&{-870}&{0}&{3045}&{435}&{435}\cr
{10}&{29}&{2}&{0}&{-1}&{0}&{3}&{7}\cr
{290}&{1160}&{145}&{3045}&{0}&{-145}&{0}&{145}\cr
{15}&{87}&{15}&{435}&{3}&{0}&{-6}&{3}\cr
{0}&{29}&{9}&{435}&{7}&{145}&{3}&{-2}\cr
\endpmatrix$.}

\vbox{\noindent
$N=109\ $ $d=435=3\cdot 5\cdot 29,\,\eta=6,\,h=1$:
$\langle 435 \rangle \oplus \langle -1 \rangle\oplus \langle -1 \rangle $
\nobreak
\newline
$\left(\smallmatrix
{8}&{24}&{11}&{3}&{1}&{2}&{84}&{8}&{125}&{3}&{5}&{1}
&{0}&{0}&{1}&{119}&{13}&{223}\cr
{-25}&{-87}&{-45}&{-14}&{-6}&{-15}&{-725}&{-71}&{-1131}
&{-30}&{-58}&{-15}&{-1}&{1}&{0}&{-290}&{-34}&{-609}\cr
{-165}&{-493}&{-225}&{-61}&{-20}&{-39}&{-1595}&{-151}&{-2349}
&{-55}&{-87}&{-15}&{1}&{0}&{-21}&{-2465}&{-269}&{-4611}\cr
\endsmallmatrix\right)$
\nobreak
\newline
$\pmatrix
A&B\\
B&A\\
\endpmatrix
\ \text{where\ \ }
\pmatrix
A\\
B\\
\endpmatrix =
\left(\smallmatrix
{-10}&{0}&{30}&{25}&{30}&{150}&{11020}&{1150}&{19140}\cr
{0}&{-58}&{0}&{29}&{58}&{348}&{27550}&{2900}&{48546}\cr
{30}&{0}&{-15}&{0}&{15}&{120}&{10440}&{1110}&{18705}\cr
{25}&{29}&{0}&{-2}&{1}&{21}&{2175}&{235}&{4002}\cr
{30}&{58}&{15}&{1}&{-1}&{0}&{290}&{34}&{609}\cr
{150}&{348}&{120}&{21}&{0}&{-6}&{0}&{6}&{174}\cr
{11020}&{27550}&{10440}&{2175}&{290}&{0}&{-290}&{0}&{870}\cr
{1150}&{2900}&{1110}&{235}&{34}&{6}&{0}&{-2}&{0}\cr
{19140}&{48546}&{18705}&{4002}&{609}&{174}&{870}&{0}&{-87}\cr
{615}&{1595}&{630}&{140}&{25}&{15}&{145}&{5}&{0}\cr
{1595}&{4263}&{1740}&{406}&{87}&{87}&{1885}&{145}&{1914}\cr
{630}&{1740}&{735}&{180}&{45}&{60}&{1740}&{150}&{2175}\cr
{140}&{406}&{180}&{47}&{14}&{24}&{870}&{80}&{1218}\cr
{25}&{87}&{45}&{14}&{6}&{15}&{725}&{71}&{1131}\cr
{15}&{87}&{60}&{24}&{15}&{51}&{3045}&{309}&{5046}\cr
{145}&{1885}&{1740}&{870}&{725}&{3045}&{206335}&{21315}&{352350}\cr
{5}&{145}&{150}&{80}&{71}&{309}&{21315}&{2207}&{36540}\cr
{0}&{1914}&{2175}&{1218}&{1131}&{5046}&{352350}&{36540}&{605607}\cr
\endsmallmatrix\right)$.}

\vbox{\noindent
$N=110\ $ $d=455=5\cdot 7\cdot 13,\,\eta=5,\,h=1$:
$\langle 7 \rangle \oplus \langle -26 \rangle\oplus \langle -10 \rangle
(0,1/2,1/2)$
\nobreak
\newline
$\left(\smallmatrix{325}&{325}&{16}&{65}&{2}&{1}&{13}&{0}&{0}&{25}&
{9}&{260}&{48}&{199}&{9087}&{1025}\cr
{{-67}\over{2}}&{-35}&{-2}&{-10}&{{-1}\over{2}}&{{-1}\over{2}}
&{-7}&{0}&{1}&{0}&{{-1}\over{2}}&{{-45}\over{2}}&{{-9}\over{2}}
&{{-39}\over{2}}&{-903}&{{-205}\over{2}}\cr
{{-533}\over{2}}&{-266}&{-13}&{-52}&{{-3}\over{2}}&{{-1}\over{2}}
&{0}&{1}&{0}&{-21}&{{-15}\over{2}}&{{-429}\over{2}}&{{-79}\over{2}}
&{{-327}\over{2}}&{-7462}&{{-1683}\over{2}}\cr
\endsmallmatrix\right)$
\nobreak
\newline
$\pmatrix
A&B\\
B&A\\
\endpmatrix
\ \text{where\ \ }
\pmatrix
A\\
B\\
\endpmatrix =
\left(\smallmatrix
{-26}&{0}&{13}&{585}&{117}&{507}&{23478}&{2665}\cr
{0}&{-35}&{0}&{455}&{105}&{490}&{23205}&{2660}\cr
{13}&{0}&{-2}&{0}&{3}&{21}&{1092}&{130}\cr
{585}&{455}&{0}&{-65}&{0}&{65}&{4095}&{520}\cr
{117}&{105}&{3}&{0}&{-1}&{0}&{91}&{15}\cr
{507}&{490}&{21}&{65}&{0}&{-2}&{0}&{5}\cr
{23478}&{23205}&{1092}&{4095}&{91}&{0}&{-91}&{0}\cr
{2665}&{2660}&{130}&{520}&{15}&{5}&{0}&{-10}\cr
{871}&{910}&{52}&{260}&{13}&{13}&{182}&{0}\cr
{910}&{1015}&{70}&{455}&{35}&{70}&{2275}&{210}\cr
{52}&{70}&{7}&{65}&{7}&{19}&{728}&{75}\cr
{260}&{455}&{65}&{910}&{130}&{455}&{19565}&{2145}\cr
{13}&{35}&{7}&{130}&{21}&{80}&{3549}&{395}\cr
{13}&{70}&{19}&{455}&{80}&{322}&{14560}&{1635}\cr
{182}&{2275}&{728}&{19565}&{3549}&{14560}&{662571}&{74620}\cr
{0}&{210}&{75}&{2145}&{395}&{1635}&{74620}&{8415}\cr
\endsmallmatrix\right)$.}

\vbox{\noindent
$N=111\ $ $d=465=3\cdot 5\cdot 31,\,\eta=5,\,h=1$:
$\langle 5\rangle \oplus \langle -31 \rangle\oplus \langle -3 \rangle $
\nobreak
\newline
$\left(\smallmatrix{0}&{0}&{1}\cr{0}&{1}&{0}\cr{3}&{0}&{-5}\cr
{31}&{-8}&{-31}\cr{17}&{-5}&{-15}\cr{16}&{-5}&{-13}\cr
{403}&{-130}&{-310}\cr{18}&{-6}&{-13}\cr{93}&{-32}&{-62}\cr
{42}&{-15}&{-25}\cr{62}&{-23}&{-31}\cr{13}&{-5}&{-5}\cr
{5}&{-2}&{-1}\cr{62}&{-25}&{0}\cr\endsmallmatrix\right)
\left(\smallmatrix
{-3}&{0}&{15}&{93}&{45}&{39}&{930}&{39}&{186}&{75}&{93}&{15}&{3}&{0}\cr
{0}&{-31}&{0}&{248}&{155}&{155}&{4030}&{186}&{992}
&{465}&{713}&{155}&{62}&{775}\cr
{15}&{0}&{-30}&{0}&{30}&{45}&{1395}&{75}&{465}&{255}
&{465}&{120}&{60}&{930}\cr
{93}&{248}&{0}&{-62}&{0}&{31}&{1395}&{93}&{713}&{465}
&{1023}&{310}&{186}&{3410}\cr
{45}&{155}&{30}&{0}&{-5}&{0}&{155}&{15}&{155}&{120}&{310}
&{105}&{70}&{1395}\cr
{39}&{155}&{45}&{31}&{0}&{-2}&{0}&{3}&{62}&{60}&{186}&{70}
&{51}&{1085}\cr
{930}&{4030}&{1395}&{1395}&{155}&{0}&{-155}&{0}&{775}&{930}
&{3410}&{1395}&{1085}&{24180}\cr
{39}&{186}&{75}&{93}&{15}&{3}&{0}&{-3}&{0}&{15}&{93}&{45}&{39}
&{930}\cr
{186}&{992}&{465}&{713}&{155}&{62}&{775}&{0}&{-31}&{0}
&{248}&{155}&{155}&{4030}\cr
{75}&{465}&{255}&{465}&{120}&{60}&{930}&{15}&{0}&{-30}
&{0}&{30}&{45}&{1395}\cr
{93}&{713}&{465}&{1023}&{310}&{186}&{3410}&{93}&{248}
&{0}&{-62}&{0}&{31}&{1395}\cr
{15}&{155}&{120}&{310}&{105}&{70}&{1395}&{45}&{155}&{30}
&{0}&{-5}&{0}&{155}\cr
{3}&{62}&{60}&{186}&{70}&{51}&{1085}&{39}&{155}&{45}&{31}
&{0}&{-2}&{0}\cr
{0}&{775}&{930}&{3410}&{1395}&{1085}&{24180}&{930}&{4030}
&{1395}&{1395}&{155}&{0}&{-155}\cr\endsmallmatrix\right)$.}

\vbox{\noindent
$N=112\ $ $d=483=3\cdot 7\cdot 23,\,\eta =7,h=1$:
$\langle 6\rangle \oplus \langle -46 \rangle\oplus \langle -7 \rangle
(1/2,1/2,0)$
\nobreak
\newline
$\left(\smallmatrix{7}&{-2}&{-4}\cr{46}&{-14}&{-23}\cr
{{21}\over{2}}&{{-7}\over{2}}&{-4}\cr
{{23}\over{2}}&{{-9}\over{2}}&{0}\cr{0}&{0}&{1}\cr
{0}&{1}&{0}\cr{1}&{0}&{-1}\cr{{5}\over{2}}&{{-1}\over{2}}&{-2}\cr
{391}&{-101}&{-253}\cr{{511}\over{2}}&{{-133}\over{2}}&{-164}\cr
{{2599}\over{2}}&{{-681}\over{2}}&{-828}\cr{399}&{-105}&{-253}\cr
{{2185}\over{2}}&{{-577}\over{2}}&{-690}\cr
{{283}\over{2}}&{{-75}\over{2}}&{-89}\cr
\endsmallmatrix\right)
\left(\smallmatrix
{-2}&{0}&{7}&{69}&{28}&{92}&{14}&{3}&{46}&{21}&{69}&{14}&{23}&{1}\cr
{0}&{-23}&{0}&{276}&{161}&{644}&{115}&{46}&{2139}
&{1288}&{6072}&{1771}&{4646}&{575}\cr
{7}&{0}&{-14}&{0}&{28}&{161}&{35}&{21}&{1288}&{798}
&{3864}&{1148}&{3059}&{385}\cr
{69}&{276}&{0}&{-138}&{0}&{207}&{69}&{69}&{6072}
&{3864}&{19182}&{5796}&{15663}&{2001}\cr
{28}&{161}&{28}&{0}&{-7}&{0}&{7}&{14}&{1771}&{1148}
&{5796}&{1771}&{4830}&{623}\cr
{92}&{644}&{161}&{207}&{0}&{-46}&{0}&{23}&{4646}
&{3059}&{15663}&{4830}&{13271}&{1725}\cr
{14}&{115}&{35}&{69}&{7}&{0}&{-1}&{1}&{575}&{385}
&{2001}&{623}&{1725}&{226}\cr
{3}&{46}&{21}&{69}&{14}&{23}&{1}&{-2}&{0}&{7}&{69}
&{28}&{92}&{14}\cr
{46}&{2139}&{1288}&{6072}&{1771}&{4646}&{575}&{0}&{-23}
&{0}&{276}&{161}&{644}&{115}\cr
{21}&{1288}&{798}&{3864}&{1148}&{3059}&{385}&{7}&{0}
&{-14}&{0}&{28}&{161}&{35}\cr
{69}&{6072}&{3864}&{19182}&{5796}&{15663}&{2001}&{69}
&{276}&{0}&{-138}&{0}&{207}&{69}\cr
{14}&{1771}&{1148}&{5796}&{1771}&{4830}&{623}&{28}&{161}
&{28}&{0}&{-7}&{0}&{7}\cr
{23}&{4646}&{3059}&{15663}&{4830}&{13271}&{1725}&{92}&{644}
&{161}&{207}&{0}&{-46}&{0}\cr
{1}&{575}&{385}&{2001}&{623}&{1725}&{226}&{14}&{115}&{35}
&{69}&{7}&{0}&{-1}\cr\endsmallmatrix\right)$.}

\vbox{\noindent
$N=113\ $ $d=570=2\cdot 3\cdot 5\cdot 19,\,\eta=6,\,h=1$:
$\langle 38\rangle \oplus \langle -10 \rangle\oplus \langle -6 \rangle
(0,1/2,1/2)$
\nobreak
\newline
$\left(\smallmatrix
{0}&{0}&{15}&{1}&{6}&{5}&{42}&{7}&{9}&{15}&{213}&{7}&{15}&{5}
&{15}&{1}\cr
{0}&{1}&{0}&{{-1}\over{2}}&{-6}&{-6}&{-57}&{-10}
&{{-27}\over{2}}&{{-47}\over{2}}&{-342}&{{-23}\over{2}}&{{-51}\over{2}}
&{-9}&{{-57}\over{2}}&{-2}\cr
{1}&{0}&{-38}&{{-5}\over{2}}&{-13}&{-10}&{-76}&{-12}&{{-29}\over{2}}
&{{-45}\over{2}}&{-304}&{{-19}\over{2}}&{{-37}\over{2}}&
{-5}&{{-19}\over{2}}&{0}\cr\endsmallmatrix\right)$
\nobreak
\newline
$\left(\smallmatrix{-6}&{0}&{228}&{15}&{78}&{60}&{456}&{72}&
{87}&{135}&{1824}&{57}&{111}&{30}&{57}&{0}\cr
{0}&{-10}&{0}&{5}&{60}&{60}&{570}&{100}&{135}&{235}&{3420}
&{115}&{255}&{90}&{285}&{20}\cr
{228}&{0}&{-114}&{0}&{456}&{570}&{6612}&{1254}&{1824}&
{3420}&{52098}&{1824}&{4332}&{1710}&{6384}&{570}\cr
{15}&{5}&{0}&{-2}&{3}&{10}&{171}&{36}&{57}&{115}&{1824}
&{66}&{165}&{70}&{285}&{28}\cr
{78}&{60}&{456}&{3}&{-6}&{0}&{228}&{60}&{111}&{255}&{4332}
&{165}&{447}&{210}&{969}&{108}\cr
{60}&{60}&{570}&{10}&{0}&{-10}&{0}&{10}&{30}&{90}&{1710}&{70}
&{210}&{110}&{570}&{70}\cr
{456}&{570}&{6612}&{171}&{228}&{0}&{-114}&{0}&{57}&{285}&{6384}
&{285}&{969}&{570}&{3363}&{456}\cr
{72}&{100}&{1254}&{36}&{60}&{10}&{0}&{-2}&{0}&{20}&{570}&{28}
&{108}&{70}&{456}&{66}\cr
{87}&{135}&{1824}&{57}&{111}&{30}&{57}&{0}&{-6}&{0}&{228}
&{15}&{78}&{60}&{456}&{72}\cr
{135}&{235}&{3420}&{115}&{255}&{90}&{285}&{20}&{0}&{-10}&{0}
&{5}&{60}&{60}&{570}&{100}\cr
{1824}&{3420}&{52098}&{1824}&{4332}&{1710}&{6384}&{570}&{228}
&{0}&{-114}&{0}&{456}&{570}&{6612}&{1254}\cr
{57}&{115}&{1824}&{66}&{165}&{70}&{285}&{28}&{15}&{5}&{0}&{-2}
&{3}&{10}&{171}&{36}\cr
{111}&{255}&{4332}&{165}&{447}&{210}&{969}&{108}&{78}&{60}
&{456}&{3}&{-6}&{0}&{228}&{60}\cr
{30}&{90}&{1710}&{70}&{210}&{110}&{570}&{70}&{60}&{60}&{570}
&{10}&{0}&{-10}&{0}&{10}\cr
{57}&{285}&{6384}&{285}&{969}&{570}&{3363}&{456}&{456}&{570}
&{6612}&{171}&{228}&{0}&{-114}&{0}\cr
{0}&{20}&{570}&{28}&{108}&{70}&{456}&{66}&{72}&{100}&{1254}
&{36}&{60}&{10}&{0}&{-2}\cr\endsmallmatrix\right)$.}

\vbox{\noindent
$N=114\ $ $d=615=3\cdot 5\cdot 41,\,\eta =0,\,h=1$:
$\langle 3 \rangle \oplus \langle -41 \rangle\oplus \langle -5 \rangle$
\nobreak
\newline
$\left(\smallmatrix{0}&{0}&{1}&{16}&{82}&{9}&{410}&
{25}&{20}&{164}&{27}&{80}&{246}&{15}&{410}&{11}\cr
{0}&{1}&{0}&{-3}&{-17}&{-2}&{-95}&{-6}&{-5}&{-42}&{-7}
&{-21}&{-65}&{-4}&{-110}&{-3}\cr
{1}&{0}&{-1}&{-9}&{-41}&{-4}&{-164}&{-9}&{-6}&{-41}&{-6}
&{-15}&{-41}&{-2}&{-41}&{0}\cr\endsmallmatrix\right)$
\nobreak
\newline
$\left(\smallmatrix
{-5}&{0}&{5}&{45}&{205}&{20}&{820}&{45}&{30}&{205}
&{30}&{75}&{205}&{10}&{205}&{0}\cr
{0}&{-41}&{0}&{123}&{697}&{82}&{3895}&{246}&{205}
&{1722}&{287}&{861}&{2665}&{164}&{4510}&{123}\cr
{5}&{0}&{-2}&{3}&{41}&{7}&{410}&{30}&{30}&{287}&{51}&{165}
&{533}&{35}&{1025}&{33}\cr
{45}&{123}&{3}&{-6}&{0}&{6}&{615}&{57}&{75}&{861}&{165}&{582}
&{1968}&{138}&{4305}&{159}\cr
{205}&{697}&{41}&{0}&{-82}&{0}&{1025}&{123}&{205}&{2665}&{533}
&{1968}&{6806}&{492}&{15785}&{615}\cr
{20}&{82}&{7}&{6}&{0}&{-1}&{0}&{3}&{10}&{164}&{35}&{138}&{492}
&{37}&{1230}&{51}\cr
{820}&{3895}&{410}&{615}&{1025}&{0}&{-205}&{0}&{205}&{4510}
&{1025}&{4305}&{15785}&{1230}&{42230}&{1845}\cr
{45}&{246}&{30}&{57}&{123}&{3}&{0}&{-6}&{0}&{123}&{33}&{159}
&{615}&{51}&{1845}&{87}\cr
{30}&{205}&{30}&{75}&{205}&{10}&{205}&{0}&{-5}&{0}&{5}&{45}
&{205}&{20}&{820}&{45}\cr
{205}&{1722}&{287}&{861}&{2665}&{164}&{4510}&{123}&{0}&{-41}
&{0}&{123}&{697}&{82}&{3895}&{246}\cr
{30}&{287}&{51}&{165}&{533}&{35}&{1025}&{33}&{5}&{0}&{-2}&{3}&{41}
&{7}&{410}&{30}\cr
{75}&{861}&{165}&{582}&{1968}&{138}&{4305}&{159}&{45}&{123}&{3}
&{-6}&{0}&{6}&{615}&{57}\cr
{205}&{2665}&{533}&{1968}&{6806}&{492}&{15785}&{615}&{205}&{697}
&{41}&{0}&{-82}&{0}&{1025}&{123}\cr
{10}&{164}&{35}&{138}&{492}&{37}&{1230}&{51}&{20}&{82}&{7}&{6}
&{0}&{-1}&{0}&{3}\cr
{205}&{4510}&{1025}&{4305}&{15785}&{1230}&{42230}&{1845}&{820}
&{3895}&{410}&{615}&{1025}&{0}&{-205}&{0}\cr
{0}&{123}&{33}&{159}&{615}&{51}&{1845}&{87}&{45}&{246}&{30}&{57}
&{123}&{3}&{0}&{-6}\cr\endsmallmatrix\right)$.}

\vbox{\noindent
$N=115\ $ $d=645=3\cdot 5\cdot 43,\,\eta=3,\,h=1$:
$\langle 10 \rangle \oplus \langle -86 \rangle\oplus \langle -3 \rangle
(1/2,1/2,0)$
\nobreak
\newline
$\left(\smallmatrix{{165}\over{2}}&{{645}\over{2}}&{{43}\over{2}}&
{516}&{22}&{4}&{{387}\over{2}}&{{3}\over{2}}&{{43}\over{2}}
&{0}&{0}&{1}&{129}&{{31}\over{2}}&{11}&{1419}&{36}&{1591}\cr
{{-33}\over{2}}&{{-131}\over{2}}&{{-9}\over{2}}&{-111}&{-5}
&{-1}&{{-105}\over{2}}&{{-1}\over{2}}&{{-15}\over{2}}&{0}&{1}&
{0}&{-18}&{{-5}\over{2}}&{-2}&{-270}&{-7}&{-315}\cr
{-122}&{-473}&{-31}&{-731}&{-30}&{-5}&{-215}&{-1}&{0}&{1}&{0}&
{-2}&{-215}&{-25}&{-17}&{-2150}&{-54}&{-2365}\cr
\endsmallmatrix\right)$
\nobreak
\newline
$\pmatrix
A&B\\
B&A\\
\endpmatrix
\ \text{where\ \ }
\pmatrix
A\\
B\\
\endpmatrix =
\left(\smallmatrix
{-3}&{0}&{6}&{645}&{75}&{51}&{6450}&{162}&{7095}\cr
{0}&{-86}&{0}&{1548}&{215}&{172}&{23220}&{602}&{27090}\cr
{6}&{0}&{-2}&{0}&{5}&{8}&{1290}&{36}&{1720}\cr
{645}&{1548}&{0}&{-129}&{0}&{129}&{25800}&{774}&{39345}\cr
{75}&{215}&{5}&{0}&{-10}&{0}&{645}&{25}&{1505}\cr
{51}&{172}&{8}&{129}&{0}&{-1}&{0}&{2}&{215}\cr
{6450}&{23220}&{1290}&{25800}&{645}&{0}&{-1290}&{0}&{7740}\cr
{162}&{602}&{36}&{774}&{25}&{2}&{0}&{-2}&{0}\cr
{7095}&{27090}&{1720}&{39345}&{1505}&{215}&{7740}&{0}&{-215}\cr
{366}&{1419}&{93}&{2193}&{90}&{15}&{645}&{3}&{0}\cr
{1419}&{5633}&{387}&{9546}&{430}&{86}&{4515}&{43}&{645}\cr
{93}&{387}&{29}&{774}&{40}&{10}&{645}&{9}&{215}\cr
{2193}&{9546}&{774}&{22317}&{1290}&{387}&{29670}&{516}&{16125}\cr
{90}&{430}&{40}&{1290}&{85}&{30}&{2580}&{50}&{1720}\cr
{15}&{86}&{10}&{387}&{30}&{13}&{1290}&{28}&{1075}\cr
{645}&{4515}&{645}&{29670}&{2580}&{1290}&{139965}&{3225}&{130935}\cr
{3}&{43}&{9}&{516}&{50}&{28}&{3225}&{77}&{3225}\cr
{0}&{645}&{215}&{16125}&{1720}&{1075}&{130935}&{3225}&{138890}\cr
\endsmallmatrix\right)$.}

\vbox{\noindent
$N=116\ $ $d=651=3\cdot 7\cdot 31,\,\eta=4,\,h=1$:
$\langle 14 \rangle \oplus \langle -31 \rangle\oplus \langle -6 \rangle
(1/2,0,1/2)$
\nobreak
\newline
$\left(\smallmatrix{129}&{992}&{13}&{{3}\over{2}}&{31}&{0}&
{0}&{{3}\over{2}}&{{31}\over{2}}&{{21}\over{2}}&{{899}\over{2}}&
{{67}\over{2}}&{{183}\over{2}}&{{25885}\over{2}}&{{2511}\over{2}}&
{5766}&{2091}&{1426}\cr
{-78}&{-602}&{-8}&{-1}&{-21}&{0}&{1}&{0}&{-8}&{-6}&{-266}&{-20}&
{-55}&{-7791}&{-756}&{-3473}&{-1260}&{-860}\cr
{-86}&{-651}&{-8}&{{-1}\over{2}}&{0}&{1}&{0}&{{-7}\over{2}}&
{{-31}\over{2}}&{{-17}\over{2}}&{{-651}\over{2}}&{{-47}\over{2}}&
{{-125}\over{2}}&{{-17577}\over{2}}&{{-1703}\over{2}}&{-3906}&
{-1414}&{-961}\cr\endsmallmatrix\right)$
\nobreak
\newline
$\pmatrix
A&B\\
B&A\\
\endpmatrix
\ \text{where\ \ }
\pmatrix
A\\
B\\
\endpmatrix =
\left(\smallmatrix
{-6}&{0}&{6}&{33}&{5208}&{516}&{2418}&{903}&{651}\cr
{0}&{-434}&{0}&{217}&{38626}&{3906}&{18662}&{7161}&{5425}\cr
{6}&{0}&{-2}&{1}&{434}&{48}&{248}&{105}&{93}\cr
{33}&{217}&{1}&{-1}&{0}&{3}&{31}&{21}&{31}\cr
{5208}&{38626}&{434}&{0}&{-217}&{0}&{651}&{651}&{1519}\cr
{516}&{3906}&{48}&{3}&{0}&{-6}&{0}&{21}&{93}\cr
{2418}&{18662}&{248}&{31}&{651}&{0}&{-31}&{0}&{248}\cr
{903}&{7161}&{105}&{21}&{651}&{21}&{0}&{-42}&{0}\cr
{651}&{5425}&{93}&{31}&{1519}&{93}&{248}&{0}&{-62}\cr
{69}&{651}&{15}&{9}&{651}&{51}&{186}&{42}&{0}\cr
{651}&{7161}&{217}&{217}&{21917}&{1953}&{8246}&{2604}&{1302}\cr
{15}&{217}&{9}&{13}&{1519}&{141}&{620}&{210}&{124}\cr
{9}&{217}&{13}&{29}&{3906}&{375}&{1705}&{609}&{403}\cr
{651}&{21917}&{1519}&{3906}&{545104}&{52731}&{241521}&{87234}&{59024}\cr
{51}&{1953}&{141}&{375}&{52731}&{5109}&{23436}&{8484}&{5766}\cr
{186}&{8246}&{620}&{1705}&{241521}&{23436}&{107663}&{39060}&{26660}\cr
{42}&{2604}&{210}&{609}&{87234}&{8484}&{39060}&{14217}&{9765}\cr
{0}&{1302}&{124}&{403}&{59024}&{5766}&{26660}&{9765}&{6789}\cr
\endsmallmatrix\right)$.}

\vbox{\noindent
$N=117\ $ $d=795=3\cdot 5\cdot 53,\,\eta =6,\,h=1$:
$\langle 3 \rangle \oplus \langle -106 \rangle\oplus \langle -10 \rangle
(0,1/2,1/2)$
\nobreak
\newline
$\left(\smallmatrix
{143}&{101}&{1060}&{4}&{265}&{3}&{53}&{0}&{0}&{5}&{7}&{49}&{1590}
&{46}&{23585}&{747}&{31747}&{2200}&{5300}&{895}\cr
{{-29}\over{2}}&{{-21}\over{2}}&{{-225}\over{2}}
&{{-1}\over{2}}&{{-75}\over{2}}&{{-1}\over{2}}
&{-9}&{0}&{1}&{0}&{{-1}\over{2}}&{{-9}\over{2}}
&{{-305}\over{2}}&{{-9}\over{2}}&{{-4695}\over{2}}
&{{-149}\over{2}}&{-3171}&{-220}&{-531}&{-90}\cr
{{-125}\over{2}}&{{-87}\over{2}}&{{-901}\over{2}}&{{-3}\over{2}}
&{{-159}\over{2}}&{{-1}\over{2}}&{0}&{1}&{0}&{-3}&{{-7}\over{2}}
&{{-45}\over{2}}&{{-1431}\over{2}}&{{-41}\over{2}}&{{-20829}\over{2}}
&{{-659}\over{2}}&{-13992}&{-969}&{-2332}&{-393}\cr
\endsmallmatrix\right)$
\nobreak
\newline
$\pmatrix
A&B\\
B&A\\
\endpmatrix
\,\text{where}\,
\pmatrix
A\\
B\\
\endpmatrix =
\left(\smallmatrix
{-2}&{3}&{265}&{10}&{6360}&{206}&{8904}&{625}&{1537}&{270}\cr
{3}&{-6}&{0}&{3}&{3975}&{135}&{6042}&{435}&{1113}&{210}\cr
{265}&{0}&{-265}&{0}&{37365}&{1325}&{61215}&{4505}&{11925}&{2385}\cr
{10}&{3}&{0}&{-1}&{0}&{2}&{159}&{15}&{53}&{15}\cr
{6360}&{3975}&{37365}&{0}&{-1590}&{0}&{6360}&{795}&{3975}&{1590}\cr
{206}&{135}&{1325}&{2}&{0}&{-2}&{0}&{5}&{53}&{30}\cr
{8904}&{6042}&{61215}&{159}&{6360}&{0}&{-159}&{0}&{954}&{795}\cr
{625}&{435}&{4505}&{15}&{795}&{5}&{0}&{-10}&{0}&{30}\cr
{1537}&{1113}&{11925}&{53}&{3975}&{53}&{954}&{0}&{-106}&{0}\cr
{270}&{210}&{2385}&{15}&{1590}&{30}&{795}&{30}&{0}&{-15}\cr
{47}&{42}&{530}&{5}&{795}&{19}&{636}&{35}&{53}&{0}\cr
{42}&{51}&{795}&{12}&{3180}&{90}&{3498}&{225}&{477}&{60}\cr
{530}&{795}&{14310}&{265}&{89040}&{2650}&{107325}&{7155}&{16165}&{2385}\cr
{5}&{12}&{265}&{6}&{2385}&{73}&{3021}&{205}&{477}&{75}\cr
{795}&{3180}&{89040}&{2385}&{1139235}&{35775}&{1510500}
&{104145}&{248835}&{41340}\cr
{19}&{90}&{2650}&{73}&{35775}&{1127}&{47700}&{3295}&{7897}&{1320}\cr
{636}&{3498}&{107325}&{3021}&{1510500}&{47700}&{2022639}
&{139920}&{336126}&{56445}\cr
{35}&{225}&{7155}&{205}&{104145}&{3295}&{139920}&{9690}
&{23320}&{3930}\cr
{53}&{477}&{16165}&{477}&{248835}&{7897}&{336126}&{23320}&{56286}&{9540}\cr
{0}&{60}&{2385}&{75}&{41340}&{1320}&{56445}&{3930}&{9540}&{1635}\cr
\endsmallmatrix\right)$.}

\vbox{\noindent
$N=118\ $ $d=1155=3\cdot 5\cdot 7\cdot 11,\,\eta=2,h=0$:
$\langle 14 \rangle \oplus \langle -22 \rangle\oplus \langle -15 \rangle
(1/2,1/2,0)$
\nobreak
\newline
$\pmatrix
{0}&{0}&{1}&{{165}\over{2}}&{{9}\over{2}}&{{55}\over{2}}
&{{39}\over{2}}&{{25}\over{2}}&{132}&{{1}\over{2}}\cr
{0}&{1}&{0}&{{-35}\over{2}}&{{-3}\over{2}}&{{-25}\over{2}}
&{{-21}\over{2}}&{{-15}\over{2}}&{-84}&{{-1}\over{2}}\cr
{1}&{0}&{-1}&{-77}&{-4}&{-22}&{-14}&{-8}&{-77}&{0}\cr
\endpmatrix$
\nobreak
\newline
$\pmatrix
{-15}&{0}&{15}&{1155}&{60}&{330}&{210}&{120}&{1155}&{0}\cr
{0}&{-22}&{0}&{385}&{33}&{275}&{231}&{165}&{1848}&{11}\cr
{15}&{0}&{-1}&{0}&{3}&{55}&{63}&{55}&{693}&{7}\cr
{1155}&{385}&{0}&{-385}&{0}&{1540}&{2310}&{2310}&{31185}&{385}\cr
{60}&{33}&{3}&{0}&{-6}&{0}&{42}&{60}&{924}&{15}\cr
{330}&{275}&{55}&{1540}&{0}&{-110}&{0}&{110}&{2310}&{55}\cr
{210}&{231}&{63}&{2310}&{42}&{0}&{-42}&{0}&{462}&{21}\cr
{120}&{165}&{55}&{2310}&{60}&{110}&{0}&{-10}&{0}&{5}\cr
{1155}&{1848}&{693}&{31185}&{924}&{2310}&{462}&{0}&{-231}&{0}\cr
{0}&{11}&{7}&{385}&{15}&{55}&{21}&{5}&{0}&{-2}\cr
\endpmatrix$.}

\vbox{\noindent
$N=119\ $ $d=1155=3\cdot 5\cdot 7\cdot 11,\,\eta=8,\,h=1$:
$\langle 1155 \rangle \oplus \langle -1 \rangle\oplus \langle -1\rangle$
\nobreak
\newline
$\left(\smallmatrix{0}&{0}&{1}&{4}&{2}&{17}&{43}&{5}&{61}&{1}&
{1}&{6}&{7}&{2}&{11}&{23}&{2}&{16}\cr
{-1}&{1}&{0}&{-33}&{-20}&{-189}&{-495}&{-60}&{-770}&{-14}&{-16}&{-105}&
{-132}&{-40}&{-231}&{-495}&{-45}&{-385}\cr
{1}&{0}&{-35}&{-132}&{-65}&{-546}&{-1375}&{-159}&{-1925}&{-31}&{-30}&
{-175}&{-198}&{-55}&{-294}&{-605}&{-51}&{-385}\cr
\endsmallmatrix\right)$
\nobreak
\newline
$\left(\smallmatrix
{-2}&{1}&{35}&{99}&{45}&{357}&{880}&{99}&{1155}&{17}&{14}&{70}
&{66}&{15}&{63}&{110}&{6}&{0}\cr
{1}&{-1}&{0}&{33}&{20}&{189}&{495}&{60}&{770}&{14}&{16}&{105}&{132}
&{40}&{231}&{495}&{45}&{385}\cr
{35}&{0}&{-70}&{0}&{35}&{525}&{1540}&{210}&{3080}&{70}&{105}&{805}
&{1155}&{385}&{2415}&{5390}&{525}&{5005}\cr
{99}&{33}&{0}&{-33}&{0}&{231}&{825}&{132}&{2310}&{66}&{132}&{1155}
&{1848}&{660}&{4389}&{10065}&{1023}&{10395}\cr
{45}&{20}&{35}&{0}&{-5}&{0}&{55}&{15}&{385}&{15}&{40}&{385}&{660}
&{245}&{1680}&{3905}&{405}&{4235}\cr
{357}&{189}&{525}&{231}&{0}&{-42}&{0}&{21}&{1155}&{63}&{231}&{2415}
&{4389}&{1680}&{11802}&{27720}&{2919}&{31185}\cr
{880}&{495}&{1540}&{825}&{55}&{0}&{-55}&{0}&{1540}&{110}&{495}&{5390}
&{10065}&{3905}&{27720}&{65395}&{6930}&{74690}\cr
{99}&{60}&{210}&{132}&{15}&{21}&{0}&{-6}&{0}&{6}&{45}&{525}&{1023}&{405}
&{2919}&{6930}&{741}&{8085}\cr
{1155}&{770}&{3080}&{2310}&{385}&{1155}&{1540}&{0}&{-770}&{0}&{385}
&{5005}&{10395}&{4235}&{31185}&{74690}&{8085}&{89705}\cr
{17}&{14}&{70}&{66}&{15}&{63}&{110}&{6}&{0}&{-2}&{1}&{35}&{99}&{45}
&{357}&{880}&{99}&{1155}\cr
{14}&{16}&{105}&{132}&{40}&{231}&{495}&{45}&{385}&{1}&{-1}&{0}&{33}
&{20}&{189}&{495}&{60}&{770}\cr
{70}&{105}&{805}&{1155}&{385}&{2415}&{5390}&{525}&{5005}&{35}&{0}
&{-70}&{0}&{35}&{525}&{1540}&{210}&{3080}\cr
{66}&{132}&{1155}&{1848}&{660}&{4389}&{10065}&{1023}&{10395}&{99}
&{33}&{0}&{-33}&{0}&{231}&{825}&{132}&{2310}\cr
{15}&{40}&{385}&{660}&{245}&{1680}&{3905}&{405}&{4235}&{45}&{20}&{35}
&{0}&{-5}&{0}&{55}&{15}&{385}\cr
{63}&{231}&{2415}&{4389}&{1680}&{11802}&{27720}&{2919}&{31185}&{357}
&{189}&{525}&{231}&{0}&{-42}&{0}&{21}&{1155}\cr
{110}&{495}&{5390}&{10065}&{3905}&{27720}&{65395}&{6930}&{74690}&{880}
&{495}&{1540}&{825}&{55}&{0}&{-55}&{0}&{1540}\cr
{6}&{45}&{525}&{1023}&{405}&{2919}&{6930}&{741}&{8085}&{99}&{60}&{210}
&{132}&{15}&{21}&{0}&{-6}&{0}\cr
{0}&{385}&{5005}&{10395}&{4235}&{31185}&{74690}&{8085}&{89705}&{1155}
&{770}&{3080}&{2310}&{385}&{1155}&{1540}&{0}&{-770}\cr
\endsmallmatrix\right)$.}

\vbox{\noindent
$N=120\ $ $d=1155=3\cdot 5\cdot 7\cdot 11,\,\eta=14,\,h=0$:
$\langle 22 \rangle \oplus \langle -15 \rangle\oplus \langle -14\rangle
(1/2,0,1/2)$
\nobreak
\newline
$\pmatrix{0}&{0}&{{7}\over{2}}&{{5}\over{2}}&{7}&{{15}\over{2}}
&{{7}\over{2}}&{{525}\over{2}}&{{5}\over{2}}&{9}\cr
{0}&{1}&{0}&{-2}&{-7}&{-8}&{-4}&{-308}&{-3}&{-11}\cr
{1}&{0}&{{-11}\over{2}}&{{-5}\over{2}}&{-5}&{{-9}\over{2}}
&{{-3}\over{2}}&{{-165}\over{2}}&{{-1}\over{2}}&{0}\cr
\endpmatrix$
\nobreak
\newline
$\pmatrix
{-14}&{0}&{77}&{35}&{70}&{63}&{21}&{1155}&{7}&{0}\cr
{0}&{-15}&{0}&{30}&{105}&{120}&{60}&{4620}&{45}&{165}\cr
{77}&{0}&{-154}&{0}&{154}&{231}&{154}&{13860}&{154}&{693}\cr
{35}&{30}&{0}&{-10}&{0}&{15}&{20}&{2310}&{30}&{165}\cr
{70}&{105}&{154}&{0}&{-7}&{0}&{14}&{2310}&{35}&{231}\cr
{63}&{120}&{231}&{15}&{0}&{-6}&{3}&{1155}&{21}&{165}\cr
{21}&{60}&{154}&{20}&{14}&{3}&{-2}&{0}&{2}&{33}\cr
{1155}&{4620}&{13860}&{2310}&{2310}&{1155}&{0}&{-2310}&{0}&{1155}\cr
{7}&{45}&{154}&{30}&{35}&{21}&{2}&{0}&{-1}&{0}\cr
{0}&{165}&{693}&{165}&{231}&{165}&{33}&{1155}&{0}&{-33}\cr
\endpmatrix$.}

\vbox{\noindent
$N=121\ $ $d=1365=3\cdot 5\cdot 7\cdot 13,\,\eta=5,\,h=1$:
$\langle 5 \rangle \oplus \langle -39 \rangle\oplus \langle -7\rangle$
\nobreak
\newline
$\left(\smallmatrix{90}&{49}&{33}&{3}&{13}&{0}&{0}&{1}
&{182}&{9}&{28}&{132}&{63}&{1417}&{462}&{1287}&{142}&{4823}\cr
{-25}&{-14}&{-10}&{-1}&{-5}&{0}&{1}&{0}&{-35}&{-2}&{-7}&{-35}
&{-17}&{-385}&{-126}&{-352}&{-39}&{-1330}\cr
{-48}&{-25}&{-15}&{-1}&{0}&{1}&{0}&{-1}&{-130}&{-6}&{-17}&{-75}
&{-35}&{-780}&{-253}&{-702}&{-77}&{-2600}\cr\endsmallmatrix\right)$
\nobreak
\newline
$\pmatrix
A&B\\
B&A\\
\endpmatrix
\ \text{where\ \ }
\pmatrix
A\\
B\\
\endpmatrix =
\left(\smallmatrix
{-3}&{0}&{60}&{39}&{975}&{336}&{975}&{114}&{4095}\cr
{0}&{-14}&{0}&{14}&{455}&{175}&{546}&{70}&{2730}\cr
{60}&{0}&{-30}&{0}&{195}&{105}&{390}&{60}&{2730}\cr
{39}&{14}&{0}&{-1}&{0}&{7}&{39}&{8}&{455}\cr
{975}&{455}&{195}&{0}&{-130}&{0}&{195}&{65}&{5005}\cr
{336}&{175}&{105}&{7}&{0}&{-7}&{0}&{7}&{910}\cr
{975}&{546}&{390}&{39}&{195}&{0}&{-39}&{0}&{1365}\cr
{114}&{70}&{60}&{8}&{65}&{7}&{0}&{-2}&{0}\cr
{4095}&{2730}&{2730}&{455}&{5005}&{910}&{1365}&{0}&{-455}\cr
{84}&{63}&{75}&{15}&{195}&{42}&{78}&{3}&{0}\cr
{63}&{63}&{105}&{28}&{455}&{119}&{273}&{21}&{455}\cr
{75}&{105}&{255}&{90}&{1755}&{525}&{1365}&{135}&{4095}\cr
{15}&{28}&{90}&{37}&{780}&{245}&{663}&{70}&{2275}\cr
{195}&{455}&{1755}&{780}&{17030}&{5460}&{15015}&{1625}&{54145}\cr
{42}&{119}&{525}&{245}&{5460}&{1771}&{4914}&{539}&{18200}\cr
{78}&{273}&{1365}&{663}&{15015}&{4914}&{13728}&{1521}&{51870}\cr
{3}&{21}&{135}&{70}&{1625}&{539}&{1521}&{171}&{5915}\cr
{0}&{455}&{4095}&{2275}&{54145}&{18200}&{51870}&{5915}&{207480}\cr
\endsmallmatrix\right)$.}

\vbox{\noindent
$N=122\ $ $d=1365=3\cdot 5\cdot 7\cdot 13,\,\eta=15,\,h=1$:
$\langle 6 \rangle \oplus \langle -65 \rangle\oplus \langle -14\rangle
(1/2,0,1/2)$
\nobreak
\newline
$\left(\smallmatrix
{{223}\over{2}}&{{5005}\over{2}}&{37}&{140}&{13}&{0}&{0}&
{{1}\over{2}}&{91}&{{15}\over{2}}&{{49}\over{2}}&{31}&{1820}
&{58}&{{1715}\over{2}}&{{1209}\over{2}}\cr
{-33}&{-742}&{-11}&{-42}&{-4}&{0}&{1}&{0}&{-21}&{-2}&{-7}&{-9}&{-532}
&{-17}&{-252}&{-178}\cr
{{-33}\over{2}}&{{-715}\over{2}}&{-5}&{-15}&{0}&{1}&{0}&{{-1}\over{2}}
&{-39}&{{-5}\over{2}}&{{-11}\over{2}}&{-6}&{-325}&{-10}&
{{-285}\over{2}}&{{-195}\over{2}}\cr\endsmallmatrix\right.$\hfill
\nobreak
\newline
\line{\hfil$\left.\smallmatrix
{{1995}\over{2}}&{8645}&{{1709}\over{2}}&{25844}&{515}&{308}\cr
{-294}&{-2549}&{-252}&{-7623}&{-152}&{-91}\cr
{{-317}\over{2}}&{-1365}&{{-269}\over{2}}&{-4056}&{-80}&{-47}\cr
\endsmallmatrix\right)$}
\nobreak
\newline
$\pmatrix
A&B\\
B&A\\
\endpmatrix$\  where
\nobreak
\newline
$\pmatrix
A\\
B\\
\endpmatrix =
\left(\smallmatrix
{-3}&{0}&{3}&{105}&{117}&{231}&{2145}&{219}&{6825}&{150}&{105}\cr
{0}&{-910}&{0}&{1365}&{2275}&{5005}&{48230}&{5005}&{158340}&{3640}&{2730}\cr
{3}&{0}&{-1}&{0}&{26}&{70}&{715}&{76}&{2457}&{60}&{49}\cr
{105}&{1365}&{0}&{-210}&{0}&{210}&{2730}&{315}&{10920}&{315}&{315}\cr
{117}&{2275}&{26}&{0}&{-26}&{0}&{260}&{39}&{1638}&{65}&{91}\cr
{231}&{5005}&{70}&{210}&{0}&{-14}&{0}&{7}&{546}&{35}&{77}\cr
{2145}&{48230}&{715}&{2730}&{260}&{0}&{-65}&{0}&{1365}&{130}&{455}\cr
{219}&{5005}&{76}&{315}&{39}&{7}&{0}&{-2}&{0}&{5}&{35}\cr
{6825}&{158340}&{2457}&{10920}&{1638}&{546}&{1365}&{0}&{-273}&{0}&{819}\cr
{150}&{3640}&{60}&{315}&{65}&{35}&{130}&{5}&{0}&{-10}&{0}\cr
{105}&{2730}&{49}&{315}&{91}&{77}&{455}&{35}&{819}&{0}&{-7}\cr
{48}&{1365}&{27}&{210}&{78}&{84}&{585}&{51}&{1365}&{15}&{0}\cr
{1365}&{42315}&{910}&{8190}&{3640}&{4550}&{34580}&{3185}&{90090}&{1365}&
{455}\cr
{27}&{910}&{21}&{210}&{104}&{140}&{1105}&{104}&{3003}&{50}&{21}\cr
{210}&{8190}&{210}&{2415}&{1365}&{1995}&{16380}&{1575}&{46410}&{840}&
{420}\cr
{78}&{3640}&{104}&{1365}&{871}&{1365}&{11570}&{1131}&{33852}&{650}&{364}\cr
{84}&{4550}&{140}&{1995}&{1365}&{2219}&{19110}&{1883}&{56784}&{1120}&
{658}\cr
{585}&{34580}&{1105}&{16380}&{11570}&{19110}&{165685}&{16380}&{495495}&
{9880}&{5915}\cr
{51}&{3185}&{104}&{1575}&{1131}&{1883}&{16380}&{1622}&{49140}&
{985}&{595}\cr
{1365}&{90090}&{3003}&{46410}&{33852}&{56784}&{495495}&{49140}&{1490853}&
{30030}&{18291}\cr
{15}&{1365}&{50}&{840}&{650}&{1120}&{9880}&{985}&{30030}&{615}&{385}\cr
{0}&{455}&{21}&{420}&{364}&{658}&{5915}&{595}&{18291}&{385}&{252}\cr
\endsmallmatrix\right)$.}

\newpage

\centerline{\bf Table 2.}
\centerline{The list of odd hyperbolic lattices of rank three}
\centerline{and with even square-free determinant (i.e. non-main)}
\centerline{which are reflective of elliptic or parabolic type.}

\vskip30pt

\noindent
$N^\prime=1\ $
$d=2,\,\text{odd},\,\eta=0,\,h=0$:
$\langle 1 \rangle \oplus \langle -1 \rangle \oplus \langle -2 \rangle$.
It is equivariantly equivalent to
$d=1,\,\eta=0,\,h=0$: $U\oplus \langle -1 \rangle$.

\vskip3pt

\noindent
$N^\prime =2\ $ $d=6=2\cdot 3,\,\text{odd},\,\eta=0,\,h=0$:
$\langle 1 \rangle \oplus \langle -1\rangle \oplus \langle -6 \rangle$.
It is equivariantly equivalent to
$d=3,\,\eta=1,\,h=0$: $U\oplus \langle -3 \rangle$.

\vskip3pt

\vbox{\noindent
$N^\prime=3\ $ $d=6=2\cdot 3,\,\text{odd},\,\eta=1,\,h=0$:
$\langle 6 \rangle \oplus \langle -1 \rangle \oplus \langle -1 \rangle$
($=\widetilde{(3,0,h=0)}$)
\nobreak
\newline
$\pmatrix{0}&{-1}&{1}\cr{0}&{1}&{0}\cr
{1}&{0}&{-3}\cr{1}&{-2}&{-2}\cr\endpmatrix
\hskip20pt
\pmatrix{-2}&{1}&{3}&{0}\cr{1}&{-1}&{0}&{2}\cr
{3}&{0}&{-3}&{0}\cr{0}&{2}&{0}&{-2}\cr\endpmatrix$.}

\vskip3pt

\vbox{\noindent
$N^\prime =4\ $
$d=10=2\cdot 5,\,\text{odd},\,\eta=0,\,h=1$:
$\langle 2 \rangle \oplus \langle -5 \rangle \oplus \langle -1\rangle$
(=$\widetilde{(5,1,h=0)}$)
\nobreak
\newline
$\pmatrix{0}&{0}&{1}\cr{0}&{1}&{0}\cr{1}&{0}&{-2}\cr
{2}&{-1}&{-2}\cr{10}&{-6}&{-5}\cr{3}&{-2}&{0}\cr\endpmatrix
\hskip20pt
\pmatrix
{-1}&{0}&{2}&{2}&{5}&{0}\cr
{0}&{-5}&{0}&{5}&{30}&{10}\cr
{2}&{0}&{-2}&{0}&{10}&{6}\cr
{2}&{5}&{0}&{-1}&{0}&{2}\cr
{5}&{30}&{10}&{0}&{-5}&{0}\cr
{0}&{10}&{6}&{2}&{0}&{-2}\cr\endpmatrix$.}

\vskip3pt

\noindent
$N^\prime =5\ $ $d=10=2\cdot 5,\,\text{odd},\,\eta=1,\,h=0$:
$\langle 1\rangle \oplus \langle-1\rangle \oplus \langle-10\rangle$.
It is equivariantly equivalent to
$d=5,\,\eta=0,\,h=0$:  $U\oplus \langle-5\rangle$.

\vskip3pt

\vbox{\noindent
$N^\prime =6\ $ $d=14=2 \cdot 7,\,\text{odd},\,\eta=0,\,h=1$:
$\langle 14 \rangle \oplus \langle -1 \rangle \oplus \langle -1\rangle$
(=$\widetilde{(7,0,h=1)}$)
\nobreak
\newline
$\pmatrix
{0}&{-1}&{1}\cr{0}&{1}&{0}\cr{1}&{0}&{-4}\cr
{2}&{-3}&{-7}\cr{4}&{-9}&{-12}\cr{3}&{-8}&{-8}\cr\endpmatrix
\hskip20pt
\pmatrix
{-2}&{1}&{4}&{4}&{3}&{0}\cr{1}&{-1}&{0}&{3}&{9}&{8}\cr
{4}&{0}&{-2}&{0}&{8}&{10}\cr{4}&{3}&{0}&{-2}&{1}&{4}\cr
{3}&{9}&{8}&{1}&{-1}&{0}\cr{0}&{8}&{10}&{4}&{0}&{-2}\cr
\endpmatrix$.}

\vskip3pt

\noindent
$N^\prime =7\ $ $d=14=2\cdot 7,\text{odd},\,\eta=1,\,h=0$:
$\langle 1 \rangle \oplus \langle -1 \rangle \oplus \langle -14 \rangle$.
It is equivariantly equivalent to
$d=7,\,\eta=1,\,h=0$: $U\oplus \langle -7\rangle$.

\vskip3pt

\noindent
$N^\prime =8\ $ $d=22=2 \cdot 11,\,\text{odd},\,\eta=0,\,h=1$:
$\langle 1\rangle \oplus \langle -1\rangle \oplus \langle -22\rangle$.
It is equivariantly equivalent to
$d=11,\,\eta=1,\,h=1$: $U\oplus \langle -11\rangle$.

\vskip3pt

\vbox{\noindent
$N^\prime =9\ $ $d=22=2\cdot 11,\,\text{odd},\,\eta=1,\,h=1$:
$\langle 22\rangle \oplus \langle -1\rangle \oplus \langle -1\rangle$
(=$\widetilde{(11,0,h=0)}$)
\nobreak
\newline
$\pmatrix{0}&{-1}&{1}\cr{0}&{1}&{0}\cr{7}&{0}&{-33}\cr
{3}&{-2}&{-14}\cr{2}&{-3}&{-9}\cr{2}&{-5}&{-8}\cr
{15}&{-44}&{-55}\cr{3}&{-10}&{-10}\cr\endpmatrix
\hskip20pt
\pmatrix{-2}&{1}&{33}&{12}&{6}&{3}&{11}&{0}\cr
{1}&{-1}&{0}&{2}&{3}&{5}&{44}&{10}\cr
{33}&{0}&{-11}&{0}&{11}&{44}&{495}&{132}\cr
{12}&{2}&{0}&{-2}&{0}&{10}&{132}&{38}\cr
{6}&{3}&{11}&{0}&{-2}&{1}&{33}&{12}\cr
{3}&{5}&{44}&{10}&{1}&{-1}&{0}&{2}\cr
{11}&{44}&{495}&{132}&{33}&{0}&{-11}&{0}\cr
{0}&{10}&{132}&{38}&{12}&{2}&{0}&{-2}\cr\endpmatrix$.}

\vskip3pt

\noindent
$N^\prime =10\ $ $d=26=2\cdot 13,\,\text{odd},\,\eta=1,\,h=1$:
$\langle 1\rangle \oplus \langle -1\rangle \oplus \langle -26\rangle$.
It is equivariantly equivalent to
$d=13,\,\eta=0,\,h=1$: $U\oplus\langle-13\rangle$.

\vskip3pt

\vbox{\noindent
$N^\prime =11\ $ $d=30=2\cdot 3\cdot 5,\,\text{odd},\eta=0,\,h=0$:
$\langle 30 \rangle \oplus  \langle -1\rangle \oplus \langle -1\rangle$
(=$\widetilde{(15,3,h=0)}$)
\nobreak
\newline
$\pmatrix{0}&{-1}&{1}\cr{0}&{1}&{0}\cr{1}&{0}&{-6}\cr
{2}&{-5}&{-10}\cr{1}&{-4}&{-4}\cr\endpmatrix
\hskip20pt
\pmatrix
{-2}&{1}&{6}&{5}&{0}\cr{1}&{-1}&{0}&{5}&{4}\cr
{6}&{0}&{-6}&{0}&{6}\cr{5}&{5}&{0}&{-5}&{0}\cr
{0}&{4}&{6}&{0}&{-2}\cr\endpmatrix$.}

\vskip3pt

\noindent
$N^\prime =12\ $ $d=30=2\cdot 3\cdot 5,\,\text{odd},\,\eta=1,\,h=0$:
$\langle 1 \rangle\oplus \langle -1\rangle \oplus \langle -30\rangle$.
It is equivariantly equivalent to
$d=15,\,\eta=2,\,h=0$: $U\oplus \langle -15\rangle$.

\vskip3pt

\noindent
$N^\prime =13\ $ $d=30=2\cdot 3\cdot 5,\,\text{odd},\,\eta=2,\,h=0$:
$\langle 10\rangle \oplus \langle -6\rangle \oplus \langle -2 \rangle
(1/2,1/2,0)$. It is equivariantly equivalent to
$d=15,\,\eta=1,\,h=0$: $\langle 5 \rangle \oplus \langle -3\rangle
\oplus \langle -1\rangle$.

\vskip3pt

\vbox{\noindent
$N^\prime =14\ $ $d=30=2\cdot 3\cdot 5,\,\text{odd},\,\eta=3,\,h=0$:
$\langle 6 \rangle \oplus \langle -10\rangle \oplus \langle -2\rangle
(1/2,1/2,0)$
($=\widetilde{(15,0,h=0)}$)
\nobreak
\newline
$\pmatrix{0}&{0}&{1}\cr{0}&{1}&{0}\cr{1}&{0}&{-2}\cr
{10}&{-4}&{-15}\cr{{5}\over{2}}&{{-3}\over{2}}&{-3}\cr
{{1}\over{2}}&{{-1}\over{2}}&{0}\cr\endpmatrix
\hskip20pt
\pmatrix
{-2}&{0}&{4}&{30}&{6}&{0}\cr
{0}&{-10}&{0}&{40}&{15}&{5}\cr
{4}&{0}&{-2}&{0}&{3}&{3}\cr
{30}&{40}&{0}&{-10}&{0}&{10}\cr
{6}&{15}&{3}&{0}&{-3}&{0}\cr
{0}&{5}&{3}&{10}&{0}&{-1}\cr\endpmatrix$.}

\vskip3pt

\noindent
$N^\prime =15\ $ $d=34=2\cdot 17,\,\text{odd},\,\eta=0,\,h=0$:
$\langle 1\rangle \oplus \langle-1 \rangle \oplus \langle -34\rangle$.
It is equivariantly equivalent to
$d=17,\,\eta=0,\,h=0$: $U\oplus \langle -17\rangle$.

\vskip3pt

\noindent
$N^\prime =16\ $ $d=42=2\cdot 3\cdot 7,\,\text{odd},\,\eta=0,\,h=1$:
$\langle 1 \rangle\oplus \langle -1 \rangle \oplus \langle -42\rangle$.
It is equivariantly equivalent to
$d=21,\,\eta=1,\,h=1$: $U\oplus \langle -21\rangle$.

\vskip3pt

\vbox{\noindent
$N^\prime =17\ $ $d=42=2\cdot 3\cdot 7,\,\text{odd},\,\eta=1,\,h=1$:
$\langle 7 \rangle \oplus \langle -3 \rangle \oplus \langle -2 \rangle$
(=$\widetilde{(21,0,h=0)}$)
\nobreak\newline
$\pmatrix{0}&{0}&{1}\cr{0}&{1}&{0}\cr{1}&{0}&{-2}\cr
{12}&{-7}&{-21}\cr{2}&{-2}&{-3}\cr{3}&{-4}&{-3}\cr
{2}&{-3}&{-1}\cr{9}&{-14}&{0}\cr\endpmatrix
\hskip20pt
\pmatrix
{-2}&{0}&{4}&{42}&{6}&{6}&{2}&{0}\cr
{0}&{-3}&{0}&{21}&{6}&{12}&{9}&{42}\cr
{4}&{0}&{-1}&{0}&{2}&{9}&{10}&{63}\cr
{42}&{21}&{0}&{-21}&{0}&{42}&{63}&{462}\cr
{6}&{6}&{2}&{0}&{-2}&{0}&{4}&{42}\cr
{6}&{12}&{9}&{42}&{0}&{-3}&{0}&{21}\cr
{2}&{9}&{10}&{63}&{4}&{0}&{-1}&{0}\cr
{0}&{42}&{63}&{462}&{42}&{21}&{0}&{-21}\cr
\endpmatrix$.}

\vskip3pt

\vbox{\noindent
$N^\prime =18\ $ $d=42=2\cdot 3\cdot 7,\,\text{odd},\,\eta=2,h=0$:
$\langle 2\rangle \oplus \langle -14 \rangle \oplus \langle -6 \rangle
(0,1/2,1/2)$
(=$\widetilde{(21,3,h=0)}$)
\nobreak\newline
$\pmatrix{0}&{0}&{1}\cr{0}&{1}&{0}\cr{3}&{0}&{-2}\cr
{7}&{{-3}\over{2}}&{{-7}\over{2}}\cr{3}&{-1}&{-1}\cr
{4}&{{-3}\over{2}}&{{-1}\over{2}}\cr{21}&{-8}&{0}\cr\endpmatrix
\hskip20pt
\pmatrix
{-6}&{0}&{12}&{21}&{6}&{3}&{0}\cr
{0}&{-14}&{0}&{21}&{14}&{21}&{112}\cr
{12}&{0}&{-6}&{0}&{6}&{18}&{126}\cr
{21}&{21}&{0}&{-7}&{0}&{14}&{126}\cr
{6}&{14}&{6}&{0}&{-2}&{0}&{14}\cr
{3}&{21}&{18}&{14}&{0}&{-1}&{0}\cr
{0}&{112}&{126}&{126}&{14}&{0}&{-14}\cr
\endpmatrix$.}

\vskip3pt

\noindent
$N^\prime =19\ $ $d=42=2\cdot 3\cdot 7,\,\text{odd},\,\eta=3,\,h=0$:
$\langle 6\rangle \oplus \langle -14\rangle \oplus \langle -2 \rangle
(1/2,0,1/2)$. It is equivariantly equivalent to
$d=21,\,\eta=2,\,h=0$:
$\langle 3 \rangle \oplus \langle -7\rangle \oplus \langle -1 \rangle$.

\vskip3pt

\vbox{\noindent
$N^\prime =20\ $ $d=66=2 \cdot 3\cdot 11,\,\text{odd},\eta=0,\,h=1$:
$\langle 2 \rangle \oplus \langle -22\rangle \oplus \langle -6 \rangle
(1/2,1/2,0)$
(=$\widetilde{(33,3,h=1)}$)
\nobreak\newline
$\pmatrix
{0}&{0}&{1}\cr{0}&{1}&{0}\cr{3}&{0}&{-2}\cr{66}&{-10}&{-33}\cr
{{17}\over{2}}&{{-3}\over{2}}&{-4}\cr{30}&{-6}&{-13}\cr
{55}&{-12}&{-22}\cr{12}&{-3}&{-4}\cr{44}&{-12}&{-11}\cr
{{3}\over{2}}&{{-1}\over{2}}&{0}\cr\endpmatrix
\pmatrix
{-6}&{0}&{12}&{198}&{24}&{78}&{132}&{24}&{66}&{0}\cr
{0}&{-22}&{0}&{220}&{33}&{132}&{264}&{66}&{264}&{11}\cr
{12}&{0}&{-6}&{0}&{3}&{24}&{66}&{24}&{132}&{9}\cr
{198}&{220}&{0}&{-22}&{0}&{66}&{264}&{132}&{990}&{88}\cr
{24}&{33}&{3}&{0}&{-1}&{0}&{11}&{9}&{88}&{9}\cr
{78}&{132}&{24}&{66}&{0}&{-6}&{0}&{12}&{198}&{24}\cr
{132}&{264}&{66}&{264}&{11}&{0}&{-22}&{0}&{220}&{33}\cr
{24}&{66}&{24}&{132}&{9}&{12}&{0}&{-6}&{0}&{3}\cr
{66}&{264}&{132}&{990}&{88}&{198}&{220}&{0}&{-22}&{0}\cr
{0}&{11}&{9}&{88}&{9}&{24}&{33}&{3}&{0}&{-1}\cr
\endpmatrix$.}

\vskip3pt

\noindent
$N^\prime =21\ $ $d=66=2\cdot 3\cdot 11,\,\text{odd},\eta=1,\,h=1$:
$\langle 1 \rangle \oplus \langle -1\rangle \oplus \langle -66\rangle$.
It is equivariantly equivalent to
$d=33,\,\eta=2,\,h=1$: $U\oplus \langle -33\rangle$.

\vskip3pt

\noindent
$N^\prime =22\ $ $d=66=2\cdot 3\cdot 11,\,\text{odd},\,\eta=2,\,h=0$:
$\langle 22 \rangle \oplus \langle -6\rangle \oplus \langle -2\rangle
(1/2,0,1/2)$.
It is equivariantly equivalent to
$d=33,\,\eta=1,\,h=0$:
$\langle 11\rangle \oplus \langle -3\rangle \oplus \langle -1\rangle$.

\vskip3pt

\vbox{\noindent
$N^\prime =23\ $ $d=70=2\cdot 5\cdot 7,\,\text{odd},\,\eta=1,\,h=1$:
$\langle 14 \rangle \oplus \langle -10\rangle \oplus \langle -2\rangle
(1/2,1/2,0)$
(=$\widetilde{(35,0,h=0)}$)
\nobreak\newline
$\left(\smallmatrix
{50}&{-11}&{-130}\cr{7}&{-2}&{-18}\cr{10}&{-4}&{-25}\cr
{{1}\over{2}}&{{-1}\over{2}}&{-1}\cr
{{5}\over{2}}&{{-7}\over{2}}&{0}\cr{0}&{0}&{1}\cr{0}&{1}&{0}\cr
{3}&{0}&{-8}\cr{40}&{-6}&{-105}\cr{{39}\over{2}}&{{-7}\over{2}}&{-51}\cr
{{1395}\over{2}}&{{-273}\over{2}}&{-1820}\cr{130}&{-26}&{-339}\cr
\endsmallmatrix\right)
\hskip20pt
\left(\smallmatrix
{-10}&{0}&{60}&{35}&{1365}&{260}&{110}&{20}&{40}&{5}&{35}&{0}\cr
{0}&{-2}&{0}&{3}&{175}&{36}&{20}&{6}&{20}&{5}&{105}&{16}\cr
{60}&{0}&{-10}&{0}&{210}&{50}&{40}&{20}&{110}&{40}&{1190}&{210}\cr
{35}&{3}&{0}&{-1}&{0}&{2}&{5}&{5}&{40}&{17}&{560}&{102}\cr
{1365}&{175}&{210}&{0}&{-35}&{0}&{35}&{105}&{1190}&{560}&{19635}&{3640}\cr
{260}&{36}&{50}&{2}&{0}&{-2}&{0}&{16}&{210}&{102}&{3640}&{678}\cr
{110}&{20}&{40}&{5}&{35}&{0}&{-10}&{0}&{60}&{35}&{1365}&{260}\cr
{20}&{6}&{20}&{5}&{105}&{16}&{0}&{-2}&{0}&{3}&{175}&{36}\cr
{40}&{20}&{110}&{40}&{1190}&{210}&{60}&{0}&{-10}&{0}&{210}&{50}\cr
{5}&{5}&{40}&{17}&{560}&{102}&{35}&{3}&{0}&{-1}&{0}&{2}\cr
{35}&{105}&{1190}&{560}&{19635}&{3640}&{1365}&{175}&{210}&{0}&{-35}&{0}\cr
{0}&{16}&{210}&{102}&{3640}&{678}&{260}&{36}&{50}&{2}&{0}&{-2}\cr
\endsmallmatrix\right)$.}

\vskip3pt

\vbox{\noindent
$N^\prime =24\ $ $d=70=2\cdot 5\cdot 7,\,\text{odd},\,\eta=2,\,h=1$:
$\langle 70 \rangle \oplus \langle -1\rangle \oplus \langle -1\rangle$
(=$\widetilde{(35,3,h=0)}$)
\nobreak\newline
$\left(\smallmatrix{0}&{-1}&{1}\cr{0}&{1}&{0}\cr{5}&{0}&{-42}\cr
{6}&{-5}&{-50}\cr{11}&{-14}&{-91}\cr{9}&{-14}&{-74}\cr
{6}&{-11}&{-49}\cr{2}&{-5}&{-16}\cr{9}&{-28}&{-70}\cr
{4}&{-15}&{-30}\cr{3}&{-14}&{-21}\cr{1}&{-6}&{-6}\cr
\endsmallmatrix\right)
\hskip20pt
\left(\smallmatrix
{-2}&{1}&{42}&{45}&{77}&{60}&{38}&{11}&{42}&{15}&{7}&{0}\cr
{1}&{-1}&{0}&{5}&{14}&{14}&{11}&{5}&{28}&{15}&{14}&{6}\cr
{42}&{0}&{-14}&{0}&{28}&{42}&{42}&{28}&{210}&{140}&{168}&{98}\cr
{45}&{5}&{0}&{-5}&{0}&{10}&{15}&{15}&{140}&{105}&{140}&{90}\cr
{77}&{14}&{28}&{0}&{-7}&{0}&{7}&{14}&{168}&{140}&{203}&{140}\cr
{60}&{14}&{42}&{10}&{0}&{-2}&{0}&{6}&{98}&{90}&{140}&{102}\cr
{38}&{11}&{42}&{15}&{7}&{0}&{-2}&{1}&{42}&{45}&{77}&{60}\cr
{11}&{5}&{28}&{15}&{14}&{6}&{1}&{-1}&{0}&{5}&{14}&{14}\cr
{42}&{28}&{210}&{140}&{168}&{98}&{42}&{0}&{-14}&{0}&{28}&{42}\cr
{15}&{15}&{140}&{105}&{140}&{90}&{45}&{5}&{0}&{-5}&{0}&{10}\cr
{7}&{14}&{168}&{140}&{203}&{140}&{77}&{14}&{28}&{0}&{-7}&{0}\cr
{0}&{6}&{98}&{90}&{140}&{102}&{60}&{14}&{42}&{10}&{0}&{-2}\cr
\endsmallmatrix\right)$.}

\vskip3pt

\noindent
$N^\prime =25\ $ $d=70=2\cdot 5\cdot 7,\,\text{odd},\,\eta=3,\,h=1$:
$\langle 2 \rangle \oplus \langle -14\rangle \oplus
\langle -10 \rangle (1/2,1/2,0)$. It is equivariantly equivalent to
$d=35,\,\eta=2,\,h=1$: $\langle 1\rangle \oplus
\langle -7\rangle\oplus \langle -5\rangle$.

\vskip3pt

\noindent
$N^\prime =26\ $ $d=78=2\cdot 3\cdot 13,\,\text{odd},\,\eta=1,\,h=1$:
$\langle 6\rangle \oplus \langle -13\rangle \oplus \langle -1\rangle$.
It is equivariantly equivalent to
$d=39,\,\eta=2,\,h=1$:
$\langle 3 \rangle \oplus \langle -26\rangle \oplus
\langle -2\rangle(0,1/2,1/2)$.

\vskip3pt

\vbox{\noindent
\noindent
$N^\prime =27\ $ $d=78=2\cdot 3\cdot 13,\,\text{odd},\eta=3,\,h=0$:
$\langle 6 \rangle\oplus \langle -26\rangle \oplus \langle -2 \rangle
(0,1/2,1/2)$
(=$\widetilde{(39,0,h=0)}$)
\nobreak\newline
$\pmatrix{0}&{0}&{1}\cr{0}&{1}&{0}\cr{1}&{0}&{-2}\cr
{4}&{-1}&{-6}\cr{39}&{-12}&{-52}\cr
{4}&{{-3}\over{2}}&{{-9}\over{2}}\cr
{1}&{{-1}\over{2}}&{{-1}\over{2}}\cr\endpmatrix
\hskip20pt
\pmatrix
{-2}&{0}&{4}&{12}&{104}&{9}&{1}\cr
{0}&{-26}&{0}&{26}&{312}&{39}&{13}\cr
{4}&{0}&{-2}&{0}&{26}&{6}&{4}\cr
{12}&{26}&{0}&{-2}&{0}&{3}&{5}\cr
{104}&{312}&{26}&{0}&{-26}&{0}&{26}\cr
{9}&{39}&{6}&{3}&{0}&{-3}&{0}\cr
{1}&{13}&{4}&{5}&{26}&{0}&{-1}\cr
\endpmatrix$.}

\vskip3pt

\noindent
$N^\prime =28\ $ $d=102=2\cdot 3\cdot 17,\,\text{odd},\,\eta=0,\,h=1$:
$\langle 34\rangle \oplus \langle -6\rangle \oplus \langle -2\rangle
(1/2,1/2,0)$. It is equivariantly equivalent to
$d=51,\,\eta=1,\,h=1$:
$\langle 17\rangle \oplus \langle -3\rangle \oplus \langle -1\rangle$.

\vskip3pt

\noindent
$N^\prime =29\ $ $d=114=2\cdot 3\cdot 19,\,\eta=1,\,h=1$:
$\langle 6 \rangle \oplus \langle -38 \rangle \oplus
\langle -2\rangle (1/2,0,1/2)$. It is equivariantly equivalent to
$d=57,\,\eta=2,\,h=1$:
$\langle 3 \rangle \oplus \langle -19\rangle \oplus
\langle -1\rangle$.

\vskip3pt

\vbox{\noindent
$N^\prime =30\ $ $d=138=2\cdot 3\cdot 23,\,\text{odd},\eta=2,\,h=1$:
$\langle 2\rangle \oplus \langle -46\rangle \oplus \langle -6\rangle
(0,1/2,1/2)$
(=$\widetilde{(69,3,h=0)}$)
\nobreak\newline
$\left(\smallmatrix
{72}&{276}&{69}&{115}&{5}&{138}&{12}&{115}&{0}&{0}&
{3}&{23}&{7}&{1104}&{126}&{1541}\cr
{-12}&{-47}&{-12}&{{-41}\over{2}}&{-1}&{{-57}\over{2}}&{{-5}\over{2}}
&{-24}&{0}&{1}&{0}&{{-5}\over{2}}&{-1}&{{-357}\over{2}}
&{{-41}\over{2}}&{-252}\cr
{-25}&{-92}&{-22}&{{-69}\over{2}}&{-1}&{{-23}\over{2}}&{{-1}\over{2}}
&{0}&{1}&{0}&{-2}&{{-23}\over{2}}&{-3}&{{-805}\over{2}}&{{-91}\over{2}}
&{-552}\cr\endsmallmatrix\right)$
\nobreak\newline
$\left(\smallmatrix
{-6}&{0}&{12}&{69}&{18}&{2415}&{273}&{3312}&{150}&{552}
&{132}&{207}&{6}&{69}&{3}&{0}\cr
{0}&{-46}&{0}&{115}&{46}&{8211}&{943}&{11592}&{552}&{2162}
&{552}&{943}&{46}&{1311}&{115}&{1104}\cr
{12}&{0}&{-6}&{0}&{6}&{1794}&{210}&{2622}&{132}&{552}&{150}
&{276}&{18}&{690}&{66}&{690}\cr
{69}&{115}&{0}&{-23}&{0}&{2484}&{299}&{3818}&{207}&{943}&{276}
&{552}&{46}&{2277}&{230}&{2530}\cr
{18}&{46}&{6}&{0}&{-2}&{0}&{2}&{46}&{6}&{46}&{18}&{46}&{6}&{414}
&{44}&{506}\cr
{2415}&{8211}&{1794}&{2484}&{0}&{-69}&{0}&{276}&{69}&{1311}&{690}
&{2277}&{414}&{42918}&{4761}&{56856}\cr
{273}&{943}&{210}&{299}&{2}&{0}&{-1}&{0}&{3}&{115}&{66}&{230}&{44}
&{4761}&{530}&{6348}\cr
{3312}&{11592}&{2622}&{3818}&{46}&{276}&{0}&{-46}&{0}&{1104}&{690}
&{2530}&{506}&{56856}&{6348}&{76222}\cr
{150}&{552}&{132}&{207}&{6}&{69}&{3}&{0}&{-6}&{0}&{12}&{69}&{18}
&{2415}&{273}&{3312}\cr
{552}&{2162}&{552}&{943}&{46}&{1311}&{115}&{1104}&{0}&{-46}&{0}
&{115}&{46}&{8211}&{943}&{11592}\cr
{132}&{552}&{150}&{276}&{18}&{690}&{66}&{690}&{12}&{0}&{-6}&{0}
&{6}&{1794}&{210}&{2622}\cr
{207}&{943}&{276}&{552}&{46}&{2277}&{230}&{2530}&{69}&{115}&{0}&
{-23}&{0}&{2484}&{299}&{3818}\cr
{6}&{46}&{18}&{46}&{6}&{414}&{44}&{506}&{18}&{46}&{6}&{0}&{-2}&{0}
&{2}&{46}\cr
{69}&{1311}&{690}&{2277}&{414}&{42918}&{4761}&{56856}&{2415}&{8211}
&{1794}&{2484}&{0}&{-69}&{0}&{276}\cr
{3}&{115}&{66}&{230}&{44}&{4761}&{530}&{6348}&{273}&{943}&{210}&{299}
&{2}&{0}&{-1}&{0}\cr
{0}&{1104}&{690}&{2530}&{506}&{56856}&{6348}&{76222}&{3312}&{11592}
&{2622}&{3818}&{46}&{276}&{0}&{-46}\cr
\endsmallmatrix\right)$.}

\vskip3pt

\vbox{\noindent
$N^\prime =31\ $ $d=210=2\cdot 3\cdot 5\cdot 7,\,\text{odd},
\eta=0,\,h=0$:
$\langle 7\rangle \oplus \langle -5\rangle \oplus\langle -6\rangle$
(=$\widetilde{(105,3,h=0)}$)
\nobreak\newline
$\pmatrix{0}&{0}&{1}\cr{0}&{1}&{0}\cr{6}&{0}&{-7}\cr
{2}&{-1}&{-2}\cr{40}&{-28}&{-35}\cr{3}&{-3}&{-2}\cr{5}&{-7}&{0}\cr
\endpmatrix
\hskip20pt
\pmatrix
{-6}&{0}&{42}&{12}&{210}&{12}&{0}\cr
{0}&{-5}&{0}&{5}&{140}&{15}&{35}\cr
{42}&{0}&{-42}&{0}&{210}&{42}&{210}\cr
{12}&{5}&{0}&{-1}&{0}&{3}&{35}\cr
{210}&{140}&{210}&{0}&{-70}&{0}&{420}\cr
{12}&{15}&{42}&{3}&{0}&{-6}&{0}\cr
{0}&{35}&{210}&{35}&{420}&{0}&{-70}\cr
\endpmatrix$.}

\vskip3pt

\noindent
$N^\prime =32\ $ $d=210=2\cdot 3\cdot 5\cdot 7,\text{odd},\eta=1,\,h=1$:
$\langle 14 \rangle \oplus \langle -30\rangle \oplus\langle -2 \rangle
(1/2,0,1/2)$. It is equivariantly equivalent to
$d=105,\eta=2,\,h=1$:
$\langle 7 \rangle \oplus \langle -15\rangle \oplus
\langle -1 \rangle$.

\vskip3pt

\noindent
$N^\prime =33\ $ $d=210=2\cdot 3\cdot 5\cdot 7,\,\text{odd},\,
\eta=2,\,h=0$:
$\langle 2 \rangle \oplus \langle -42 \rangle \oplus \langle -10\rangle
(0,1/2,1/2)$. It is equivariantly equivalent to
$d=105=3\cdot 5\cdot 7,\,\eta=1,\,h=0$:
$\langle 1 \rangle \oplus \langle -21\rangle \oplus \langle -5\rangle$.

\vskip3pt

\noindent
$N^\prime =34\ $ $d=210=2\cdot 3\cdot 5\cdot 7,\,\text{odd},\,\eta=4,\,h=0$:
$\langle 5 \rangle \oplus \langle -7\rangle \oplus
\langle -6\rangle$. It
is equivariantly equivalent to
$d=105,\,\eta=7,\,h=0$:
$\langle 10\rangle \oplus \langle -14\rangle \oplus
\langle -3\rangle(1/2,1/2,0)$.

\vskip3pt

\vbox{\noindent
$N^\prime =35\ $ $d=210=2.3.5.7,\,\text{odd},\,\eta=6,\,h=1$:
$\langle 10\rangle\oplus \langle -14 \rangle\oplus
\langle -6 \rangle(0,1/2,1/2)$ r
(=$\widetilde{(105,5,h=0)}$)
\nobreak\newline
$\left(\smallmatrix
{0}&{0}&{3}&{2}&{13}&{35}&{81}&{105}&{24}&{28}
&{21}&{4}&{7}&{7}&{9}&{7}\cr
{0}&{1}&{0}&{{-1}\over{2}}&{-5}&{{-31}\over{2}}&{{-75}\over{2}}
&{-50}&{-12}&{-15}&{-12}&{{-5}\over{2}}&{-5}&{{-11}\over{2}}
&{{-15}\over{2}}&{-6}\cr
{1}&{0}&{-4}&{{-5}\over{2}}&{-15}&{{-77}\over{2}}&{{-175}\over{2}}
&{-112}&{-25}&{-28}&{-20}&{{-7}\over{2}}&{-5}&{{-7}\over{2}}
&{{-5}\over{2}}&{0}\cr\endsmallmatrix\right)$
\nobreak\newline
$\left(\smallmatrix
{-6}&{0}&{24}&{15}&{90}&{231}&{525}&{672}&{150}&{168}
&{120}&{21}&{30}&{21}&{15}&{0}\cr
{0}&{-14}&{0}&{7}&{70}&{217}&{525}&{700}&{168}&{210}
&{168}&{35}&{70}&{77}&{105}&{84}\cr
{24}&{0}&{-6}&{0}&{30}&{126}&{330}&{462}&{120}&{168}
&{150}&{36}&{90}&{126}&{210}&{210}\cr
{15}&{7}&{0}&{-1}&{0}&{14}&{45}&{70}&{21}&{35}&{36}&{10}
&{30}&{49}&{90}&{98}\cr
{90}&{70}&{30}&{0}&{-10}&{0}&{30}&{70}&{30}&{70}&{90}&{30}
&{110}&{210}&{420}&{490}\cr
{231}&{217}&{126}&{14}&{0}&{-7}&{0}&{28}&{21}&{77}&{126}&{49}
&{210}&{448}&{945}&{1148}\cr
{525}&{525}&{330}&{45}&{30}&{0}&{-15}&{0}&{15}&{105}&{210}&{90}
&{420}&{945}&{2040}&{2520}\cr
{672}&{700}&{462}&{70}&{70}&{28}&{0}&{-14}&{0}&{84}&{210}&{98}
&{490}&{1148}&{2520}&{3150}\cr
{150}&{168}&{120}&{21}&{30}&{21}&{15}&{0}&{-6}&{0}&{24}&{15}&{90}
&{231}&{525}&{672}\cr
{168}&{210}&{168}&{35}&{70}&{77}&{105}&{84}&{0}&{-14}&{0}&{7}&{70}
&{217}&{525}&{700}\cr
{120}&{168}&{150}&{36}&{90}&{126}&{210}&{210}&{24}&{0}&{-6}&{0}&{30}
&{126}&{330}&{462}\cr
{21}&{35}&{36}&{10}&{30}&{49}&{90}&{98}&{15}&{7}&{0}&{-1}&{0}&{14}
&{45}&{70}\cr
{30}&{70}&{90}&{30}&{110}&{210}&{420}&{490}&{90}&{70}&{30}&{0}
&{-10}&{0}&{30}&{70}\cr
{21}&{77}&{126}&{49}&{210}&{448}&{945}&{1148}&{231}&{217}&{126}
&{14}&{0}&{-7}&{0}&{28}\cr
{15}&{105}&{210}&{90}&{420}&{945}&{2040}&{2520}&{525}&{525}&{330}
&{45}&{30}&{0}&{-15}&{0}\cr
{0}&{84}&{210}&{98}&{490}&{1148}&{2520}&{3150}&{672}&{700}&{462}
&{70}&{70}&{28}&{0}&{-14}\cr\endsmallmatrix\right)$.}

\vskip3pt

\noindent
$N^\prime =36\ $ $d=330=2\cdot 3\cdot 5\cdot 11,\,\text{odd},\,
\eta=3,\,h=1$:
$\langle 6\rangle \oplus \langle -10\rangle \oplus \langle -22\rangle
(1/2,1/2,0)$. It is equivariantly equivalent to
$d=165,\,\eta=4,\,h=1$: $\langle 3\rangle \oplus
\langle -5\rangle\oplus \langle -11\rangle$.

\vskip3pt

\vbox{\noindent
$N^\prime=37\ $ $d=390=2\cdot 3\cdot 5\cdot 13,\,
\text{odd},\,\eta=4,\,h=1$:
$\langle 30 \rangle\oplus \langle -26 \rangle \oplus \langle -2 \rangle
(0,1/2,1/2)$
(=$\widetilde{(195,3,h=0)}$)
\nobreak\newline
$\left(\smallmatrix
{0}&{0}&{1}&{26}&{1}&{13}&{3}&{13}&{2}&{26}&{21}&{364}
&{9}&{65}&{7}&{13}\cr
{0}&{1}&{0}&{{-15}\over{2}}&{{-1}\over{2}}&{-9}&{{-5}\over{2}}
&{-12}&{-2}&{-27}&{-22}&{{-765}\over{2}}&{{-19}\over{2}}
&{-69}&{{-15}\over{2}}&{-14}\cr
{1}&{0}&{-4}&{{-195}\over{2}}&{{-7}\over{2}}&{-39}&{{-15}\over{2}}
&{-26}&{-3}&{-26}&{-18}&{{-585}\over{2}}&{{-13}\over{2}}&{-39}
&{{-5}\over{2}}&{0}\cr\endsmallmatrix\right)$
\nobreak\newline
$\left(\smallmatrix
{-2}&{0}&{8}&{195}&{7}&{78}&{15}&{52}&{6}&{52}
&{36}&{585}&{13}&{78}&{5}&{0}\cr
{0}&{-26}&{0}&{195}&{13}&{234}&{65}&{312}&{52}&{702}&{572}&{9945}
&{247}&{1794}&{195}&{364}\cr
{8}&{0}&{-2}&{0}&{2}&{78}&{30}&{182}&{36}&{572}&{486}&{8580}&{218}
&{1638}&{190}&{390}\cr
{195}&{195}&{0}&{-195}&{0}&{780}&{390}&{2730}&{585}&{9945}&{8580}
&{152295}&{3900}&{29640}&{3510}&{7410}\cr
{7}&{13}&{2}&{0}&{-1}&{0}&{5}&{52}&{13}&{247}&{218}&{3900}&{101}
&{780}&{95}&{208}\cr
{78}&{234}&{78}&{780}&{0}&{-78}&{0}&{234}&{78}&{1794}&{1638}
&{29640}&{780}&{6162}&{780}&{1794}\cr
{15}&{65}&{30}&{390}&{5}&{0}&{-5}&{0}&{5}&{195}&{190}&{3510}&{95}
&{780}&{105}&{260}\cr
{52}&{312}&{182}&{2730}&{52}&{234}&{0}&{-26}&{0}&{364}&{390}&{7410}
&{208}&{1794}&{260}&{702}\cr
{6}&{52}&{36}&{585}&{13}&{78}&{5}&{0}&{-2}&{0}&{8}&{195}&{7}&{78}
&{15}&{52}\cr
{52}&{702}&{572}&{9945}&{247}&{1794}&{195}&{364}&{0}&{-26}&{0}&{195}
&{13}&{234}&{65}&{312}\cr
{36}&{572}&{486}&{8580}&{218}&{1638}&{190}&{390}&{8}&{0}&{-2}&{0}&{2}
&{78}&{30}&{182}\cr
{585}&{9945}&{8580}&{152295}&{3900}&{29640}&{3510}&{7410}&{195}&{195}
&{0}&{-195}&{0}&{780}&{390}&{2730}\cr
{13}&{247}&{218}&{3900}&{101}&{780}&{95}&{208}&{7}&{13}&{2}&{0}&{-1}
&{0}&{5}&{52}\cr
{78}&{1794}&{1638}&{29640}&{780}&{6162}&{780}&{1794}&{78}&{234}&{78}
&{780}&{0}&{-78}&{0}&{234}\cr
{5}&{195}&{190}&{3510}&{95}&{780}&{105}&{260}&{15}&{65}&{30}&{390}
&{5}&{0}&{-5}&{0}\cr
{0}&{364}&{390}&{7410}&{208}&{1794}&{260}&{702}&{52}&{312}&{182}
&{2730}&{52}&{234}&{0}&{-26}\cr\endsmallmatrix\right)$.}

\vskip3pt

\vbox{\noindent
$N^\prime=38\ $ $d=570=2\cdot 3\cdot 5\cdot 19,\,
\text{odd},\,\eta=2,\,h=1$:
$\langle 2 \rangle\oplus \langle -190 \rangle \oplus \langle -6 \rangle
(0,1/2,1/2)$
(=$\widetilde{(285,5,h=0)}$)
\nobreak\newline
$\left(\smallmatrix
{158}&{2850}&{11}&{15}&{19}&{0}&{0}&{3}&{38}&{48}&{855}\cr
{{-27}\over{2}}&{{-489}\over{2}}&{-1}&{{-3}\over{2}}&{-2}
&{0}&{1}&{0}&{-2}&{-3}&{-56}\cr
{{-101}\over{2}}&{{-1805}\over{2}}&{-3}&{{-5}\over{2}}
&{0}&{1}&{0}&{-2}&{-19}&{-22}&{-380}\cr
\endsmallmatrix\right.$\hfill
\nobreak\newline
\line{\hfil$\left.\smallmatrix
{22}&{570}&{13}&{165}&{893}&{288}&{2280}&{285}&{874}&{600}&{8265}\cr
{{-3}\over{2}}&{{-81}\over{2}}&{-1}&{{-27}\over{2}}
&{-74}&{-24}&{-191}&{-24}&{-74}&{-51}&{-704}\cr
{{-19}\over{2}}&{{-475}\over{2}}&{-5}&{{-115}\over{2}}
&{-304}&{-97}&{-760}&{-94}&{-285}&{-194}&{-2660}\cr
\endsmallmatrix\right)$}
\nobreak\newline
$\pmatrix
A&B\\
B&A\\
\endpmatrix$\  where
\nobreak
\newline
$\pmatrix
A\\
B\\
\endpmatrix =
\left(\smallmatrix
{-1}&{0}&{2}&{135}&{874}&{303}&{2565}&{342}&{1121}&{807}&{11400}\cr
{0}&{-285}&{0}&{2280}&{15390}&{5415}&{46455}&{6270}&{20805}
&{15105}&{214320}\cr
{2}&{0}&{-2}&{0}&{38}&{18}&{190}&{30}&{114}&{90}&{1330}\cr
{135}&{2280}&{0}&{-15}&{0}&{15}&{285}&{60}&{285}&{255}&{3990}\cr
{874}&{15390}&{38}&{0}&{-38}&{0}&{380}&{114}&{684}&{684}&{11210}\cr
{303}&{5415}&{18}&{15}&{0}&{-6}&{0}&{12}&{114}&{132}&{2280}\cr
{2565}&{46455}&{190}&{285}&{380}&{0}&{-190}&{0}&{380}&{570}&{10640}\cr
{342}&{6270}&{30}&{60}&{114}&{12}&{0}&{-6}&{0}&{24}&{570}\cr
{1121}&{20805}&{114}&{285}&{684}&{114}&{380}&{0}&{-38}&{0}&{380}\cr
{807}&{15105}&{90}&{255}&{684}&{132}&{570}&{24}&{0}&{-6}&{0}\cr
{11400}&{214320}&{1330}&{3990}&{11210}&{2280}&
{10640}&{570}&{380}&{0}&{-190}\cr
{226}&{4275}&{28}&{90}&{266}&{57}&{285}&{18}&{19}&{3}&{0}\cr
{4275}&{81510}&{570}&{1995}&{6270}&{1425}&{7695}&{570}
&{855}&{285}&{2280}\cr
{28}&{570}&{6}&{30}&{114}&{30}&{190}&{18}&{38}&{18}&{190}\cr
{90}&{1995}&{30}&{240}&{1140}&{345}
&{2565}&{300}&{855}&{555}&{7410}\cr
{266}&{6270}&{114}&{1140}&{5814}&{1824}&{14060}
&{1710}&{5092}&{3420}&{46550}\cr
{57}&{1425}&{30}&{345}&{1824}&{582}&{4560}&{564}&{1710}
&{1164}&{15960}\cr
{285}&{7695}&{190}&{2565}&{14060}&{4560}&{36290}&{4560}
&{14060}&{9690}&{133760}\cr
{18}&{570}&{18}&{300}&{1710}&{564}&{4560}&{582}&{1824}
&{1272}&{17670}\cr
{19}&{855}&{38}&{855}&{5092}&{1710}&{14060}&{1824}&{5814}
&{4104}&{57380}\cr
{3}&{285}&{18}&{555}&{3420}&{1164}&{9690}&{1272}&{4104}&{2922}
&{41040}\cr
{0}&{2280}&{190}&{7410}&{46550}&{15960}&{133760}&{17670}&{57380}
&{41040}&{577790}\cr
\endsmallmatrix\right)$.}

\newpage

\centerline{\bf Table 3}
\centerline{The list of main hyperbolic lattices of the rank $3$}
\centerline{with square-free determinant $d\le 100000$ and $h\le 1$}

\vskip20pt

{\settabs 12 \columns
\+$n=1,\,d=1,\,\eta=0,\,h=0$  &&&&&&&
$U\oplus \langle-1\rangle$  &&&&&& $er$\cr
\+$n=2,\,d=2,\,\eta=0,\,h=0$  &&&&&&& $U\oplus\langle-2\rangle$
&&&&&& $er$\cr
\+$n=3,\,d=3,\,\eta=0,\,h=0$  &&&&&&& $\langle3\rangle\oplus\langle-1\rangle
                                   \oplus\langle-1\rangle$   &&&&&& $er$\cr
\+$n=4,\,d=3,\,\eta=1,\,h=0$  &&&&&&& $U\oplus\langle-3\rangle$ &&&&&& $er$\cr
\+$n=5,\,d=5,\,\eta=0,\,h=0$  &&&&&&& $U\oplus\langle-5\rangle$ &&&&&& $er$\cr
\+$n=6,\,d=5,\,\eta=1,\,h=0$  &&&&&&& $\langle1\rangle\oplus\langle-10\rangle
                                       \oplus\langle-2\rangle(0,1/2,1/2)$
                                                              &&&&&& $er$\cr
\+$n=7,\,d=6=2\cdot3,\,\eta=0,\,h=0$ &&&&&&&
$U\oplus\langle-6\rangle$ &&&&&& $er$\cr
\+$n=8,\,d=7,\,\eta=0,\,h=1$ &&&&&&&
$\langle7\rangle\oplus\langle-1\rangle\oplus\langle-1\rangle$ &&&&&& $er$\cr
\+$n=9,\,d=7,\,\eta=1,\,h=0$ &&&&&&& $U\oplus\langle-7\rangle$ &&&&&& $er$\cr
\+$n=10,\,d=10=2\cdot 5,\,\eta=1,\,h=0$ &&&&&&&
$U\oplus\langle-10\rangle$ &&&&&& $er$\cr
\+$n=11,\,d=11,\,\eta=0,\,h=0$ &&&&&&&
$\langle11\rangle\oplus\langle-1\rangle\oplus\langle-1\rangle$ &&&&&& $er$\cr
\+$n=12,\,d=11,\,\eta=1,\,h=1$ &&&&&&&
$U\oplus\langle-11\rangle$ &&&&&& $er$\cr
\+$n=13,\,d=13,\,\eta=0,\,h=1$ &&&&&&&
$U\oplus\langle-13\rangle$ &&&&&& $er$\cr
\+$n=14,\,d=13,\,\eta=1,\,h=1$ &&&&&&&
$\langle26\rangle\oplus\langle-2\rangle\oplus\langle-1\rangle
(1/2,1/2,0)$ &&&&&& $er$\cr
\+$n=15,\,d=14=2\cdot7,\,\eta=1,\,h=0$ &&&&&&&
$U\oplus\langle-14\rangle$ &&&&&&  $er$\cr
\+$n=16,\,d=15=3\cdot 5,\,\eta=0,\,h=0$ &&&&&&&
$\langle3\rangle\oplus\langle-5\rangle\oplus\langle-1\rangle$ &&&&&&  $er$\cr
\+$n=17,\,d=15=3\cdot 5,\,\eta=1,\,h=0$ &&&&&&&
$\langle5\rangle\oplus\langle-3\rangle\oplus\langle-1\rangle$ &&&&&&  $er$\cr
\+$n=18,\,d=15=3\cdot 5,\,\eta=2,\,h=0$ &&&&&&&
$U\oplus\langle-15\rangle$ &&&&&& $er$\cr
\+$n=19,\,d=15=3\cdot 5,\,\eta=3,\,h=0$ &&&&&&&
$\langle15\rangle\oplus\langle-1\rangle\oplus\langle-1\rangle$ &&&&&& $er$\cr
\+$n=20,\,d=17,\,\eta=0,\,h=0$ &&&&&&&
$U\oplus\langle-17\rangle$ &&&&&& $er$\cr
\+$n=21,\,d=19,\,\eta=0,\,h=1$ &&&&&&&
$\langle19\rangle\oplus\langle-1\rangle\oplus\langle-1\rangle$ &&&&&& $er$\cr
\+$n=22,\,d=21=3\cdot 7,\,\eta=0,\,h=0$ &&&&&&&
$\langle14\rangle\oplus\langle-6\rangle\oplus\langle-1\rangle
(1/2,1/2,0)$ &&&&&& $er$\cr
\+$n=23,\,d=21=3\cdot 7,\,\eta=1,\,h=1$ &&&&&&&
$U\oplus\langle-21\rangle$ &&&&&& $er$\cr
\+$n=24,\,d=21=3\cdot 7,\,\eta=2,\,h=0$ &&&&&&&
$\langle3\rangle\oplus\langle-7\rangle\oplus\langle-1\rangle$ &&&&&& $er$\cr
\+$n=25,\,d=21=3\cdot 7,\,\eta=3,\,h=0$ &&&&&&&
$\langle1\rangle\oplus\langle-7\rangle\oplus\langle-3\rangle$ &&&&&& $er$\cr
\+$n=26,\,d=22=2\cdot 11,\,\eta=0,\,h=1$ &&&&&&&
$U\oplus\langle-22\rangle$ &&&&&&  $er$\cr
\+$n=27,\,d=23,\,\eta=0,\,h=1$ &&&&&&&
$\langle23\rangle\oplus\langle-1\rangle\oplus\langle-1\rangle$ &&&&&& $er$\cr
\+$n=28,\,d=26=2\cdot 13,\,\eta=1,\,h=0$ &&&&&&&
$U\oplus\langle-26\rangle$ &&&&&& $er$\cr
\+$n=29,\,d=29,\eta=0,\,h=1$ &&&&&&&
$U\oplus\langle-29\rangle$ &&&&&& $nr$\cr
\+$n=30,\,d=29,\,\eta=1,\,h=1$ &&&&&&&
$\langle1\rangle\oplus\langle-58\rangle\oplus\langle-2\rangle
(0,1/2,1/2)$ &&&&&& $er$\cr
\+$n=31,\,d=30=2\cdot 3\cdot 5,\,\eta=1,\,h=1$ &&&&&&&
$U\oplus\langle-30\rangle$ &&&&&& $er$\cr
\+$n=32,\,d=30=2\cdot 3\cdot 5\cdot \eta=2,\,h=0$ &&&&&&&
$\langle10\rangle\oplus\langle-6\rangle\oplus\langle-2\rangle
(0,1/2,1/2)$ &&&&&& $er$\cr
\+$n=33,\,d=33=3\cdot 11,\,\eta=0,\,h=1$ &&&&&&&
$\langle3\rangle\oplus\langle-22\rangle\oplus\langle-2\rangle
(0,1/2,1/2)$ &&&&&& $er$\cr
\+$n=34,\,d=33=3\cdot 11,\,\eta=1,\,h=0$ &&&&&&&
$\langle11\rangle\oplus\langle-3\rangle\oplus\langle-1\rangle$ &&&&&& $er$\cr
\+$n=35,\,d=33=3\cdot 11,\,\eta=2,\,h=1$ &&&&&&&
$U\oplus\langle-33\rangle$ &&&&&& $er$\cr
\+$n=36,\,d=33=3\cdot 11,\,\eta=3,\,h=1$ &&&&&&&
$\langle1\rangle\oplus\langle-11\rangle\oplus\langle-3\rangle$ &&&&&& $er$\cr
\+$n=37,\,d=34=2\cdot17,\,\eta=0,\,h=1$ &&&&&&&
$U\oplus\langle-34\rangle$ &&&&&& $er$\cr
\+$n=38,\,d=35=5\cdot 7,\,\eta=0,\,h=0$ &&&&&&&
$\langle7\rangle\oplus\langle-5\rangle\oplus\langle-1\rangle$ &&&&&& $er$\cr
\+$n=39,\,d=35=5\cdot 7,\,\eta=2,\,h=1$ &&&&&&&
$\langle1\rangle\oplus\langle-7\rangle\oplus\langle-5\rangle$ &&&&&&  $er$\cr
\+$n=40,\,d=35=5\cdot 7,\,\eta=3,\,h=0$ &&&&&&&
$\langle35\rangle\oplus\langle-1\rangle\oplus\langle-1\rangle$ &&&&&& $er$\cr
\+$n=41,\,d=38=2\cdot 19,\,\eta=0,\,h=1$ &&&&&&&
$U\oplus\langle-38\rangle$ &&&&&& $er$\cr
\+$n=42,\,d=39=3\cdot 13,\,\eta=0,\,h=0$ &&&&&&&
$\langle3\rangle\oplus\langle-13\rangle\oplus\langle-1\rangle$ &&&&&& $er$\cr
\+$n=43,\,d=39=3\cdot 13,\,\eta=2,\,h=1$ &&&&&&&
$\langle3\rangle\oplus\langle-26\rangle\oplus\langle-2\rangle
(0,1/2,1/2)$ &&&&&& $er$\cr
\+$n=44,\,d=39=3\cdot 13,\,\eta=3,\,h=1$ &&&&&&&
$\langle2\rangle\oplus\langle-26\rangle\oplus\langle-3\rangle
(1/2,1/2,0)$ &&&&&& $er$\cr
\+$n=45,\,d=41,\,\eta=0,\,h=1$ &&&&&&&
$U\oplus\langle-41\rangle$ &&&&&&  $nr$\cr
\+$n=46,\,d=42=2\cdot 3\cdot 7,\,\eta=0,\,h=0$ &&&&&&&
$U\oplus\langle-42\rangle$ &&&&&& $er$\cr
\+$n=47,\,d=42=2\cdot 3\cdot 7,\,\eta=3,\,h=1$ &&&&&&&
$\langle6\rangle\oplus\langle-14\rangle\oplus\langle-2\rangle
(1/2,1/2,0)$ &&&&&& $er$\cr
\+$n=48,\,d=51=3\cdot 17,\,\eta=0,\,h=0$ &&&&&&&
$\langle3\rangle\oplus\langle-17\rangle\oplus\langle-1\rangle$ &&&&&& $er$\cr
\+$n=49,\,d=51=3\cdot 17,\,\eta=1,\,h=1$ &&&&&&&
$\langle17\rangle\oplus\langle-3\rangle\oplus\langle-1\rangle$ &&&&&& $er$\cr
\+$n=50,\,d=51=3\cdot 17,\,\eta=3,\,h=1$ &&&&&&&
$\langle51\rangle\oplus\langle-1\rangle\oplus\langle-1\rangle$ &&&&&& $er$\cr
\+$n=51,\,d=55=5\cdot 11,\,\eta=0,\,h=1$ &&&&&&&
$\langle11\rangle\oplus\langle-5\rangle\oplus\langle-1\rangle$ &&&&&& $er$\cr
\+$n=52,\,d=55=5\cdot 11,\,\eta=3,\,h=1$ &&&&&&&
$\langle2\rangle\oplus\langle-11\rangle\oplus\langle-10\rangle
(1/2,0,1/2)$ &&&&&& $er$\cr
\+$n=53,\,d=57=3\cdot 19,\,\eta=0,\,h=1$ &&&&&&&
$\langle1\rangle\oplus\langle-38\rangle\oplus\langle-6\rangle
(0,1/2,1/2)$ &&&&&& $er$\cr
\+$n=54,\,d=57=3\cdot 19,\,\eta=1,\,h=1$ &&&&&&&
$U\oplus\langle-57\rangle$ &&&&&& $nr$\cr
\+$n=55,\,d=57=3\cdot 19,\,\eta=2,\,h=1$ &&&&&&&
$\langle3\rangle\oplus\langle-19\rangle\oplus\langle-1\rangle$ &&&&&& $er$\cr
\+$n=56,\,d=65=5\cdot 13,\,\eta=0,\,h=0$ &&&&&&&
$\langle1\rangle\oplus\langle-13\rangle\oplus\langle-5\rangle$ &&&&&& $hr$\cr
\+$n=57,\,d=65=5\cdot 13,\,\eta=2,\,h=1$ &&&&&&&
$\langle5\rangle\oplus\langle-26\rangle\oplus\langle-2\rangle
(0,1/2,1/2)$ &&&&&& $er$\cr
\+$n=58,\,d=65=5\cdot 13,\,\eta=3,\,h=1$ &&&&&&&
$U\oplus\langle-65\rangle$ &&&&&& $nr$\cr
\+$n=59,\,d=66=2\cdot 3\cdot 11,\,\eta=1,\,h=1$ &&&&&&&
$U\oplus\langle-66\rangle$ &&&&&& $er$\cr
\+$n=60,\,d=66=2\cdot 3\cdot 11,\,\eta=2,\,h=1$ &&&&&&&
$\langle22\rangle\oplus\langle-6\rangle\oplus\langle-2\rangle
(0,1/2,1/2)$ &&&&&& $er$\cr
\+$n=61,\,d=69=3\cdot 23,\,\eta=3,\,h=0$ &&&&&&&
$\langle1\rangle\oplus\langle-23\rangle\oplus\langle-3\rangle$ &&&&&& $er$\cr
\+$n=62,\,d=70=2\cdot 5\cdot 7,\,\eta=3,\,h=1$ &&&&&&&
$\langle2\rangle\oplus\langle-14\rangle\oplus\langle-10\rangle
(0,1/2,1/2)$ &&&&&& $er$\cr
\+$n=63,\,d=71,\,\eta=0,\,h=1$ &&&&&&&
$\langle71\rangle\oplus\langle-1\rangle\oplus\langle-1\rangle$ &&&&&& $nr$\cr
\+$n=64,\,d=74=2\cdot 37,\,\eta=1,\,h=0$ &&&&&&&
$U\oplus\langle-74\rangle$ &&&&&& $hr$\cr
\+$n=65,\,d=77=7\cdot 11,\,\eta=1,\,h=0$ &&&&&&&
$U\oplus\langle-77\rangle$ &&&&&& $nr$\cr
\+$n=66,\,d=77=7\cdot 11,\,\eta=3,\,h=0$ &&&&&&&
$\langle1\rangle\oplus\langle-11\rangle\oplus\langle-7\rangle$ &&&&&& $er$\cr
\+$n=67,\,d=78=2\cdot 3\cdot 13,\,\eta=2,\,h=0$ &&&&&&&
$U\oplus\langle-78\rangle$ &&&&&& $er$\cr
\+$n=68,\,d=85=5\cdot 17,\,\eta=0,\,h=1$ &&&&&&&
$\langle1\rangle\oplus\langle-17\rangle\oplus\langle-5\rangle$ &&&&&& $nr$\cr
\+$n=69,\,d=85=5\cdot 17,\,\eta=1,\,h=1$ &&&&&&&
$\langle1\rangle\oplus\langle-34\rangle\oplus\langle-10\rangle
(0,1/2,1/2)$ &&&&&& $er$\cr
\+$n=70,\,d=86=2\cdot 43,\,\eta=0,\,h=1$ &&&&&&&
$U\oplus\langle-86\rangle$ &&&&&& $nr$\cr
\+$n=71,\,d=87=3\cdot 29,\,\eta=3,\,h=1$ &&&&&&&
$\langle2\rangle\oplus\langle-58\rangle\oplus\langle-3\rangle
(1/2,1/2,0)$ &&&&&& $er$\cr
\+$n=72,\,d=91=7\cdot 13,\,\eta=3,\,h=1$ &&&&&&&
$\langle1\rangle\oplus\langle-26\rangle\oplus\langle-14\rangle
(0,1/2,1/2)$ &&&&&& $er$\cr
\+$n=73,\,d=93=3\cdot 31,\,\eta=0,\,h=1$ &&&&&&&
$\langle31\rangle\oplus\langle-6\rangle\oplus\langle-2\rangle
(0,1/2,1/2)$ &&&&&& $er$\cr
\+$n=74,\,d=95=5\cdot 19,\,\eta=0,\,h=0$ &&&&&&&
$\langle19\rangle\oplus\langle-5\rangle\oplus\langle-1\rangle$ &&&&&& $hr$\cr
\+$n=75,\,d=95=5\cdot 19,\,\eta=3,\,h=1$ &&&&&&&
$\langle2\rangle\oplus\langle-19\rangle\oplus\langle-10\rangle
(1/2,0,1/2)$ &&&&&& $er$\cr
\+$n=76,\,d=102=2\cdot 3\cdot 17,\,\eta=0,\,h=1$ &&&&&&&
$\langle2\rangle\oplus\langle-34\rangle\oplus\langle-6\rangle
(0,1/2,1/2)$ &&&&&& $er$\cr
\+$n=77,\,d=105=3\cdot 5\cdot 7,\,\eta=0,\,h=1$ &&&&&&&
$\langle2\rangle\oplus\langle-42\rangle\oplus\langle-5\rangle
(1/2,1/2,0)$ &&&&&& $er$\cr
\+$n=78,\,d=105=3\cdot 5\cdot 7,\,\eta=1,\,h=0$ &&&&&&&
$\langle1\rangle\oplus\langle-21\rangle\oplus\langle-5\rangle$ &&&&&& $er$\cr
\+$n=79,\,d=105=3\cdot 5\cdot 7,\,\eta=2,\,h=1$ &&&&&&&
$\langle7\rangle\oplus\langle-15\rangle\oplus\langle-1\rangle$ &&&&&& $er$\cr
\+$n=80,\,d=105=3\cdot 5\cdot 7,\,\eta=3,\,h=0$ &&&&&&&
$\langle14\rangle\oplus\langle-10\rangle\oplus\langle-3\rangle
(1/2,1/2,0)$ &&&&&& $er$\cr
\+$n=81,\,d=105=3\cdot 5\cdot 7,\,\eta=4,\,h=1$ &&&&&&&
$U\oplus\langle-105\rangle$ &&&&&& $nr$\cr
\+$n=82,\,d=105=3\cdot 5\cdot 7,\,\eta=5,\,h=0$ &&&&&&&
$\langle5\rangle\oplus\langle-7\rangle\oplus\langle-3\rangle$ &&&&&& $er$\cr
\+$n=83,\,d=105=3\cdot 5\cdot 7,\,\eta=6,\,h=0$ &&&&&&&
$\langle1\rangle\oplus\langle-15\rangle\oplus\langle-7\rangle$ &&&&&& $er$\cr
\+$n=84,\,d=105=3\cdot 5\cdot 7,\,\eta=7,\,h=0$ &&&&&&&
$\langle10\rangle\oplus\langle-14\rangle\oplus\langle-3\rangle
(1/2,1/2,0)$ &&&&&& $er$\cr
\+$n=85,\,d=110=2\cdot 5\cdot 11,\,\eta=1,\,h=1$ &&&&&&&
$U\oplus\langle-110\rangle$ &&&&&& $er$\cr
\+$n=86,\,d=111=3\cdot 37,\,\eta=0,\,h=1$ &&&&&&&
$\langle3\rangle\oplus\langle-37\rangle\oplus\langle-1\rangle$ &&&&&& $er$\cr
\+$n=87,\,d=114=2\cdot 3\cdot 19,\,\eta=2,\,h=0$ &&&&&&&
$U\oplus\langle-114\rangle$ &&&&&& $hr$\cr
\+$n=88,\,d=119=7\cdot 17,\,\eta=3,\,h=1$ &&&&&&&
$\langle119\rangle\oplus\langle-1\rangle\oplus\langle-1\rangle$ &&&&&& $nr$\cr
\+$n=89,\,d=130=2\cdot 5\cdot 13,\,\eta=0,\,h=1$ &&&&&&&
$U\oplus\langle-130\rangle$ &&&&&& $nr$\cr
\+$n=90,\,d=130=2\cdot 5\cdot 13,\,\eta=3,\,h=1$ &&&&&&&
$\langle10\rangle\oplus\langle-26\rangle\oplus\langle-2\rangle
(1/2,1/2,0)$ &&&&&& $er$\cr
\+$n=91,\,d=134=2\cdot 67,\,\eta=0,\,h=1$ &&&&&&&
$U\oplus\langle-134\rangle$ &&&&&& $nr$\cr
\+$n=92,\,d=141=3\cdot 47,\,\eta=3,\,h=1$ &&&&&&&
$\langle1\rangle\oplus\langle-47\rangle\oplus\langle-3\rangle$ &&&&&& $er$\cr
\+$n=93,\,d=143=11\cdot 13,\,\eta=0,\,h=1$ &&&&&&&
$\langle11\rangle\oplus\langle-13\rangle\oplus\langle-1\rangle$ &&&&&& $nr$\cr
\+$n=94,\,d=146=2\cdot 73,\,\eta=0,\,h=1$ &&&&&&&
$U\oplus\langle-146\rangle$ &&&&&& $nr$\cr
\+$n=95,\,d=154=2\cdot 7\cdot 11,\,\eta=3,\,h=1$ &&&&&&&
$U\oplus\langle-154\rangle$ &&&&&& $nr$\cr
\+$n=96,\,d=155=5\cdot 31,\,\eta=3,\,h=1$ &&&&&&&
$\langle2\rangle\oplus\langle-31\rangle\oplus\langle-10\rangle
(1/2,0,1/2)$ &&&&&& $er$\cr
\+$n=97,\,d=161=7\cdot 23,\,\eta=1,\,h=1$ &&&&&&&
$U\oplus\langle-161\rangle$ &&&&&& $nr$\cr
\+$n=98,\,d=161=7\cdot 23,\,\eta=3,\,h=1$ &&&&&&&
$\langle1\rangle\oplus\langle-23\rangle\oplus\langle-7\rangle$ &&&&&& $nr$\cr
\+$n=99,\,d=165=3\cdot 5\cdot 11,\,\eta=0,\,h=0$ &&&&&&&
$\langle5\rangle\oplus\langle-6\rangle\oplus\langle-22\rangle
(0,1/2,1/2)$ &&&&&& $er$\cr
\+$n=100,\,d=165=3\cdot 5\cdot 11,\,\eta=1,\,h=1$ &&&&&&&
$\langle11\rangle\oplus\langle-5\rangle\oplus\langle-3\rangle$ &&&&&& $nr$\cr
\+$n=101,\,d=165=3\cdot 5\cdot 11,\,\eta=3,\,h=1$ &&&&&&&
$\langle15\rangle\oplus\langle-22\rangle\oplus\langle-2\rangle
(0,1/2,1/2)$ &&&&&& $er$\cr
\+$n=102,\,d=165=3\cdot 5\cdot11,\,\eta=4,\,h=1$ &&&&&&&
$\langle3\rangle\oplus\langle-5\rangle\oplus\langle-11\rangle$ &&&&&& $er$\cr
\+$n=103,\,d=165=3\cdot 5\cdot 11,\,\eta=5,\,h=1$ &&&&&&&
$\langle1\rangle\oplus\langle-55\rangle\oplus\langle-3\rangle$ &&&&&& $er$\cr
\+$n=104,\,d=165=3\cdot 5\cdot 11,\,\eta=6,\,h=0$ &&&&&&&
$\langle1\rangle\oplus\langle-15\rangle\oplus\langle-11\rangle$ &&&&&& $er$\cr
\+$n=105,\,d=170=2\cdot 5\cdot 17,\,\eta=1,\,h=1$ &&&&&&&
$\langle2\rangle\oplus\langle-10\rangle\oplus\langle-34\rangle
(1/2,1/2,0)$ &&&&&& $er$\cr
\+$n=106,\,d=170=2\cdot 5\cdot 17,\,\eta=2,\,h=0$ &&&&&&&
$U\oplus\langle-170\rangle$ &&&&&& $nr$\cr
\+$n=107,\,d=182=2\cdot 7\cdot 13,\,\eta=3,\,h=1$ &&&&&&&
$\langle2\rangle\oplus\langle-26\rangle\oplus\langle-14\rangle
(0,1/2,1/2)$ &&&&&& $nr$\cr
\+$n=108,\,d=186=2\cdot 3\cdot 31,\,\eta=0,\,h=0$ &&&&&&&
$U\oplus\langle-186\rangle$ &&&&&& $nr$\cr
\+$n=109,\,d=194=2\cdot 97,\,\eta=0,\,h=1$  &&&&&&&
$U\oplus\langle-194\rangle$ &&&&&& $nr$\cr
\+$n=110,\,d=195=3\cdot 5\cdot 13,\,\eta=3,\,h=0$ &&&&&&&
$\langle15\rangle\oplus\langle-13\rangle\oplus\langle-1\rangle$ &&&&&& $er$\cr
\+$n=111,\,d=195=3\cdot 5\cdot 13,\,\eta=5,\,h=1$ &&&&&&&
$\langle6\rangle\oplus\langle-26\rangle\oplus\langle-5\rangle
(1/2,1/2,0)$ &&&&&& $er$\cr
\+$n=112,\,d=195=3\cdot 5\cdot 13,\,\eta=6,\,h=0$ &&&&&&&
$\langle3\rangle\oplus\langle-65\rangle\oplus\langle-1\rangle$ &&&&&& $er$\cr
\+$n=113,\,d=203=7\cdot 29,\,\eta=3,\,h=1$  &&&&&&&
$\langle2\rangle\oplus\langle-58\rangle\oplus\langle-7\rangle
(1/2,1/2,0)$ &&&&&& $nr$\cr
\+$n=114,\,d=205=5\cdot 41,\,\eta=1,\,h=1$  &&&&&&&
$\langle2\rangle\oplus\langle-41\rangle\oplus\langle-10\rangle
(1/2,0,1/2)$ &&&&&& $er$\cr
\+$n=115,\,d=210=2\cdot 3\cdot 5\cdot 7,\,\eta=2,\,h=0$ &&&&&&&
$\langle2\rangle\oplus\langle-10\rangle\oplus\langle-42\rangle
(1/2,1/2,0)$ &&&&&& $er$\cr
\+$n=116,\,d=210=2\cdot 3\cdot 5\cdot 7,\,\eta=4,\,h=1$ &&&&&&&
$\langle30\rangle\oplus\langle-14\rangle\oplus\langle-2\rangle
(0,1/2,1/2)$ &&&&&& $er$\cr
\+$n=117,\,d=215=5\cdot 43,\,\eta=0,\,h=1$ &&&&&&&
$\langle43\rangle\oplus\langle-5\rangle\oplus\langle-1\rangle$ &&&&&& $nr$\cr
\+$n=118,\,d=219=3\cdot 73,\,\eta=0,\,h=0$ &&&&&&&
$\langle3\rangle\oplus\langle-73\rangle\oplus\langle-1\rangle$ &&&&&& $er$\cr
\+$n=119,\,d=231=3\cdot 7\cdot 11,\,\eta=1,\,h=1$  &&&&&&&
$\langle14\rangle\oplus\langle-22\rangle\oplus\langle-3\rangle
(1/2,1/2,0)$ &&&&&& $er$\cr
\+$n=120,\,d=231=3\cdot 7\cdot 11,\,\eta=2,\,h=0$  &&&&&&&
$\langle3\rangle\oplus\langle-22\rangle\oplus\langle-14\rangle
(0,1/2,1/2)$ &&&&&& $hr$\cr
\+$n=121,\,d=231=3\cdot 7\cdot 11,\,\eta=7,\,h=0$ &&&&&&&
$\langle6\rangle\oplus\langle-14\rangle\oplus\langle-11\rangle
(1/2,1/2,0)$ &&&&&& $er$\cr
\+$n=122,\,d=246=2\cdot 3\cdot 41,\,\eta=0,\,h=1$ &&&&&&&
$\langle2\rangle\oplus\langle-82\rangle\oplus\langle-6\rangle
(0,1/2,1/2)$ &&&&&& $nr$\cr
\+$n=123,\,d=255=3\cdot 5\cdot 17,\,\eta=0,\,h=1$ &&&&&&&
$\langle3\rangle\oplus\langle-17\rangle\oplus\langle-5\rangle$ &&&&&& $er$\cr
\+$n=124,\,d=255=3\cdot 5\cdot 17,\,\eta=5,\,h=1$ &&&&&&&
$\langle51\rangle\oplus\langle-5\rangle\oplus\langle-1\rangle$ &&&&&& $nr$\cr
\+$n=125,\,d=255=3\cdot 5\cdot 17,\,\eta=6,\,h=0$ &&&&&&&
$\langle3\rangle\oplus\langle-85\rangle\oplus\langle-1\rangle$ &&&&&& $er$\cr
\+$n=126,\,d=259=7\cdot 37,\,\eta=3,\,h=1$ &&&&&&&
$\langle2\rangle\oplus\langle-74\rangle\oplus\langle-7\rangle
(1/2,1/2,0)$ &&&&&& $nr$\cr
\+$n=127,\,d=266=2\cdot 7\cdot 19,\,\eta=0,\,h=1$ &&&&&&&
$U\oplus\langle-266\rangle$ &&&&&& $nr$\cr
\+$n=128,\,d=266=2\cdot 7\cdot 19,\,\eta=3,\,h=1$ &&&&&&&
$\langle38\rangle\oplus\langle-14\rangle\oplus\langle-2\rangle
(0,1/2,1/2)$ &&&&&& $nr$\cr
\+$n=129,\,d=273=3\cdot 7\cdot 13,\,\eta=0,\,h=1$  &&&&&&&
$\langle14\rangle\oplus\langle-13\rangle\oplus\langle-6\rangle
(1/2,0,1/2)$ &&&&&& $er$\cr
\+$n=130,\,d=273=3\cdot 7\cdot 13,\,\eta=6,\,h=1$ &&&&&&&
$\langle26\rangle\oplus\langle-7\rangle\oplus\langle-6\rangle
(1/2,0,1/2)$ &&&&&& $er$\cr
\+$n=131,\,d=285=3\cdot 5\cdot 19,\,\eta=1,\,h=1$ &&&&&&&
$\langle19\rangle\oplus\langle-5\rangle\oplus\langle-3\rangle$ &&&&&& $nr$\cr
\+$n=132,\,d=285=3\cdot 5\cdot 19,\,\eta=3,\,h=1$ &&&&&&&
$\langle10\rangle\oplus\langle-38\rangle\oplus\langle-3\rangle
(1/2,1/2,0)$ &&&&&& $er$\cr
\+$n=133,\,d=285=3\cdot 5\cdot 19,\,\eta=5,\,h=0$ &&&&&&&
$\langle1\rangle\oplus\langle-95\rangle\oplus\langle-3\rangle$ &&&&&& $er$\cr
\+$n=134,\,d=285=3\cdot 5\cdot 19,\,\eta=6,\,h=1$ &&&&&&&
$\langle57\rangle\oplus\langle-10\rangle\oplus\langle-2\rangle
(0,1/2,1/2)$ &&&&&& $nr$\cr
\+$n=135,\,d=291=3\cdot 97,\,\eta=0,\,h=1$ &&&&&&&
$\langle3\rangle\oplus\langle-97\rangle\oplus\langle-1\rangle$ &&&&&& $er$\cr
\+$n=136,\,d=299=13\cdot 23,\,\eta=3,\,h=1$ &&&&&&&
$\langle2\rangle\oplus\langle-26\rangle\oplus\langle-23\rangle
(1/2,1/2,0)$ &&&&&& $nr$\cr
\+$n=137,\,d=326=2\cdot 163,\,\eta=0,\,h=1$ &&&&&&&
$U\oplus\langle-326\rangle$ &&&&&& $nr$\cr
\+$n=138,\,d=330=2\cdot 3\cdot 5\cdot 11,\,\eta=0,\,h=1$ &&&&&&&
$U\oplus\langle-330\rangle$ &&&&&& $nr$\cr
\+$n=139,\,d=330=2\cdot 3\cdot 5\cdot 11,\,\eta=3,\,h=0$ &&&&&&&
$\langle6\rangle\oplus\langle-110\rangle\oplus\langle-2\rangle
(1/2,1/2,0)$ &&&&&& $er$\cr
\+$n=140,\,d=335=5\cdot 67,\,\eta=0,\,h=1$ &&&&&&&
$\langle67\rangle\oplus\langle-5\rangle\oplus\langle-1\rangle$ &&&&&& $nr$\cr
\+$n=141,\,d=345=3\cdot 5\cdot 23,\,\eta=6,\,h=0$ &&&&&&&
$\langle1\rangle\oplus\langle-23\rangle\oplus\langle-15\rangle$ &&&&&& $er$\cr
\+$n=142,\,d=354=2\cdot 3\cdot 59,\,\eta=2,\,h=1$ &&&&&&&
$\langle118\rangle\oplus\langle-6\rangle\oplus\langle-2\rangle
(0,1/2,1/2)$ &&&&&& $nr$\cr
\+$n=143,\,d=357=3\cdot 7\cdot 17,\,\eta=0,\,h=1$ &&&&&&&
$\langle14\rangle\oplus\langle-17\rangle\oplus\langle-6\rangle
(1/2,0,1/2)$ &&&&&& $nr$\cr
\+$n=144,\,d=357=3\cdot 7\cdot 17,\,\eta=3,\,h=1$ &&&&&&&
$\langle17\rangle\oplus\langle-7\rangle\oplus\langle-3\rangle$ &&&&&& $er$\cr
\+$n=145,\,d=357=3\cdot 7\cdot 17,\,\eta=5,\,h=1$ &&&&&&&
$\langle1\rangle\oplus\langle-119\rangle\oplus\langle-3\rangle$ &&&&&& $er$\cr
\+$n=146,\,d=371=7\cdot 53,\,\eta=3,\,h=1$ &&&&&&&
$\langle2\rangle\oplus\langle-106\rangle\oplus\langle-7\rangle
(1/2,1/2,0)$ &&&&&& $nr$\cr
\+$n=147,\,d=385=5\cdot 7\cdot 11,\,\eta=3,\,h=1$ &&&&&&&
$\langle11\rangle\oplus\langle-14\rangle\oplus\langle-10\rangle
(0,1/2,1/2)$ &&&&&& $nr$\cr
\+$n=148,\,d=385=5\cdot 7\cdot 11,\,\eta=6,\,h=1$ &&&&&&&
$\langle5\rangle\oplus\langle-11\rangle\oplus\langle-7\rangle$ &&&&&& $er$\cr
\+$n=149,\,d=386=2\cdot 193,\,\eta=0,\,h=1$ &&&&&&&
$U\oplus\langle-386\rangle$ &&&&&& $nr$\cr
\+$n=150,\,d=390=2\cdot 3\cdot 5\cdot 13,\,\eta=6,\,h=1$ &&&&&&&
$\langle10\rangle\oplus\langle-26\rangle\oplus\langle-6\rangle
(1/2,1/2,0)$ &&&&&& $er$\cr
\+$n=151,\,d=399=3\cdot 7\cdot 19,\,\eta=2,\,h=1$ &&&&&&&
$\langle3\rangle\oplus\langle-38\rangle\oplus\langle-14\rangle
(0,1/2,1/2)$ &&&&&& $nr$\cr
\+$n=152,\,d=399=3\cdot 7\cdot 19,\,\eta=4,\,h=0$ &&&&&&&
$\langle14\rangle\oplus\langle-19\rangle\oplus\langle-6\rangle
(1/2,0,1/2)$ &&&&&& $er$\cr
\+$n=153,\,d=407=11\cdot 37,\,\eta=0,\,h=1$ &&&&&&&
$\langle11\rangle\oplus\langle-37\rangle\oplus\langle-1\rangle$ &&&&&& $nr$\cr
\+$n=154,\,d=429=3\cdot 11\cdot 13,\,\eta=3,\,h=1$ &&&&&&&
$\langle1\rangle\oplus\langle-39\rangle\oplus\langle-11\rangle$ &&&&&& $er$\cr
\+$n=155,\,d=435=3\cdot 5\cdot 29,\,\eta=0,\,h=0$ &&&&&&&
$\langle3\rangle\oplus\langle-29\rangle\oplus\langle-5\rangle$ &&&&&&  $er$\cr
\+$n=156,\,d=435=3\cdot 5\cdot 29,\,\eta=6,\,h=1$ &&&&&&&
$\langle435\rangle\oplus\langle-1\rangle\oplus\langle-1\rangle$
&&&&&&  $er$\cr
\+$n=157,\,d=455=5\cdot 7\cdot 13,\,\eta=0,\,h=1$ &&&&&&&
$\langle7\rangle\oplus\langle-13\rangle\oplus\langle-5\rangle$ &&&&&& $nr$\cr
\+$n=158,\,d=455=5\cdot 7\cdot 13,\,\eta=5,\,h=1$ &&&&&&&
$\langle7\rangle\oplus\langle-26\rangle\oplus\langle-10\rangle
(0,1/2,1/2)$ &&&&&& $er$\cr
\+$n=159,\,d=465=3\cdot 5\cdot 31,\,\eta=5,\,h=1$ &&&&&&&
$\langle5\rangle\oplus\langle-31\rangle\oplus\langle-3\rangle$ &&&&&& $er$\cr
\+$n=160,\,d=483=3\cdot 7\cdot 23,\,\eta=7,\,h=1$ &&&&&&&
$\langle6\rangle\oplus\langle-46\rangle\oplus\langle-7\rangle
(1/2,1/2,0)$ &&&&&& $er$\cr
\+$n=161,\,d=506=2\cdot 11\cdot 23,\,\eta=0,\,h=1$ &&&&&&&
$U\oplus\langle-506\rangle$ &&&&&& $nr$\cr
\+$n=162,\,d=530=2\cdot 5\cdot 53,\,\eta=0,\,h=1$ &&&&&&&
$U\oplus\langle-530\rangle$ &&&&&& $nr$\cr
\+$n=163,\,d=534=2\cdot 3\cdot 89,\,\eta=0,\,h=1$ &&&&&&&
$\langle2\rangle\oplus\langle-178\rangle\oplus\langle-6\rangle
(0,1/2,1/2)$ &&&&&& $nr$\cr
\+$n=164,\,d=546=2\cdot 3\cdot 7\cdot 13,\,\eta=2,\,h=0$ &&&&&&&
$U\oplus\langle-546\rangle$ &&&&&& $nr$\cr
\+$n=165,\,d=546=2\cdot 3\cdot 7\cdot 13,\,\eta=4,\,h=0$ &&&&&&&
$\langle14\rangle\oplus\langle-26\rangle\oplus\langle-6\rangle
(1/2,0,1/2)$ &&&&&& $nr$\cr
\+$n=166,\,d=561=3\cdot 11\cdot 17,\,\eta=2,\,h=1$ &&&&&&&
$\langle3\rangle\oplus\langle-11\rangle\oplus\langle-17\rangle$ &&&&&& $nr$\cr
\+$n=167,\,d=570=2\cdot 3\cdot 5\cdot 19,\,\eta=6,\,h=1$ &&&&&&&
$\langle38\rangle\oplus\langle-10\rangle\oplus\langle-6\rangle
(0,1/2,1/2)$ &&&&&& $er$\cr
\+$n=168,\,d=602=2\cdot 7\cdot 43,\,\eta=3,\,h=1$ &&&&&&&
$U\oplus\langle-602\rangle$ &&&&&& $nr$\cr
\+$n=169,\,d=615=3\cdot 5\cdot 41,\,\eta=0,\,h=1$ &&&&&&&
$\langle3\rangle\oplus\langle-41\rangle\oplus\langle-5\rangle$ &&&&&&  $er$\cr
\+$n=170,\,d=645=3\cdot 5\cdot 43,\,\eta=3,\,h=1$ &&&&&&&
$\langle10\rangle\oplus\langle-86\rangle\oplus\langle-3\rangle
(1/2,1/2,0)$ &&&&&& $er$\cr
\+$n=171,\,d=645=3\cdot 5\cdot 43,\,\eta=6,\,h=1$ &&&&&&&
$\langle10\rangle\oplus\langle-6\rangle\oplus\langle-43\rangle
(1/2,1/2,0)$ &&&&&& $nr$\cr
\+$n=172,\,d=651=3\cdot 7\cdot 31,\,\eta=2,\,h=1$ &&&&&&&
$\langle62\rangle\oplus\langle-7\rangle\oplus\langle-6\rangle
(1/2,0,1/2)$ &&&&&& $nr$\cr
\+$n=173,\,d=651=3\cdot 7\cdot 31,\,\eta=4,\,h=1$ &&&&&&&
$\langle14\rangle\oplus\langle-31\rangle\oplus\langle-6\rangle
(1/2,0,1/2)$ &&&&&& $er$\cr
\+$n=174,\,d=714=2\cdot 3\cdot 7\cdot 17,\,\eta=0,\,h=1$ &&&&&&&
$\langle14\rangle\oplus\langle-34\rangle\oplus\langle-6\rangle
(1/2,0,1/2)$ &&&&&& $nr$\cr
\+$n=175,\,d=714=2\cdot 3\cdot 7\cdot 17,\,\eta=6,\,h=1$ &&&&&&&
$\langle102\rangle\oplus\langle-14\rangle\oplus\langle-2\rangle
(0,1/2,1/2)$ &&&&&& $nr$\cr
\+$n=176,\,d=759=3\cdot 11\cdot 23,\,\eta=4,\,h=1$ &&&&&&&
$\langle11\rangle\oplus\langle-69\rangle\oplus\langle-1\rangle$ &&&&&& $nr$\cr
\+$n=177,\,d=777=3\cdot 7\cdot 37,\,\eta=6,\,h=1$ &&&&&&&
$\langle74\rangle\oplus\langle-6\rangle\oplus\langle-7\rangle
(1/2,1/2,0)$ &&&&&& $nr$\cr
\+$n=178,\,d=795=3\cdot 5\cdot 53,\,\eta=6,\,h=1$ &&&&&&&
$\langle3\rangle\oplus\langle-106\rangle\oplus\langle-10\rangle
(0,1/2,1/2)$ &&&&&&  $er$\cr
\+$n=179,\,d=805=5\cdot 7\cdot 23,\,\eta=3,\,h=1$ &&&&&&&
$\langle1\rangle\oplus\langle-115\rangle\oplus\langle-7\rangle$ &&&&&& $nr$\cr
\+$n=180,\,d=854=2\cdot 7\cdot 61,\,\eta=3,\,h=1$ &&&&&&&
$\langle2\rangle\oplus\langle-122\rangle\oplus\langle-14\rangle
(0,1/2,1/2)$ &&&&&&  $nr$\cr
\+$n=181,\,d=897=3\cdot 13\cdot 23,\,\eta=5,\,h=1$ &&&&&&&
$\langle13\rangle\oplus\langle-23\rangle\oplus\langle-3\rangle$ &&&&&& $nr$\cr
\+$n=182,\,d=897=3\cdot 13\cdot 23,\,\eta=7,\,h=1$ &&&&&&&
$\langle6\rangle\oplus\langle-26\rangle\oplus\langle-23\rangle
(1/2,1/2,0)$ &&&&&& $nr$\cr
\+$n=183,\,d=930=2\cdot 3\cdot 5\cdot 31,\,\eta=2,\,h=1$ &&&&&&&
$\langle2\rangle\oplus\langle-186\rangle\oplus\langle-10\rangle
(1/2,0,1/2)$ &&&&&& $nr$\cr
\+$n=184,\,d=935=5\cdot 11\cdot 17,\,\eta=6,\,h=1$ &&&&&&&
$\langle187\rangle\oplus\langle-5\rangle\oplus\langle-1\rangle$ &&&&&& $nr$\cr
\+$n=185,\,d=966=2\cdot 3\cdot 7\cdot 23,\,\eta=2,\,h=1$ &&&&&&&
$\langle46\rangle\oplus\langle-14\rangle\oplus\langle-6\rangle
(1/2,1/2,0)$ &&&&&& $nr$\cr
\+$n=186,\,d=1106=2\cdot 7\cdot 79,\,\eta=1,\,h=1$ &&&&&&&
$U\oplus\langle-1106\rangle$ &&&&&& $nr$\cr
\+$n=187,\,d=1155=3\cdot 5\cdot 7\cdot 11,\,\eta=2,\,h=0$ &&&&&&&
$\langle14\rangle\oplus\langle-22\rangle\oplus\langle-15\rangle
(1/2,1/2,0)$ &&&&&& $er$\cr
\+$n=188,\,d=1155=3\cdot 5\cdot 7\cdot 11,\,\eta=4,\,h=1$ &&&&&&&
$\langle3\rangle\oplus\langle-77\rangle\oplus\langle-5\rangle$ &&&&&& $nr$\cr
\+$n=189,\,d=1155=3\cdot 5\cdot 7\cdot 11,\,\eta=8,\,h=1$  &&&&&&&
$\langle1155\rangle\oplus\langle-1\rangle\oplus\langle-1\rangle$
&&&&&& $er$\cr
\+$n=190,\,d=1155=3\cdot 5\cdot 7\cdot 11,\,\eta=14,\,h=0$ &&&&&&&
$\langle22\rangle\oplus\langle-15\rangle\oplus\langle-14\rangle
(1/2,0,1/2)$ &&&&&& $er$\cr
\+$n=191,\,d=1254=2\cdot 3\cdot 11\cdot 19,\,\eta=4,\,h=1$ &&&&&&&
$\langle38\rangle\oplus\langle-22\rangle\oplus\langle-6\rangle
(1/2,1/2,0)$ &&&&&& $nr$\cr
\+$n=192,\,d=1365=3\cdot 5\cdot 7\cdot 13,\,\eta=1,\,h=1$ &&&&&&&
$\langle455\rangle\oplus\langle-3\rangle\oplus\langle-1\rangle$ &&&&&& $nr$\cr
\+$n=193,\,d=1365=3\cdot 5\cdot 7\cdot 13,\,\eta=5,\,h=1$ &&&&&&&
$\langle5\rangle\oplus\langle-39\rangle\oplus\langle-7\rangle$ &&&&&& $er$\cr
\+$n=194,\,d=1365=3\cdot 5\cdot 7\cdot 13,\,\eta=15,\,h=1$ &&&&&&&
$\langle6\rangle\oplus\langle-65\rangle\oplus\langle-14\rangle
(1/2,0,1/2)$ &&&&&& $er$\cr
\+$n=195,\,d=1394=2\cdot 17\cdot 41,\,\eta=3,\,h=1$ &&&&&&&
$U\oplus\langle-1394\rangle$ &&&&&& $nr$\cr
\+$n=196,\,d=1659=3\cdot 7\cdot 79,\,\eta=2,\,h=1$ &&&&&&&
$\langle79\rangle\oplus\langle-14\rangle\oplus\langle-6\rangle
(0,1/2,1/2)$ &&&&&& $nr$\cr
\+$n=197,\,d=1785=3\cdot 5\cdot 7\cdot 17,\,\eta=4,\,h=1$ &&&&&&&
$\langle30\rangle\oplus\langle-17\rangle\oplus\langle-14\rangle
(1/2,0,1/2)$ &&&&&& $nr$\cr
\+$n=198,\,d=1785=3\cdot 5\cdot 7\cdot 17,\,\eta=6,\,h=1$  &&&&&&&
$\langle17\rangle\oplus\langle-15\rangle\oplus\langle-7\rangle$ &&&&&& $nr$\cr
\+$n=199,\,d=2145=3\cdot 5\cdot 11\cdot 13,\,\eta=10,\,h=1$ &&&&&&&
$\langle1\rangle\oplus\langle-143\rangle\oplus\langle-15\rangle$
&&&&&& $nr$\cr
\+$n=200,\,d=2210=2\cdot 5\cdot 13\cdot 17,\,\eta=6,\,h=1$ &&&&&&&
$\langle2\rangle\oplus\langle-170\rangle\oplus\langle-26\rangle
(1/2,0,1/2)$ &&&&&& $nr$\cr
\+$n=201,\,d=2310=2\cdot 3\cdot 5\cdot 7\cdot 11,\,\eta=1,\,h=1$ &&&&&&&
$\langle6\rangle\oplus\langle-70\rangle\oplus\langle-22\rangle
(1/2,0,1/2)$ &&&&&& $nr$\cr
\+$n=202,\,d=2730=2\cdot 3\cdot 5\cdot 7\cdot 13,\,\eta=6,\,h=0$ &&&&&&&
$\langle2\rangle\oplus\langle-546\rangle\oplus\langle-10\rangle
(1/2,0,1/2)$ &&&&&& $nr$\cr
\+$n=203,\,d=3311=7\cdot 11\cdot 43,\,\eta=1,\,h=1$ &&&&&&&
$\langle43\rangle\oplus\langle-77
\rangle\oplus\langle-1\rangle$ &&&&&& $nr$\cr
\+$n=204,\,d=3570=2\cdot 3\cdot 5\cdot 7\cdot 17,\,\eta=7,\,h=1$ &&&&&&&
$\langle2\rangle\oplus\langle-714\rangle\oplus\langle-10\rangle
(1/2,0,1/2)$ &&&&&& $nr$\cr
\+$n=205,\,d=3990=2\cdot 3\cdot 5\cdot 7\cdot 19,\,\eta=4,\,h=1$ &&&&&&&
$\langle30\rangle\oplus\langle-38\rangle\oplus\langle-14\rangle
(1/2,1/2,0)$ &&&&&& $nr$\cr
\+$n=206,\,d=4466=2\cdot 7\cdot 11\cdot 29,\,\eta=1,\,h=1$ &&&&&&&
$U\oplus\langle-4466\rangle$ &&&&&& $nr$\cr}

\enddocument
\end

\enddocument
\end